\documentclass[11pt]{book}
\usepackage{a4}

\usepackage[english]{babel}
\usepackage{latexsym}
\usepackage{amsmath}
\usepackage{amscd}
\usepackage[small,nohug,heads=littlevee]{diagrams}
\diagramstyle[labelstyle=\scriptstyle]
\usepackage{dsfont}
\usepackage{bbm}
\usepackage{amssymb}
\usepackage{amsfonts}
\usepackage{eufrak}
\usepackage{euscript}
\usepackage[dvips]{graphicx}
\usepackage{theorem}

\newtheorem{defi}{Definition}[chapter]

\newtheorem{theo}[defi]{Theorem}
\newtheorem{lem}[defi]{Lemma}
\newtheorem{cor}[defi]{Corollary}
\newtheorem{prop}[defi]{Proposition}

\newtheorem{examp}[defi]{Example}
\newtheorem{examps}[defi]{Examples}
\def\R{\text{I\!\!\! R}}

\def\C{\mathbb{C}}

\def\N{\mathbb{N}}

\def\TT{\mathbb{T}}
\def\Z{\mathbb{Z}}
\def\1{\mathbbm{1}}

\def\Aa{\mathfrak{A}}
\def\Ba{\mathfrak{B}}
\def\Ca{\mathfrak{C}}
\def\Da{\mathfrak{D}}
\def\Fa{\mathfrak{F}}

\def\Ha{\mathfrak{H}}

\def\Ha{\mathfrak{H}}

\def\Ja{\mathfrak{J}}
\def\Ka{\mathfrak{K}}
\def\La{\mathfrak{L}}
\def\Lainf{\mathfrak{L}^\infty}
\def\Ma{\mathfrak{M}}

\def\Na{\mathfrak{N}}
\def\Pa{\mathfrak{P}}

\def\Sa{\mathfrak{S}}
\def\Ua{\mathfrak{U}}
\def\Za{\mathfrak{Z}}
\def\Fu{\mathfrak{Fu}}
\def\Alg{\mathfrak{Alg}}

\def\SA{\EuScript{A}}
\def\SB{\EuScript{B}}
\def\SF{\EuScript{F}}

\def\A{\mathcal{A}}
\def\D{\mathcal{D}}
\def\E{\mathcal{E}}
\def\L{\mathcal{L}}

\def\H{\mathcal{H}}

\def\S{\mathcal{S}}

\def\U{\mathcal{U}}
\def\X{\mathcal{X}}
\def\Y{\mathcal{Y}}

\def\Ebf{\mathbf{E}}
\def\Kbf{\mathbf{K}}
\def\Mbf{\mathbf{M}}
\def\Gr{\mathbf{G}}
\def\Ge{\tilde{\mathbf{G}}}
\def\Hr{\mathbf{H}}

\def\L{\mathcal{L}}
\def\LH{\mathcal{L}(\mathfrak{H})}
\def\LCH{\mathcal{LC}(\Ha)}

\def\rep{\mathbf{rep}}
\def\irrep{\mathbf{irrep}}

\def\o{\omega}

\def\Cs{\mathcal{C}}
\def\Cnn{\mathcal{C}_{\circ\circ}}
\def\Cg{\mathcal{C}^\infty}

\def\EA{\mathbf{E}_{\Aa}}

\def\autMa{\mathbf{Aut}(\Ma)}
\def\innMa{\mathbf{Inn}(\Ma)}
\def\spaMa{\mathbf{Spa}(\Ma)}
\def\outMa{\mathbf{Out}(\Ma)}
\def\Aut{\mathbf{Aut}}

\def\*{{}^{*}}
\def\rc{{}^\wr}

\def\im{{\text{im}}}
\def\const{{\text{const}}}
\def\supp{{\text{supp }}}
\def\Tr{\text{Tr\,}}
\def\singsupp{{\text{singsupp}}}

\def\rem{\hspace{-0.75cm} \textbf{Remark}: }

\def\proof{\hspace{-0.75cm} Proof: }
\def\proofof{\hspace{-0.75cm} Proof of }
\def\sproof{\hspace{-0.75cm} Sketch of the Proof: }
\def\qed{\begin{flushright}$\Box$\end{flushright}}

\usepackage{fancyhdr}
\begin{document}

\noindent

\begin{center}

\thispagestyle{empty}

\vspace*{3cm}
{\huge \textbf{On Galois Correspondence and\\ \vspace{6mm} Non-Commutative Martingales}}\\
\vspace{4cm} 
{\LARGE Dissertation}\\
\vspace{0.5cm}
{\Large zur Erlangung des Doktorgrades \\
der Fakult\"at f\"ur Wirtschafts- und Sozialwissenschaften\\ \vspace{0.1cm}
der Helmut-Schmidt-Universit\"at\\ \vspace{0.1cm}
Universit\"at der Bundeswehr Hamburg}\\

\vspace{2cm}

{\large vorgelegt von}\\
\vspace{1cm}
{\LARGE Timor Saffary}\\ 
\vspace{0.2cm}
{\large aus Kabul, Afghanistan}\\
\vspace{2.5cm}
{\Large  Hamburg}\\
\vspace{0.2cm}
{\Large November 2008}

\end{center}

\newpage

\thispagestyle{empty}
\vspace*{14cm}
\begin{flushleft}
\hspace*{-0.8cm}
\begin{tabular}{ll}
\textbf{Gutachter der Dissertation:} & Prof. Dr. W. Krumbholz\\
                                     & Prof. Dr. J. Michali\v cek\\
\\
\textbf{Fachpr\"ufer:}               & Prof. Dr. W. Krumbholz (Quantitative Methoden)\\
                                     & Prof. Dr. H. Hebbel (Statistik)\\ 
                                     & Prof. Dr. U. T\"ushaus (ABWL)\\ 
\\
\textbf{Datum der Pr\"ufung:}      & 10. Dezember 2008\\
\\
\textbf{Vorsitzender des Promotionsausschusses:} & Prof. Dr. W. Krumbholz\\
\\
\textbf{Dekan der Fakult\"at:} & Prof. Dr. G. Reiner\\
\end{tabular}
\end{flushleft}

\newpage

\thispagestyle{empty}
\vspace*{8cm}
\begin{center}
To my aunt\\ 
\vspace{0,5cm}
Hakimah Gizabi
\end{center}

\newpage

\thispagestyle{empty}

\begin{center}
{\large\textbf{On Galois Correspondence\\ and  Non-Commutative Martingales}}
\end{center}

\vspace{2cm}

\begin{center}
\textbf{Abstract:}\\
\end{center}
The subject of this thesis is Galois correspondence for von Neumann algebras and its interplay with non-commutative probability theory. After a brief introduction to representation theory for compact groups, in particular to Peter-Weyl theorem, and to operator algebras, including von Neumann algebras, automorphism groups, crossed products and decomposition theory, we formulate first steps of a non-commutative version of probability theory and introduce non-abelian analogues of stochastic processes and martingales. The central objects are a von Neumann algebra $\Ma$ and a compact group $\Gr$ acting on $\Ma$, for which we give in three consecutive steps, i.e. for inner, spatial and general automorphism groups one-to-one correspondences between subgroups of $\Gr$ and von Neumann subalgebras of $\Ma$. Furthermore, we identify non-abelian martingales in our approach and prove for them a convergence theorem.

\vspace{1.5cm}

\begin{center}
\textbf{Zusammenfassung:}
\end{center}
Gegenstand dieser Dissertation ist Galois-Korrespondenz f\"ur von Neumann Algebren und deren Zusammenspiel mit nichtkommutativer Wahrscheinlichkeitstheorie. Nach einer kurzen Einf\"uhrung in die Darstellungstheorie f\"ur kompakte Gruppen, insbesondere in das Theorem von Peter-Weyl, und in Operatoralgebren, einschlie\ss lich von Neumann Algebren, Gruppen von Automorphismen, semidirekte Produkte und Dekompositionstheorie, formulieren wir die ersten Schritte einer nichtkommutativen Version der Wahrscheinlichkeitstheorie und pr\"asentieren nichtabelsche Analoga f\"ur stochastische Prozesse und Martingale. Die zentralen Objekte dieser Arbeit sind eine von Neumann Algebra $\Ma$ und eine kompakte Gruppe $\Gr$, die auf $\Ma$ wirkt, f\"ur die wir in drei aufeinander aufbauenden F\"allen, d.h.~f\"ur innere, \"au\ss ere und allgemeine Gruppen von Automorphismen, bijektive Korrespondenzen zwischen Untergruppen von $\Gr$ und von Neumann Unteralgebren von $\Ma$ angeben. Dar\"uberhinaus identifizieren wir in unserem Formalismus nichtabelsche Martingale und beweisen f\"ur diese ein Konvergenztheorem.

\newpage

\setcounter{page}{1}
\setcounter{tocdepth}{2}
\tableofcontents

\setcounter{secnumdepth}{-1}
\setcounter{secnumdepth}{2}



\chapter{Introduction}

This thesis is concerned with Galois correspondence in operator algebras and non-commutative probability theory as well as with their interplay. To be more precise, we formulate a one-to-one correspondence between von Neumann subalgebras and subgroups of an automorphism group which allows the introduction of a non-commutative version of martingales in the framework given in \cite{Coja:2005}. Both subject are not only of utmost importance in mathematics but also offer increasingly far reaching applications in other disciplines.\\
The most impressive example is the interaction of operator algebras and theoretical physics, in particular quantum field theory, where von Neumann algebras are assigned the crucial r\^ole of algebras of local observables and one-parameter automorphism groups represent their dynamics \cite{Haag:1996}. Here, Galois correspondence makes an powerful and invaluable contribution to the analysis of their substructure which turns out to be decisive for understanding the underlying quantum system. \\
The intense interaction between mathematics and other scientific fields also applies to stochastic processes, in fact, their root may be traced back to botany, namely they grew out of efforts to give the Brownian motion a solid mathematical foundation. Since then, we have witnessed them penetrating more and more disciplines and transferring back the benefits of stochastic methods notably to physics, chemistry, biology and finance \cite{Freund:2000}, \cite{Gardiner:2004}, \cite{Musiela:2005}, \cite{Paul:1999}. In the second part of last century, the invention of non-commutative probability theory inspired by quantum theory was tackled, starting with the formulation of a non-commutative analogue of conditional expectation \cite{Umagaki:1954} and followed by non-abelian martingales in late seventies, see \cite{Randrianantoanina:2007} and references therein. A totally new access to this young theory is made possible through non-commutative differential geometry \cite{Dimakis:1993}. Quantisation of stochastic processes, e.g. Brownian motion, or stochastic differential equations may lead to deeper insights in elementary particle physics and disclose new non-classical effects, \cite{Parthasarathy:1992}, \cite{Accardi:1982}, \cite{Accardi:2002}. The first applications to finance, especially the investigation and reformulation of Black-Scholes model for option pricing, are proposed in \cite{Forgy:2002}, \cite{Chen:2001}.\\ \\
To start with, one of the greatest mathematical invention of 19th century was without doubt Galois theory, which connected for the first time group theory with field theory and therefore allowed the affiliation of problems from the more complex field theory to the much better understood group theory, in particular millennium old classical problems, e.g. constructibility of geometric objects with compass and straightedge, could be easily solved. It has become since then irreplaceable not only in mathematics but proved to be of utmost importance on other fields, too.\\
One of the disciplines for which Galois theory has become an increasingly important tool is operator algebra,  in particular in the classification of von Neumann algebras. For this purpose, one starts with a topological group $\Gr$ acting on a Neumann algebra $\Ma$, 
\begin{align*}
\alpha:\Gr&\longrightarrow\autMa\\
g&\mapsto\alpha^g\,,
\end{align*}
and asks for a possible one-to-one correspondence between subgroups $\Hr$ of $\Gr$ and von Neumann subalgebras $\Aa$ of $\Ma$ which consist of pointwise $\Hr$-fixed elements.
\begin{gather}
\Hr\longleftrightarrow\Aa\equiv\Ma^\Hr.
\end{gather}
Obviously, each automorphism subgroup $\alpha^\Hr$ implies a fixpoint von Neumann subalgebra, but the opposite direction is by far not trivial, since von Neumann algebras do have a much richer structure with plenty of their subalgebras not corresponding to any subgroup. There have been many attempts to introduce such a Galois correspondence for special cases, especially for factors which play a central part not only in operator theory but also in mathematical physics, to be more specific in the algebraic approach to quantum field theory \cite{Haag:1996}. The following list contains main results on this subject but it should not be regarded as exhaustive:
\begin{itemize}
\item[(i)] Nakamura and Takeda \cite{Nakamura:1960} and Suzuki \cite{Suzuki:1960} analyse finite Groups $\Gr$ with minimal action on factors of type $II_1$ $\Ma$, i.e. 
\begin{gather*}
\big(\Ma^\Gr\big)'\cap\Ma=\C\1.
\end{gather*}
\item[(ii)] Kishimoto considers a one-to-one correspondence between closed, normal subgroups of compact groups $\Gr$ and $\Gr$-invariant von Neumann subalgebras \cite{Kishimoto:1977}.
\item[(iii)] In \cite{Connes:1978} Connes and Takesaki give a Galois correspondence between closed subgroups of a locally compact, abelian group $\Gr$ and $\widehat{\Gr}$-invariant von Neumann subalgebras, where $\widehat{\Gr}$ is the dual group of $\Gr$. The same case is dealt with by Nakagami and Takesaki \cite{Nakagami:1979}.
\item[(iv)]  In \cite{Choda:1978} Choda investigates the Galois correspondence for discrete groups with outer action of $\Gr$ on von Neumann algebras via their crossed products.
\item[(vi)] Izumi, Longo and Popa deal with compact and discrete groups with minimal action and outer action, respectively, on factors with separable preduals \cite{Izumi:1998}. They also give a generalised version to compact Kac algebras with minimal action on subfactors.
\end{itemize}
Here, we are concerned with a more general setting, namely a compact group $\Gr$ acting on a von Neumann algebra $\Ma\equiv(\Ma,\Ha)$, where $\Ha$ is the underlying Hilbert space, without assuming any additional conditions neither on the action of the group nor on the nature of the von Neumann subalgebras.\\
First, we make sure that our investigation does not turn out to be void and prove at least some of the fixpoint subalgebras not to be trivial.\\
In our approach, we are dealing for technical reasons separately with inner, spatial and general automorphisms. We call a group of automorphisms 
\begin{gather}
\alpha(A)=UAU^*,\qquad A\in\Ma,
\end{gather}
inner or spatial, if the unitary operator $U\in\Ma$ or $U\in\LH$, respectively, whereas general automorphism groups are supposed not to have such a description. We start our investigation with the inner case which constitutes the main pillar of the analysis. There, we examine Galois correspondence of the following kind
\begin{gather*}
\Hr\longleftrightarrow\Aa\equiv\Ma^\Hr\quad\text{with}\quad\Aa\rc\rc=\Aa,
\end{gather*}
where $\Aa\rc:=\Aa'\cap\Ma$ denotes the relative commutant of $\Aa$. Since a spatial automorphism on $\Ma$ is inner on $\LH$ one may transfer the procedures and some results from the inner to the spatial case where relative commutants turn to usual commutants, i.e.
\begin{gather*}
\Hr\longleftrightarrow\Aa\equiv\Ma^\Hr\quad\text{with}\quad\Aa''=\Aa.
\end{gather*}
The general setting is then built on the spatial case by shifting the analysis onto the crossed product $\Ma\rtimes_\alpha\Gr$ where the existence of a spatial automorphism group 
\begin{gather*}
\pi_\alpha\circ\alpha:\Gr\longrightarrow\L\big(L^2(\Gr,\Ha)\big)
\end{gather*}
is always ensured. Thus, we may apply the results already obtained for the inner and spatial actions on $\Ma$  to the crossed product and project them back onto $\Ma$ itself.\\
One of our main tools for the analysis of all three cases is Peter-Weyl theorem which ensures one for compact groups $\Gr$ a complete orthonormal system in $L^2(\Gr)$ and, thus, provides a decomposition of the space $L^2(\Gr)$ into $\Gr$-invariant subspaces.\\ \\
Our formulation of non-commutative martingales is based on the ansatz of Coja-Oghlan and Michali\v cek \cite{Coja:2005}. A real valued random variable on classical probability space $(\Omega,\SA,P)$ from the functional analytic viewpoint is understood as 
\begin{align}\label{randomvar}
X:\Cnn(\R)&\longrightarrow L^\infty(\Omega,P),
\end{align}
namely a $*$-homomorphism from an abelian $C^*$-algebra to an abelian von Neumann algebra. Now, they introduce its non-commutative counterpart by allowing general, i.e. non-commutative, algebras on both sides of \ref{randomvar}. Random variables are then determined as mappings from a $C^*$-algebra $\Aa$ to a von Neumann algebra $\Ma$,
\begin{align*}
X:\Aa&\longrightarrow\Ma,\\
X^t:\Sa(\Ma)&\longrightarrow\Sa(\Aa),
\end{align*}
where $\Sa(\cdot)$ denotes the state space, and probability measures are represented by normal states. Thus, a non-commutative probability space is represented by the pair $(\Ma,\varphi)$ consisting of a general von Neumann algebra $\Ma$ acting on a Hilbert space $\Ha$ and a normal state $\varphi$ on $\Ma$. In this framework, conditional expectations are projections from $\Ma$ onto a von Neumann subalgebra $\Na$, i.e.
\begin{gather*}
\E:\Ma\longrightarrow\Na,
\end{gather*}
satisfying some obvious conditions, and filtrations are given as a increasing sequence of von Neumann subalgebras $(\Ma_t)_{t\in\TT}$ such that their union is dense in $\Ma$. Finally, we define a non-abelian martingale with respect to the filtration $(\Ma_t)_{t\in\TT}$ in $\Ma$ as an increasing sequence $(X_t)_{t\in\TT}$ of random variables on $(\Ma,\varphi)$ fulfilling for all $s,t\in\TT$ with $s<t$
\begin{gather*}
\E(X_t|\Ma_s)=X_s.
\end{gather*}
\\ \\
This thesis is structured as follows.\\
The second chapter will be concerned with the very basics of representation theory for compact groups. Our main aim is the formulation of Peter-Weyl theorem, the equivalent of Fourier series for compact groups, which constitutes one of the main pillars of our analysis. \\
In chapter three we give a compact abstract of operator algebras confining ourselves to the most important notions and facts of particular importance. We will start with an introduction into $C^*$-algebras, von Neumann algebras and their one-parameter automorphism groups comprising together the setting of modular theory, which is applied throughout this thesis. Crossed products are discussed which, as aforementioned,  form one of our main tools in the investigation of general compact automorphism groups. Last but not least, decomposition theory is given in a nutshell which is needed to prove the invariant spaces being not trivial.\\
Our approach to non-commutative probability theory is discussed in chapter four. The formalism of Coja-Oghlan and Michali\v cek \cite{Coja:2005} is recapitulated upon which we give non-commutative analogues of conditional expectations, stochastic independence, stochastic processes and martingales.\\
The bulk of our results is depicted in chapter five where we begin with the analysis of general invariant subspaces of von Neumann algebras $\Ma$ and their dual algebras $\Ma^*$. Next, we deal with Galois correspondence for the inner, spatial and finally general automorphism groups. The investigation of Izumi et al. \cite{Izumi:1998} is localised in the framework of our ansatz. At the end, we apply our results to non-commutative probability, namely we identify non-commutative martingales and prove a convergence theorem for them.\\ 
We conclude this thesis we a summary and outlook.



\chapter{Representations of Compact Groups}

\begin{flushright}
 \emph{Just go on...and faith will soon return.}\\
\vspace{0,5cm}
Jean Le Rond D'Alembert\\
To a friend hesitant with respect to infinitesimals.
\end{flushright}
\vspace{0,5cm}
The Fourier transform is indeed one of the most powerful concepts in mathematics. It states that any square integrable function from the one-sphere to the complex numbers 
\begin{gather*}
f:S^1\cong\R/\Z=\big\{e^{ix}|\;x\in[0,2\pi]\big\}\longrightarrow\C
\end{gather*}
can be approximated in a Fourier series
\begin{gather*}
f(x)=\sum_{k\in\Z}a_ke^{ikx},
\end{gather*}
where the Fourier coefficients are given as
\begin{gather*}
a_k:=\frac{1}{2\pi}\int_0^{2\pi}f(x)e^{-ikx}dx.
\end{gather*}
The generalisation of this statement to arbitrary compact groups is the so-called Peter-Weyl theorem, namely if $\Gr$ is a compact topological group, $\big\{(\alpha_\sigma,\Ha_\sigma)\big\}_{\sigma\in\Sigma}$ a complete family of inequivalent unitary, irreducible representations of $\Gr$, $\{e_i^\sigma\}$ an orthonormal basis of $\Ha_\sigma$, then 
\begin{gather*}
\sqrt{d_\sigma}D_{jk}^\sigma(g),
\end{gather*}
where $D_{ij}(g):=\langle\alpha^g(e_j),e_i\rangle$ are the so-called matrix elements and $d_\sigma$ the dimension of $\alpha_\sigma$, form a complete orthonormal system in $L^2(\Gr)$. Moreover, the matrix coefficients form a dense subalgebra of $\C(\Gr)$ and of $L^2(\Gr)$ with respect to the supremum norm and the $L^2$-norm, respectively. They may be seen as the analogues of $e^{ikx}$.\\
This theorem is insofar of prime importance for this investigation as it enables the decomposition of $L^2(\Gr)$ into subspaces invariant under the action of $\Gr$. \\
The following composition of basic facts and main results is excerpeted mainly from \cite{Barut-Raczka:1980}, \cite{Gaal:1973} and \cite{Hewitt:1963}.

\section{Representations of Compact Groups}

To start with, the following theorem allows one to deal with only unitary representations in the case of compact groups.

\begin{defi}
A representation $\alpha:\Gr\longrightarrow\LH$ is said to be the direct sum of subrepresentations  $\alpha_i:\Gr\longrightarrow\L(\Ha_i)$ on $\alpha$-invariant subspaces $\Ha_i$ if $\Ha=\sum_i^\oplus\Ha_i$.
\end{defi}

\begin{defi}
A representation $\alpha:\Gr\longrightarrow\LH$ is said to be irreducible, if it has no proper $\alpha$-invariant subspaces. $\alpha$ is called completely reducible if it is a direct sum of irreducible subrepresentations.
\end{defi}

\begin{theo}
Let $\alpha:\Gr\longrightarrow\LH$ be a representation of a compact group $\Gr$ on a Hilbert space $\Ha$, then there exists a new scalar product in $\Ha$ with equivalent norm to the original one and which defines a new unitary representation of the group.
\end{theo}

\proof
One introduces a new scalar product in terms of the original one $\langle\cdot,\cdot\rangle$ by
\begin{gather*}
\langle u,v\rangle_n:=\int_\Gr\big\langle\alpha^{g}(u),\alpha^{g}(v)\big\rangle dg, \quad u,v\in\Ha.
\end{gather*}
For the induced norms one obtains
\begin{align*}
\|u\|_n^2&=\langle u,v\rangle_n\int_\Gr\big\langle\alpha^{g}(u),\alpha^{g}(u)\big\rangle dg\\
&\leq\left(\sup_{g\in\Gr}\|\alpha^g\|\right)^2\int_\Gr\langle u,u\rangle dg=\left(\sup_{g\in\Gr}\|\alpha^g\|\right)^2\langle u,u\rangle\\
&=\left(\sup_{g\in\Gr}\|\alpha^g\|\right)^2\|u\|^2.
\end{align*}
On the other hand we may write
\begin{align*}
\|u\|^2&=\langle u,u\rangle=\big\langle\alpha^{g^{-1}}\circ\alpha^{g}(u),\alpha^{g^{-1}}\circ\alpha^{g}(u)\big\rangle\\
&\leq\left(\sup_{g\in\Gr}\|\alpha^g\|\right)^2\big\langle\alpha^{g}(u),\alpha^{g}(u)\big\rangle\\
&=\left(\sup_{g\in\Gr}\|\alpha^g\|\right)^2\int_\Gr\big\langle\alpha^{g}(u),\alpha^{g}(u)\big\rangle dg\\
&=\left(\sup_{g\in\Gr}\|\alpha^g\|\right)^2\langle u,u\rangle_n\\
&=\left(\sup_{g\in\Gr}\|\alpha^g\|\right)^2\|u\|_n^2.
\end{align*}
Therefore, the equivalence of the norms is established as well as that of their induced strong topologies, i.e. the map $g\mapsto\alpha^g$, $g\in\Gr$, has to be continuous with respect to the new topology.\\
Each operator $\alpha^{g_1}$ is isometric with respect to the new scalar product and its domain is the whole Hilbert space, since for all $u,v\in\Ha$ one has
\begin{align*}
\big\langle\alpha^{g_1}(u),\alpha^{g_1}(v)\big\rangle_n:&=\int_\Gr\big\langle\alpha^{g_1g_2}(u),\alpha^{g_1g_2}(v)\big\rangle dg_2\\
&=\int_\Gr\big\langle\alpha^{g}(u),\alpha^{g}(v)\big\rangle dg\\
&=\langle u,v\rangle_n,
\end{align*} 
where $g_1,g_2,g\in\Gr$. 

\qed
\begin{defi}
Two representations $\alpha_1:\Gr\longrightarrow\L(\Ha_1)$ and $\alpha_2:\Gr\longrightarrow\L(\Ha_2)$ on the Hilbert spaces $\Ha_1$ and $\Ha_2$, respectively, are said to be equivalent, denoted by $\alpha_1\simeq\alpha_2$, if there exists a unitary operator $U:\Ha_1\longrightarrow\Ha_2$ such that
\begin{gather*}
\alpha_2^g=U\alpha_2^g\,U^{-1},\quad\forall g\in\Gr.
\end{gather*}
The set of all equivalence classes is denoted by $\Sigma\equiv\Sigma(\Gr)$ and their elements by $\sigma$.
\end{defi}

\begin{lem}\label{completelyirreducible}
Any finite-dimensional unitary representation $\alpha:\Gr\longrightarrow\LH$ is completely reducible. 
\end{lem}

\proof
If $\Ha_i$ is a proper invariant subspace of $\Ha$, then so is its orthogonal complement $\Ha_i^\perp$ and one obtains the decomposition $\Ha=\Ha_1\oplus\Ha_1^\perp$. In our finite-dimensional case we may proceed with this recipe until we achieve a decomposition of $\Ha$ into irreducible $\alpha$-invariant subspaces.
\qed

\begin{prop}\label{Weyloperator}
For every unitary representation $\alpha:\Gr\longrightarrow\L(\Ha)$ of a compact group $\Gr$ and a fixed vector $u\in\Ha$ the so-called Weyl operator 
\begin{gather*}
K_uv:=\int_\Gr\langle v,\alpha^g(u)\rangle\alpha^g(u)dg, \quad\forall v\in\Ha,
\end{gather*}
has the following properties:
\begin{itemize}
\item[(i)] $K_u$ is bounded;
\item[(ii)] $K_u$ commutes with $\alpha^g$ for all elements $g$ of $\Gr$ and vectors of $\Ha$;
\item[(iii)] $K_u$ is selfadjoint;
\item[(iv)] $K_u$ is a Hilbert-Schmidt operator;
\item[(v)] $K_u$ has only $\alpha$-invariant eigenspaces $\Ha_i$;
\item[(vi)] The Hilbert space $\Ha$ is spanned by the eigenspaces $\Ha_i$ of $K_u$, i.e. $\Ha=\Ha_0+\sum_i^\oplus\Ha_i$, $\dim\Ha_i<\infty$, where $\Ha_0$ is the in general infinite-dimensional $0$-eigenspace.
\end{itemize}
\end{prop}

\proof
\begin{itemize}
\item[(i)] Boundedness may be seen as follows:
\begin{align*}
\|K_uv\|\leq\int_\Gr\big|\big\langle v,\alpha^g(u)\big\rangle\big|\|\alpha^g(u)\|dg\leq\int_\Gr\|v\|\|\alpha^g(u)\|\|\alpha^g(u)\|dg=\|v\|\|u\|^2.
\end{align*}
\item[(ii)] The commutation property is ensured for all $v\in\Ha$ and $g_1\in\Gr$ by
\begin{align*}
K_u\alpha^{g_1}(v)&=\int_\Gr\big\langle\alpha^{g_1}(v),\alpha^{g_2}(u)\big\rangle\alpha^{g_2}dg_2\\
&=\int_\Gr\big\langle\alpha^{g_1}(v),\alpha^{g_1g_2}(u)\big\rangle\alpha^{g_1g_2}dg_2\\
&=\int_\Gr\big\langle v,\alpha^{g_2}(u)\big\rangle\alpha^{g_1g_2}dg_2\\
&=\alpha^{g_1}(K_uv).
\end{align*}
\item[(iii)] $K_u^*=K_u$ follows from:
\begin{align*}
\langle K_uv,w\rangle&=\int_\Gr\big\langle v,\alpha^g(u)\big\rangle\big\langle\alpha^g(u),w\big\rangle dg\\
&=\int_\Gr\overline{\big\langle w,\alpha^g(u)\big\rangle}\big\langle v,\alpha^g(u)\big\rangle dg\\
&=\Big\langle v,\int_\Gr\big\langle w,\alpha^g(u)\big\rangle\alpha^g(u)\Big\rangle dg\\
&=\langle v,K_uw\rangle.
\end{align*} 
\item[(iv)] Let $\{e_i\}$ be a basis in $\Ha$, then one obtains:
\begin{align*}
\sum_i\|K_ue_i\|^2&=\sum_i\int_\Gr\int_\Gr\big\langle e_i,\alpha^{g_1}(u)\big\rangle\big\langle \alpha^{g_1}(u),\alpha^{g_2}(u)\big\rangle\big\langle \alpha^{g_2}(u),e_i\big\rangle dg_1dg_2\\
&=\int_\Gr\int_\Gr\sum_i\big\langle e_i,\alpha^{g_1}(u)\big\rangle\big\langle \alpha^{g_1}(u),\alpha^{g_2}(u)\big\rangle\big\langle \alpha^{g_2}(u),e_i\big\rangle dg_1dg_2\\
&=\int_\Gr\int_\Gr\sum_i\big\langle \alpha^{g_1}(u),\alpha^{g_2}(u)\big\rangle\big\langle \alpha^{g_2}(u),\alpha^{g_1}(u)\big\rangle dg_1dg_2<\infty,
\end{align*}
where Lebesgue's theorem has been applied to the second equation. Since we assumed no conditions on the basis, $K_u$ has to be a Hilbert-Schmidt operator.
\item[(v)] Since each $\Ha_i$ is finite-dimensional the statement follows directly from Lemma \ref{completelyirreducible}.
\item[(vi)] This is a direct consequence of Rellich-Hilbert-Schmidt spectral theorem, confer \cite[p.655]{Barut-Raczka:1980}.
\end{itemize}
\qed

\begin{theo}
All irreducible unitary representations of a compact group $\Gr$ are finite-dimensional.
\end{theo}

\proof
In Proposition \ref{Weyloperator} we have seen that $K_u\alpha^{g_1}(v)= \alpha^{g_1}(K_uv)$, and hence $K_u^*\alpha^{g_1}(v)= \alpha^{g_1}(K_u^*v)$ for all $g\in\Gr$ and $u,v\in\Ha$. This implies for the Weyl operator the structure $K_u=\beta(u)\1$ leading to
\begin{gather*}
\langle K_uv,v\rangle=\int_\Gr\big\langle v,\alpha^g(u)\big\rangle\big\langle\alpha^g(u),v\big\rangle dg=\beta(u)\langle v,v\rangle
\end{gather*}
and therefore to:
\begin{gather}\label{equationWeyl}
\int_\Gr\big|\big\langle\alpha^g(u),v\big\rangle\big|^2dg=\beta(u)\|v\|^2.
\end{gather}
By changing the variables in the latter equation we obtain 
\begin{align*}
\beta(v)\|u\|^2&=\int_\Gr\big|\big\langle\alpha^g(v),u\big\rangle\big|^2dg=\int_\Gr\big|\big\langle u,\alpha^g(v)\big\rangle\big|^2dg\\
&=\int_\Gr\big|\big\langle\alpha^{g^{-1}}(u),v\big\rangle\big|^2dg=\int_\Gr\big|\big\langle\alpha^g(u),v\big\rangle\big|^2dg\\
&=\beta(u)\|v\|^2\\
\Longrightarrow&\quad\beta(u)=c\|u\|,
\end{align*}
with $c=\const$ and where we have used the relation 
\begin{gather*}
\int_\Gr f\big(g^{-1}\big)dg=\int_\Gr f(g)dg.
\end{gather*}
The constant $c$ can be identified as positive by setting $u=v$ with $\|u\|=1$ in Equation \ref{equationWeyl}. Choose now an arbitrary orthonormal basis $\{e_i\}_{i\in\underline{n}}$ in $\Ha$ and set $u=e_k$, $k\in\underline{n}$, and $v=e_1$, then Equation \ref{equationWeyl} becomes
\begin{gather}
\int_\Gr\big|\big\langle\alpha^g(e_k),e_1\big\rangle\big|^2dg=\beta(e_k)\|e_1\|^2=c.
\end{gather}
We complete the proof by using Parseval's inequality:
\begin{align*}
nc&=\sum_{k=1}^n\int_\Gr\big|\big\langle\alpha^g(e_k),e_1\big\rangle\big|^2dg=\int_\Gr\sum_{k=1}^n\big|\big\langle\alpha^g(e_k),e_1\big\rangle\big|^2dg\leq\int_\Gr\|e_1\|^2dg=1,
\end{align*}
namely the dimension $n$ of the Hilbert space $\Ha$ has to be finite.

\qed

\begin{theo}
Every unitary representation $\alpha:\Gr\longrightarrow\LH$ of a compact group $\Gr$ on a Hilbert space $\Ha$ can be written as a direct sum of irreducible finite-dimensional unitary representations.
\end{theo}

\proof
First of all, the subspace $\Ha\ominus\Ha_0$ cannot be empty, because one has $\langle K_uv,v\rangle>0$ for all vectors $v\in\Ha$ not being orthogonal to $u$ and can therefore be expressed as a direct orthogonal sum of finite-dimensional invariant subspaces, i.e. $\Ha\ominus\Ha_0=\sum_i^\oplus\Ha_i$. Each $\Ha_i$ can be decomposed into a direct, orthogonal sum of irreducible and $\alpha$-invariant subspaces, as shown in the proof of Lemma \ref{completelyirreducible}. Concerning $\Ha_0$, we construct for it a Weyl operator which, due to Lemma \ref{completelyirreducible} and the properties (v) and (vi) in Proposition \ref{Weyloperator}, ensures the existence of a non-trivial, finite-dimensional, minimal and $\alpha$-invariant subspace of $\Ha_0$. Finally, the smallest subspace $\Ha_1$ of $\Ha$ which includes all minimal, orthogonal and invariant subspaces must be the whole Hilbert space itself, since unitarity of $\alpha$ and $\alpha$-invariance of $\Ha_1$ implies $\alpha$-invariance for its complement $\Ha_1^\perp$. But this enforces $\Ha_1^\perp=0$ and therefore $\Ha_1=\Ha$, otherwise $\Ha_1^\perp$ would contain a non-trivial, $\alpha$-invariant finite-dimensional subspace, contrary to its definition.

\qed

\section{Peter-Weyl Theorem}

\begin{defi}
Let $\alpha:\Gr\longrightarrow\LH$ be a $n$-dimensional representation of a group $\Gr$ on the Hilbert space $\Ha$ and $\{e_i\}_{i\in\underline{n}}$, $n<\infty$, an orthonormal basis of $\Ha$. We define the matrix element with respect to the operator $\alpha^g$, $g\in\Gr$, as
\begin{gather*}
D_{ij}(g):=\big\langle\alpha^g(e_j),e_i\big\rangle.
\end{gather*}
We will denote these elements of all finite-dimensional, continuous and unitary representations of $\Gr$ by $\Da(\Gr)$.
\end{defi}

\begin{lem}
The set $\Da(\Gr)$ of all coefficients constitutes an involutive algebra over $\C$, the so-called coefficient algebra of $\Gr$.
\end{lem} 

\proof
Let $(\alpha_1,\Ha_1)$ and $(\alpha_2,\Ha_2)$ be two continuous, unitary and irreducible representations of $\Gr$ and $\xi_1,\eta_1\in\Ha_1$, $\xi_2,\eta_2\in\Ha_2$ and $c_1,c_2\in\C$. Then, linearity is given by
\begin{gather*}
c_1\langle\alpha_1\xi_1,\eta_1\rangle+c_1\langle\alpha_2\xi_2,\eta_2\rangle=\langle\alpha\xi,\eta\rangle\in\Da(\Gr),
\end{gather*}
where we have set $\alpha:=\alpha_1\oplus\alpha_2$, $\xi:=c_1\xi_1+c_2\xi_2$ and $\eta:=\eta_1+\eta_2$.\\
For $n$-dimensional and $m$-dimensional Hilbert spaces $\Ha^n$ and $\Ha^m$, respectively, the space $\L(\Ha^n,\Ha^m)$ of of all linear mappings from $\Ha^n$ to $\Ha^m$ has dimension $nm$. Their orthonormal bases $(\xi_1,...,\xi_n)$ and $(\eta_1,...,\eta_m)$ determine an orthonormal basis for $\L(\Ha^n,\Ha^m)\cong\Ha^{nm}$ as $\{E_{ij}|\;E_{ij}\xi_k:=\delta_{ik}\eta_j,\,i\in\underline{n},\,j\in\underline{m}\}$. The tensor product of two continuous, unitary and irreducible representations $\alpha_1$ and $\alpha_2$ on $\Ha^n$ and $\Ha^m$, respectively, is defined as
\begin{gather*}
(\alpha_1\otimes\alpha_2)^g(A):=\alpha_1^gB\,\alpha_2^{g^{-1}},\qquad\forall B\in\L(\Ha^n,\Ha^m),\,\forall g\in\Gr,
\end{gather*}
which turns out to be continuous and unitary but not irreducible. If we denote by $(A_{ij})_1^g$, $i,j=1,...,n$, and $(A_{ij})_2^g$, $i,j=1,...,m$, the matrices corresponding to $\alpha_1$ and $\alpha_2$ and in terms of $(\xi_1,...,\xi_n)$ and $(\eta_1,...,\eta_m)$, respectively, then we obtain:
\begin{align*}
\big\langle(\alpha_1\otimes\alpha_2)^g(E_{ij}),E_{kl}\big\rangle&=\big((\alpha_1\otimes\alpha_2)^g(E_{ij})\big)_{kl}\\
&=\Big(\alpha_1^gE_{ij}\,\alpha_2^{g^{-1}}\Big)_{kl}\\
&=\sum_{nm}(A_{km})_1^g(E_{ij})_{mn}(A_{nl})_2^{g^{-1}}\\
&=(A_{ki})_1^g(A_{jl})_2^{g^{-1}}\\
&=\big\langle\alpha_1^g(\xi_k),\xi_i\big\rangle\big\langle\alpha_2^{g^{-1}}(\eta_j),\eta_l\big\rangle\\
&=\big\langle\alpha_1^g(\xi_k),\xi_i\big\rangle\big\langle\overline{\alpha_2}\,{}^g(\eta_l),\eta_j\big\rangle,
\end{align*}
where he have set $\overline{\alpha_2}\,{}^g:=\overline{\alpha_2^g}$ for which one has $\big\langle\overline{\alpha_2}\,{}^g(\eta_l),\eta_j\big\rangle:=\big\langle\overline{\alpha_2^{g^{-1}}}(\eta_j),\eta_l\big\rangle$. Thus, the coefficients of the tensor product $\alpha_1\otimes\alpha_2$ are identified with the product of the coefficients of $\alpha_1$ and $\alpha_2$, i.e. the set $\Da(\Gr)$ is closed with respect to multiplication. The introduction of the involution in $\Da(\Gr)$
\begin{gather*}
\big\langle\overline{\alpha}^g(\xi),\eta\big\rangle:=\overline{\big\langle\alpha^g(\xi),\eta\big\rangle},\qquad\forall\xi,\eta\in\Ha,\,\forall g\in\Gr,
\end{gather*}
completes the proof.

\qed

\begin{theo}\label{coefficients}
If $\alpha_1:\Gr\longrightarrow\L(\Ha_1)$ and $\alpha_2:\Gr\longrightarrow\L(\Ha_2)$ are two irreducible unitary representations of a compact group $\Gr$ on the Hilbert spaces $\Ha_1$ and $\Ha_2$, respectively, and $d_1$ is the dimension of $\alpha_1$, then their matrix elements satisfy the following relations:
\begin{gather*}
\int_\Gr D_{ij}^1(g)\overline{D_{kl}^2}(g)dg=\begin{cases}\frac{1}{d_1}\delta_{ik}\delta_{jl}, & \alpha_1\simeq\alpha_2 \\ 0, & \text{else}. \end{cases}
\end{gather*}
\end{theo}

\proof
Let $(e_{ij}):=\delta_i^m\delta_j^n$, $i,j,m,n=1,2,...,d_s$, then the operator
\begin{gather*}
E_{ij}:=\int_\Gr\alpha_1^ge_{ij}\alpha_2^{g^{-1}}dg
\end{gather*}
fulfills for all $g_1\in\Gr$
\begin{gather}\label{Ecommutation}
\alpha_1^{g_1}E_{ij}=E_{ij}\alpha_2^{g_1}.
\end{gather}
For non-equivalent representations $\alpha_1$ and $\alpha_2$ one obtains due to Schur's lemma 
\begin{gather}\label{Eorthogonal}
\int_\Gr D_{li}^1(g)D_{jk}^2\big(g^{-1}\big)dg=\int_\Gr D_{li}^1(g)\overline{D_{jk}^2(g)}dg=0\\
\Longleftrightarrow\quad E_{ij}=0.\notag
\end{gather}
The operator $E_{ij}$ satisfies obviously for $\alpha_1=\alpha_2$ the commutaion relation \ref{Ecommutation} which in turn implies 
\begin{gather}\label{Elambda}
E_{ij}=\lambda_{ij}\1.
\end{gather}
Therefore, the orthogonality relation \ref{Eorthogonal} are ensured for pairs $(l,i)\neq(k,j)$. For pairs $(l,i)=(k,j)$ one obtains by definition and \ref{Elambda} 
\begin{gather*}
\big(E_{ii}\big)_{ll}=\int_\Gr D_{li}^1(g)D_{jk}^2\big(g^{-1}\big)dg=\int_\Gr\big|D_{li}^1(g)\big|^2dg=\lambda_{ii}.
\end{gather*}
The constant is determined via the definition of $E_{ij}$ as
\begin{gather*}
\Tr\big(E_{ii}\big)=d_1\lambda_{ii}=\int_\Gr\Tr\Big[\alpha_1^ge_{ii}\,\alpha_2^{g^{-1}}\Big]dg=\Tr e_{ii}=1,
\end{gather*}
i.e. $\lambda_{ii}=\frac{1}{d_1}$.

\qed

\begin{defi}
The trace of a finite-dimensional representation $\alpha$ of $\Gr$ is called character of $\alpha$ and is denoted by $\chi(g):=\Tr\,\alpha^g$, $g\in\Gr$.
\end{defi}

\begin{lem}\label{character}
Let $\Gr$ be a compact group, $\chi(g)$, $\chi_1(g)$ and $\chi_2(g)$ be the characters with respect to the representations $\alpha^g$, $\alpha_1^g$ and $\alpha_2^g$, $g\in\Gr$, of the group on the Hilbert space $\Ha$, respectively. Then the following statements hold:
\begin{itemize}
\item[(i)] $\chi(g)=D_{ii}(g)$;
\item[(ii)] $\chi\big(g_2^{-1}g_1g_2\big)=\chi(g_1)$;
\item[(iii)] $\chi\big(g^{-1}\big)=\overline{\chi}(g)\,$;
\item[(iv)] $\chi_1=\chi_2$ whenever $\alpha_1\simeq\alpha_2$;
\item[(v)] \begin{gather*}
\int_\Gr\chi_1(g)\overline{\chi_2}(g)dx=\begin{cases}0, &\alpha_1\neq\alpha_2\\ 1, &\alpha_1=\alpha_2\end{cases}.
\end{gather*}
\end{itemize}
\end{lem}

\proof
Let $\{e_i\}_{i\in\underline{n}}$ be an orthonormal basis of $\Ha$.
\begin{itemize}
\item[(i)] $\chi(g)=\Tr\,\alpha^g=\big\langle\alpha^g(e_i),e_i\big\rangle=D_{ii}(g)$;
\item[(ii)] $\chi\big(g_2^{-1}g_1g_2\big)=\Tr\,\alpha^{g_2^{-1}g_1g_2}=\Tr\big(\alpha^{g_2^{-1}}\alpha^{g_1}\alpha^{g_2}\big)=\chi(g_1)$;
\item[(iii)] $\chi\big(g^{-1}\big)=\Tr\,\alpha^{g^{-1}}=\Tr\big(\alpha^{g^{-1}}\big)^*=\overline{D_{ii}}(g)=\overline{\chi}(g)\,$;
\item[(iv)] By assumption there exists a unitary operator $U$ such that $\alpha_2=U\alpha_1U^{-1}$. Thus we may conclude $\chi_2=\Tr\,\alpha_2=\Tr(U\alpha_1U^{-1})=\Tr\,\alpha_1=\chi_1$;
\item[(v)] The application of Theorem \ref{coefficients} verifies this claim:
\begin{gather*}
\int_\Gr\chi_1(g)\overline{\chi_2}(g)dx=\int_\Gr D_{ii}^1(g)\overline{D_{jj}^2}(g)dx=\begin{cases}0, &\alpha_1\neq\alpha_2\\ \frac{1}{d}, &\alpha_1=\alpha_2\end{cases}.
\end{gather*}
\end{itemize}
\qed
We have now the ability to formulate the main statement of this chapter.
\begin{theo}[Peter-Weyl]\label{PeterWeyl}
Let $\Gr$ be a compact group, $\alpha_\sigma:\Gr\longrightarrow\LH$ an irreducible unitary representation belonging to the equivalence class $\sigma\in\Sigma\equiv\Sigma(\Gr)$, $d_\sigma$ the dimension of $\alpha_\sigma$ and $D_{jk}^\sigma(g)$, $g\in\Gr$, the matrix elements of $\alpha_\sigma$, then the functions 
\begin{gather*}
\sqrt{d_\sigma}D_{jk}^\sigma(g),\qquad j,k=1,...,d_\omega,
\end{gather*}
constitute a complete orthonormal system in $L^2(\Gr)$.
\end{theo}

\proof
Throughout the proof the letters $g_1,g_2,g_3$ and $g$ will be elements of $\Gr$. Since each linear closed subspace $L$ of $L^2(\Gr)$ spanned by the functions $\sqrt{d_\sigma}D_{jk}^\sigma(g)$ are invariant under right translations $\alpha_r$, the orthogonal complement $L^\perp$ has this property, too. Let $v$ be a nontrivial element of $\L^\perp$, then 
\begin{gather*}
u(g_1):=\int_\Gr\big[\alpha_r^{g_1}v(g_2)\big]\overline{v}(g_2)dg_2
\end{gather*}
is a continuous function on $\Gr$ belonging to $L^\perp$ and $u(e)=\|v\|^2>0$. Consider the operator 
\begin{gather*}
A\psi(g_1):=\int_\Gr w(g_1g_2^{-1})\psi(g_2)dg_2,\\
w(g):=u(g)+\overline{u(g^{-1})},
\end{gather*}
which is self-adjoint, compact and $\alpha_r$-invariant:
\begin{align*}
\big[\alpha_r^{g_3}(A\psi)\big](g_1)&=\int_\Gr w\big(g_1g_3g_2^{-1}\big)\psi(g_2)dg_2\\
&=\int_\Gr w\big(g_1g_2^{-1}\big)\psi(g_3g_2)dg_2\\
&=\big[A\alpha_r^{g_3}(\psi)\big](g_1).
\end{align*}
 Since $n\neq 0$ and therefore $w\neq 0$, $A$ has to have at least one nontrivial eigenvalue $\lambda$ with the corresponding eigenspace $\Ha(\lambda$). An arbitrary eigenfunction $\psi_\lambda(g)$ of $A$ has to be an element of $L^\perp$:
\begin{align*}
\int_\Gr\psi_\lambda(g)\overline{D_{ij}^\sigma}(g)dg&=\frac{1}{\lambda}\int_\Gr(A\psi_\lambda)(g)\overline{D_{ij}^\sigma}(g)dg\\
&=\frac{1}{\lambda}\sum_k\int_\Gr w(g_1)\overline{D_{ij}^\sigma}(g_1)dg_1\int_\Gr\psi_\lambda(g)\overline{D_{ij}^\sigma}(g)dg\\
&=0.
\end{align*}
$A$ inherits its $\alpha_r$-invariance to its eigenspace $\Ha_\lambda$ and one defines a completely reducible representation by 
\begin{gather*}
\alpha^{g_1}\big(\psi_\lambda(g_2)\big):=\psi_\lambda(g_2g_1).
\end{gather*}
Let $\Ha_\sigma(\lambda)\subset\Ha(\lambda)$ be an irreducible subspace equipped with an orthonormal basis $\{e_k^\sigma\}_{k\in\underline{d_\sigma}}$ then one obtains $e_k^\sigma\in L$, since one has
\begin{gather*}
\alpha_r^{g_1}e_k^\sigma(g_2)=e_k^\sigma(g_2g_1)=D_{jk}^\sigma(g_1)e_j^\sigma(g_2)
\end{gather*}
and thanks to the continuity of the eigenfunctions
\begin{gather*}
e_k^\sigma(g_1)=D_{jk}^\sigma(g_1)e_j^\sigma(e).
\end{gather*}
Consequently, $\Ha(\lambda)=\{0\}$ contrary to the assumption and therefore $L^\perp$ must be trivial. 
\qed

\begin{rem}
Two important consequences of Theorem \ref{PeterWeyl} are that the matrix coefficients form a dense subalgebra of $\Cg(\Gr)$ and of $L^2(\Gr)$ with respect to the supremum norm and the $L^2$-norm, respectively. 
\end{rem} 



\def\const{{\text{const}}}
\def\supp{{\text{supp}}}
\def\singsupp{{\text{singsupp}}}

\chapter{Automorphism Groups on von Neumann Algebras}

\begin{flushright}
 \emph{Alles sollte so einfach wie m\"oglich gemacht sein,\\ aber nicht einfacher.}\\
\vspace{0,5cm}
Albert Einstein
\end{flushright}
\vspace{0,5cm}
In this chapter we want to give in a nutshell the basic terminology and some concepts of operator algebras which are used throughout this thesis. We are concerned with $C^*$-algebras and von Neumann algebras, investigate their automorphism groups, in particular modular automorphism groups, and conclude with non-commutative integration and crossed products.\\
In the following we present a number of definitions and statements collected mainly from \cite{Takesaki1:1970}, \cite{Takesaki2:2001}, \cite{Bratteli:1979tw} and \cite{Stratila:1981}.

\section{$C^*$-Algebras, von Neumann Algebras}

An Algebra $\Aa$ is a complex vectorspace equipped with a bilinear, associative product
\begin{equation*}
\Aa\times\Aa\longrightarrow\Aa.
\end{equation*}
The algebra is said to be commutative or abelian, if the product possesses these properties. A map
\begin{equation*}
\begin{split}
\*:\Aa&\longrightarrow\Aa\\
A&\mapsto A^{*}
\end{split}
\end{equation*}
is called involution of $\Aa$, if for all $A,B\in\Aa$ and $a,b\in\C$ the following conditions are fulfilled:
\begin{itemize}
\item $(A^{*})^{*}=A$;
\item $(aA+bB)^{*}=\overline{a}A^{*}+\overline{b}B^{*}$;
\item $(AB)^{*}=B^{*}A^{*}$.
\end{itemize}
$A^{*}$ is called the adjoint of $A$ and an algebra with such an involution involutive algebra or $\*$-algebra. The set of all self-adjoint elements $A$ of $\Aa$, i.e. $A^*=A$, will be denoted by $\Aa_{sa}$. Since the involution is norm continuous $\Aa_{sa}$ is closed and therefore a Banach space. The unital algebra contains the neutral element $\1$. If the $\*$-algebra is also a Banach space and is satisfying 
\begin{equation}
\|A^{*}\|=\|A\|,\qquad\forall A\in\Aa,
\end{equation}
then it is said to be a Banach $\*$-algebra or $B\*$-algebra.

\begin{defi}
A $\*$-algebra $\Ua$ is called $U\*$-algebra, if its unitisation $\Ua^\1$ is the linear span of its unitary group $(\Ua^\1)_\U$.
\end{defi}

The $U\*$-algebra can be equipped with a norm $\|\cdot\|_\Ua:\Ua\longrightarrow\R_+$ by the so-called unitary norm
\begin{gather}
\|\A\|_\Ua:=\inf\left\{\sum_{i=1}^n|\lambda_i|\;\big|\quad A=\sum_{i=1}^n|\lambda_iU_i,\quad n\in\N,\lambda_i\in\C,U_i\in(\Ua^\1)_\Ua\right\}.
\end{gather}

\begin{defi}
A C$\*$-algebra $\Aa$ is a Banach$\*$-algebra, i.e. $\Aa$ is a Banach algebra with involution and $\|A^*\|=\|A\|$, such that 
\begin{gather}\label{cnorm}
\|A^{*}A\|=\|A\|^{2},\quad \forall A\in\Aa.
\end{gather}
$\Aa$ is said to be simple, if the only closed ideals are $\{0\}$ and $\Aa$ itself.
\end{defi}
It can be shown that an abelian and unital C$\*$-algebra $\Aa$ is (isometrically) isomorphic to $\mathcal{C}(\X)$ where $\X$ is a compact Hausdorff space. In the case of a non-unital C$\*$-algebra, $\Aa$ is (isometrically) isomorphic to $\mathcal{C}_0(\X)$, if $\X$ is locally compact. In both cases the space $\X$ is uniquely determined up to homomorphisms. The Hausdorff space $\X$ can be chosen as the set of characters of $\Aa$.

\begin{defi}
A character of an abelian C$\*$-algebra $\Aa$ is a nonzero linear map 
\begin{align*}
\o:\Aa&\longrightarrow\C 
\end{align*}
with the property
\begin{align*}
\o(AB)=\o(A)\o(B),\quad\forall A,B\in\Aa.
\end{align*}
\end{defi}




\begin{defi}
A C$\*$-algebra $\Ma$ is called $W^*$-algebra, if it is the dual space of some Banach space.
\end{defi}
If one is concerned with operator algebras, then it is usual to call $W^*$-algebras von Neumann algebras. Sakai has proved that the abstract characterisation of von Neumann algebras given above  is equivalent to the more usual definition via the commutator. Let $\Ha$ be a Hilbert space and $\LH$ the algebra of linear bounded operators
\begin{equation*}
A:\Ha\longrightarrow\Ha
\end{equation*}
equipped with the norm
\begin{equation*}
\|A\|:=\underset{\psi\in\Ha,\|\psi\|\leq1}{\sup}\|A\psi\|<\infty.
\end{equation*}
Define the commutant $\Aa'$ of $\Aa$ as the algebra
\begin{equation*}
\Aa':=\{X\in\LH:\;XA=AX,\,A\in\Aa\}
\end{equation*}
which has the following properties: 
\begin{equation}\label{commutant}
\begin{split}
\Aa\subseteq\Aa''=\Aa^{(4)}=\Aa^{(6)}=\cdots,\\
\Aa'=\Aa'''=\Aa^{(5)}=\Aa^{(7)}=\cdots,\\
\text{where}\quad\Aa^{(n+1)}:=\big(\Aa^{(n)}\big)'.
\end{split}
\end{equation}
Then, a von Neumann algebra $\Ma$ can be characterised via its double commutant by the requirement $\Ma=\Ma''$.

\begin{theo}
Every $\*$-algebra $\Ma\subseteq\LH$ is a von Neumann algebra if and only if $\Ma=\Ma''$ holds.
\end{theo}

\proof
Confer \cite[Theorem 3.5.]{Takesaki:1970}.
\qed

\begin{examps} Let $\Ha$ be a Hilbert space.
\begin{itemize} 
\item[(i)] $\LH$ is not only a von Neumann algebra, but even a factor since  $\LH'=\C\1$, confer Definition \ref{factor}.
\item[(ii)] The $C^*$-algebra of compact operators $\LCH$ on $\Ha$ cannot be a von Neumann algebra, because $\LCH''=(\C\1)'=\LH\neq\LCH.$
\end{itemize}
\end{examps}
In this thesis we will make extensive use of the following notion.
\begin{defi}
The relative commutant of a subalgebra $\Ba$ of a $C^*$-algebra $\Aa$ is denoted by $\Ba\rc$, i.e. $\Ba\rc:=\Ba'\cap\Aa$.
\end{defi}
For further discussions on von Neumann algebras, additional topologies apart from the uniform topology are needed. There are many in general inequivalent representations of a $C^*$-algebra by bounded operators on a Hilbert space which form always a closed algebra with respect to the uniform topology. But for a detailed analysis one relies on approximations of these operators and is therefore interested in weaker topologies. For the sake of completeness, we introduce the most important ones although we will not make use of all of them explicitly.

\begin{defi}
On $\LH$ we distinguish between the locally convex operator topologies defined by the following sets of seminorms:
\begin{itemize}
\item [(i)] $\sigma$-weak: $p_{(\xi_n),(\eta_n)}(A):=\big|\underset{n}{\sum}\langle \xi_{n},A\eta_{n}\rangle\big|$ for all $\xi_n,\eta_n\in\Ha$ \\
with $\underset{n}{\sum}\big(\|\xi_{n}\|^{2}+\|\eta_{n}\|^{2}\big)<\infty$;
\item [(ii)] weak: $p_{\xi,\eta}(A):=|\langle \xi,A\eta\rangle|$ for all $\xi,\eta\in\Ha$;
\item [(iii)] strong: $p_{\xi}(A):=\|A\xi\|$ for all $\xi\in\Ha$;
\item [(iv)] $\sigma$-strong: $p_{(\xi_n)}(A):=\underset{n}{\sum}\|A\xi_{n}\|^{2}$ for all $\xi_n\in\Ha$ with $\underset{n}{\sum}\|\xi_{n}\|^{2}<\infty$;
\item [(v)] $\*$-strong: $A\mapsto p_{\xi}(A)+p_{\xi}(A^{*})$ for all $\xi\in\Ha$;
\item [(vi)] $\sigma\*$-strong: $A\mapsto\left\{\underset{n}{\sum}\|A\xi_{n}\|^{2}+\underset{n}{\sum}\|A^{*}\xi_{n}\|^{2}\right\}^{1/2}$ for all $\xi_n\in\Ha$; \\
with $\underset{n}{\sum}\|\xi_{n}\|^{2}<\infty$.
\end{itemize}
\end{defi} 
For our purposes the $\sigma$-weak and $\sigma$-strong topologies will be of particular interest, because the modular group of automorphisms is continuous with respect to them. If we denote by ``$<$'' the relation ``finer than'', then the following diagram shows the relation between the various topologies:\\

\begin{center}
uniform $<\;\sigma\*$-strong$\;<\;\sigma$-strong$\;<\;\sigma$-weak\\
${}\hspace{2,3cm}\wedge\hspace{1.7cm}\wedge\hspace{1.5cm}\wedge$\\
${}\hspace{2cm}\*$-strong\hspace{0.2cm}$<\hspace{0.4cm}$strong$\hspace{0.3cm}<\hspace{0.5cm}$weak\quad.\\ 
\end{center}

\vspace{0.5cm}
\hspace{-0.7cm}``Finer'' means the existence of more open or closed sets, less compact sets, more linear functionals and less convergent sequences. The $\sigma\*$-strong, $\sigma$-strong and $\sigma$-weak topologies allow for the same continuous linear functionals. This statement remains true if one drops the '$\sigma$-'. The main difference between the $\*$-strong and the $\sigma\*$-strong topology is the fact that the involution $A\mapsto A^{*}$ is only continuous with respect to the former topology. \\
On several occasions we will follow the notation of common literature on topolgies in locally convex spaces. For such a space $\X$ we denote the weak topology by $\sigma(\X,\X^*)$, i.e. the topology generated by the semi-norms 
\begin{gather*}
p_{X^*}(X):=|X^*(X)|,\quad X^*\in\X^*.
\end{gather*}
Analogously, the semi-norms
\begin{gather*}
p_X(X^*):=|X^*(X)|,\quad X\in\X,
\end{gather*}
generate the so-called weak${}^*$-topology denoted by $\sigma(\X^*,\X)$.

\begin{prop}
The unitary operators $\U(\Ha)$ on a Hilbert space $\Ha$ form a $\*$-strongly closed group.
\end{prop}

\proof
Choose an arbitrary sequence $(U_i)_{i\in\N}$ converging $\*$-strongly to $U\in\LH$, then one obtains for all $\xi\in\Ha$:
\begin{gather*}
\|U\xi\|=\lim_{i\rightarrow\infty}\|U_i\xi\|=\|\xi\|,
\|U^*\xi\|=\lim_{i\rightarrow\infty}\|U_i^*\xi\|=\|\xi\|.
\end{gather*}
Thus both $U$ and $U^*$ are isometries, i.e. $U$ is a unitary operator.
\qed
The strong$\*$- , the strong and the topologies coincide on $\U(\Ha)$, but in the case of infinite dimensional Hilbert spaces $\U(\Ha)$ does not to be strongly closed and, in general, the strong limit of unitary operators is an isometry, not necessarily unitary.\\
The following theorem shows the equivalence of these topologies in the case of von Neumann algebras.

\begin{theo}[Bicommutant -]\label{bicommutanttheo}
Let $\Aa\subset\LH$ be a nondegenerate $\*$-algebra, i.e. $[\Aa\Ha]:=\{\overline{\text{span}A\xi}|\;A\in\Aa,\xi\in\Ha\}=\Ha$, then the following statements are equivalent:

\begin{itemize}
\item[(i)] $\Aa=\Aa''$.
\item[(ii)]$\Aa$ is weakly closed.
\item[(iii)]$\Aa$ is strongly closed.
\item[(iv)]$\Aa$ is $\*$-strongly closed.
\item[(v)]$\Aa$ is $\sigma$-weakly closed.
\item[(vi)]$\Aa$ is $\sigma$-strongly closed.
\item[(vii)]$\Aa$ is $\sigma\*$-strongly closed.
\end{itemize}
\end{theo} 

\proof
Since the commutant of a self-adjoint set is closed with respect to all locally convex topologies, $(i)$ implies the rest of the statements. $(v)\Rightarrow(vi)\Rightarrow(vii)$ are trivial. $(vii)\Rightarrow(v)$ follows from the fact that every $\sigma\*$-strongly continuous functional is also $\sigma$-weakly continuous and  $(ii)\Rightarrow(iii)\Rightarrow(iv)\Rightarrow(vii)$ are evident. Consider now the countable infinite sum $\tilde{\Ha}:=\bigoplus_{n=1}^\infty\Ha_n$ with $\Ha_n:=\Ha$ for all $n\in\N$, then the map
\begin{align*}
\pi:\LH&\longrightarrow\L(\tilde{\Ha})\\
A&\mapsto\pi(A),\quad\pi(A)\left(\bigoplus_n\xi_n\right):=)\left(\bigoplus_n A\xi_n\right)
\end{align*}
defines a $\*$-automorphism of $\LH$ into a subalgebra of $\L(\tilde{\Ha})$, thus  $(vi)\Rightarrow(i)$ is valid.
\qed
Therefore, a von Neumann algebra is a weakly closed C$\*$-subalgebra of $\LH$. The bicommutant theorem also states that the closure of the $\*$-algebra is independent of the choice of a particular topology. An immediate consequence of which is the next conclusion. 

\begin{cor}[von Neumann density -]
The nondegenerate $\*$-algebra of operators $\Aa$ on $\Ha$ is weakly, strongly, $\*$-strongly, $\sigma$-weakly, $\sigma$-strongly and $\sigma\*$-strongly dense in $\Aa''$.
\end{cor}

\proof
Let $\overline{\Aa}$ be the closure of $\Aa$ with respect to one of the topologies, then one has $\Aa'=\overline{\Aa}'$ and therefore $\Aa''=\overline{\Aa}''=\overline{\Aa}$, where the second equation is justified by Theorem \ref{bicommutanttheo}.
\qed
This statement can be tightened and a stronger version is formulated by Kaplansky.

\begin{theo}[Kaplanskys density -]
Let $\Aa$ be a $\*$-algebra of operators on $\Ha$, $\Ma$ its weak closure and $\Aa_1$ and $\Ma_1$ their unit balls, respectively. Then $\Aa_1$ is $\sigma\*$-strongly dense in $\Ma_1$. 
\end{theo}

\proof
Since $\Aa_1$ is norm dense in the unit ball of the normal closure, we may assume without restriction that $\Aa$ is a $C^*$-algebra. Moreover, Theorem \ref{bicommutanttheo} $\Aa$ is also $\sigma\*$-strongly dense in $\Ma$ and therefore $\Aa_{sa}$, the self-adjoint operators of $\Aa$, are  $\sigma$-strongly dense in $\Ma_{sa}$. Since the function 
\begin{align*}
f:[-1,1]&\longrightarrow[-1,1]\\
x&\mapsto f(x):=2x(1+x^2)^{-1}
\end{align*}
is strictly increasing, the map
\begin{align*}
f:\LH_{sa}&\longrightarrow\LH_{sa}\\
A&\mapsto f(A):=2A(\1+A^2)^{-1}
\end{align*}
maps for all $C^*$-algebras $\Ba\subseteq\LH$ $\Ba_{sa}$ into $\Ba_{1,sa}$. For all $A,B\in\LH_{sa}$ one has
\begin{align*}
\frac{1}{2}[f(A)-f(B)]&=(\1+A^2)^{-1}\big[A(\1+B^2)-(\1+A^2)B\big](\1+B^2)^{-1}\\
&=(\1+A^2)^{-1}(A-B)(\1+B^2)^{-1}\\
&\quad+(\1+A^2)^{-1}A(B-A)B(\1+B^2)^{-1}\\
&=(\1+A^2)^{-1}(A-B)(\1+B^2)^{-1}+\frac{1}{4}f(A)(B-A)f(B),
\end{align*}
i.e. $f$ is strongly continuous since $\|(\1+A^2)^{-1}\|\leq1$ and $\|(\1+A^2)^{-1}A\|\leq1$ and $f(A)$ converges to strongly to $f(B)$ whenever $A$ tends strongly to $B$. $f$ must be also $\sigma$-strongly continuous because these two topologies coincide on the unit ball. Thus $\Aa_{1,sa}=f(\Aa_{sa})$ is $\sigma$-strongly dense in $\Ma_{1,sa}=f(\Ma_{sa})$.
\qed

\begin{defi}
Let $\Aa,\Ba$ be C$\*$-algebras then a linear map $\pi:\Aa\longrightarrow\Ba$ is called a $\*$-homomorphism if

\begin{itemize}
\item[(i)] $\pi(AB)=\pi(A)\pi(B)$ and
\item[(ii)]$\pi(A^{*})=\big(\pi(A)\big)^{*}$
\end{itemize}
hold for all $A,B\in\Aa$.
\end{defi}
The notions mono-, epi-, iso-, endo- and automorphism are introduced in the usual manner. We want to keep hold of some fundamental and frequently used properties of $\*$-homomorphisms in the below lemma.

\begin{lem}\label{homomorphism}
Let $\Aa$ and $\Ba$ be C$\*$-algebras and $\pi:\Aa\longrightarrow\Ba$ a $\*$-homomorphism. Then the following statements are valid:
\begin{itemize}
\item[(i)] $\pi$ preserves positivity: $A\in\Aa_+:=\{A\in\Aa_{sa}|\;\sigma(A)\geq 0\Rightarrow\pi(A)\geq 0$.
\item[(ii)]$\pi$ is continuous and $\|\phi(A)\|\leq\|A\|$, i.e. $\pi$ can only decrease the norm.
\item[(iii)] $\pi$ is an $\*$-isomorphism if and only if $\ker\pi:=\{A\in\Aa:\,\pi(A)=0\}=\{0\}$ is fulfilled.
\item[(iv)] An $\*$-isomorphism is automatically isometrical, i.e. norm preserving: $\|\pi(A)\|=\|A\|$.
\item[(v)] If $\pi$ is an $\*$-isomorphism, then the image $\pi(\Aa)$ of a C$\*$-algebra $\Aa$ is again a C$\*$-algebra.
\end{itemize}
\end{lem} 

\proof
\begin{itemize}
\item[(i)] An arbitrary positive element $A\in\Aa$ can be written as $A=B^*B$ for some $B\in\Aa$ and one concludes
\begin{gather*}
\pi(A)=\pi(B^*B)=\pi(B)^*\pi(B)\geq0.
\end{gather*}
\item[(ii)] Let us first assume self-adjointness for $A$, then one has
\begin{gather*}
\|\pi(A)\|=\sup\big\{|\lambda||\;\lambda\in\sigma[\pi(A)]\big\}\|
\end{gather*}
Consider the projection $P:=\pi(\1_\Aa)$, i.e. $P^2=P$, where $\1_\Aa$ is the identity in $\Aa$, If one replaces $\Ba$ by the new $C^*$-algebra $P\Ba P$, then $P$ becomes the identity of the new algebra $\Ba$ with $\pi(\Aa)\subseteq\Ba$. Consequently we obtains the following estimations:
\begin{gather*}
\|\pi(A)\|^2\leq\sup\big\{|\lambda||\;\lambda\in\sigma_\Aa(A)\big\}\leq\|A\|^2.
\end{gather*}
The general case where $A$ is not self-adjoint and thus both statements follow from the $C^*$-norm property \ref{cnorm}:
\begin{gather*}
\|\pi(A)\|^2\leq\|\pi(A^*A)\|\leq\|A^*A\|\leq\|A\|^2, \qquad\forall A\in\Aa.
\end{gather*}

\item[(iii)] The $\*$-homomorphism is a $\*$-isomorphism if and only if its range is equal to $\Ba$, i.e. $\ker\pi=\{0\}$.
\item[(iv)] As $\pi$ is an isomorphism its kernel has to be trivial and $\pi^{-1}$ exists. The statement is then a result of $(ii)$:
\begin{gather}\label{proof(iv)} 
\|A\|=\|\pi^{-1}\big(\pi(A)\big)\|\leq\|\pi(A)\|\leq\|A\|.
\end{gather} 
\item[(v)] Obviously, $\pi(\Aa)$ is a $\*$-subalgebra of $\Ba$. Because $\Aa$  is as a $C^*$-algebra closed and, due to $(iv)$, $\pi$ is an isometry, $\pi(\Aa)$ is closed, contains the identity and inherits the $C^*$-norm property from $\Aa$ and thus forms a $C^*$-algebra.
\end{itemize}
\qed
The last statement remains valid if $\pi$ is only a $\*$-homomorphism. In the case of von Neumann algebras one can state stricter properties of $\*$-homomorphism, namely each $\*$-homomorphism $\pi:\Ma\longrightarrow\Na$ between von Neumann algebras $\Ma$ and $\Na$ is $\sigma$-weakly as well as $\sigma$-strongly continuous.

\begin{defi}
A representation of a C$\*$-algebra $\Aa$ is ein pair $(\Ha,\pi)$, consisting of a complex Hilbert space $\Ha$ and a $\*$-homomorphism $\pi:\Aa\longrightarrow\LH$, and is said to be faithful if $\pi$ is a $\*$-isomorphism, and nondegenerate if 
\begin{gather*}
\{v\in\Ha:\;\pi(A)v=0,\,A\in\Aa\}=\{0\}.
\end{gather*}
A subspace $\Fa\subset\Ha$ is called invariant under $\pi(\Aa)$ if $\pi(A)\Fa\subseteq\Fa$ holds for all $A\in\Aa$. Whenever $\{0\}$ and $\Ha$ are the sole invariant and closed subspaces we call the representation $(\Ha,\pi)$  irreducible, otherwise reducible. Two representations $(\Ha_{1},\pi_{1})$ and $(\Ha_{2},\pi_{2})$ are said to be unitarily equivalent if there exists an unitary operator $U:\Ha_{1}\longrightarrow\Ha_{2}$ such that $\pi_{2}(A)=U\pi_{1}U^{*}$. If the two Hilbert spaces are connected via an $\*$-isomorphism the we speak about (quasi-)equivalence.
\end{defi}

\begin{cor}
The representation $(\Ha,\pi)$ of a C$\*$-algebra $\Aa$ is faithful if and only if one of the following equivalent conditions are satisfied:
\begin{itemize}
\item[(i)] $\ker\pi=\{0\}$.
\item[(ii)] $\|\pi(A)\|=\|A\|,\quad\forall A\in\Aa$.
\item[(iii)] $\pi(A)>0,\quad\forall A\in\Aa_{+}$.
\end{itemize}
\end{cor}

\proof
Faithfulness is obviously equivalent to $(i)$. The implication $(i)\Rightarrow(ii)$ is valid since $\ker\pi=\{0\}$ ensures the existence of the inverse of $\pi$ and one applies the Estimation \ref{proof(iv)}.\\
If we assume $(ii)$ then for strictly positive elements $A$ we get $\|\pi(A)\|=\|A\|>0$ and thus $\pi(A)>0$ or $\pi(A)<0$. But the first claim of Lemma \ref{homomorphism} states $\pi(A)\geq0$ and therefore $(iii)$ follows.\\
Let us assume $\ker\pi\neq\{0\}$, namely there exists an element $A\in\ker\pi$ with $A\neq0$ and $\pi(A)=0$, thus $\pi(A^*A)=0$. On the other side one has $\|A^*A\|\geq0$ and $\|A^*A\|=\|A\|^2\neq0$ and therefore $A^*A\geq0$. Consequently, $(i)$ is falsified. 

\qed
Each representation $(\Ha,\pi)$ of a C$\*$-algebra $\Aa$ defines a faithful representation for the quotient algebra $\Aa_{\pi}:=\Aa/\ker\pi$. The representation of a simple algebra is alway faithful.

\begin{defi}
Let $\Ma$ be a von Neumann algebra on a Hilbert space $\Ha$, then a subspace $\Ha'\subseteq\Ha$ is said to be separating for $\Ma$ if and only if $A\xi=0$ implies $A=0$ for all $A\in\Ma$ and $\xi\in\Ha'$. A cyclic  representation is a triple $(\Ha,\pi,\xi)$, where $(\Ha,\pi)$ is a representation of $\Aa$ and $\xi\in\Ha$ is a cyclic vector for $\pi$ in $\Ha$, i.e. $\{\pi(A)\xi:\;A\in\Aa\}$ is dense in $\Ha$.
\end{defi}

\begin{cor}\label{cyclicseperating}
For a von Neumann algebra $\Ma$ on a Hilbert space $\Ha$ and $\Ka\subseteq\Ha$ cyclicity of $\Ka$ for $\Ma$ is equivalent to $\Ka$ being separating for $\Ma'$.
\end{cor}

\proof
Let $\Ka$ be a cyclic set for $\Ma$, $\xi\in\Ka$ and $A'\in\Ma'$ with $A'\xi=\{0\}$, then one has for all elements of $\Ma$ 
\begin{gather*}
A'B\xi=BA'\xi=0.
\end{gather*}
This implies $A'[\Ma\Ka]=0$ and therefore $A'=0$. If, on the other hand, $\Ka$ is separating for $\Ma'$, then $P':=[\Ma\mathcal{K}]$ is a projection in $\Ma'$ satisfying 
\begin{gather*}
(\1-P')\Ka=\{0\}.
\end{gather*}
Therefore $1-P'=0$ and $[\Ma\Ka]=\Ha$, i.e. $\Ka$ is cyclic for $\Ma$.

\qed
Thus, if the vector $\xi\in\Ha$ is cyclic and separating for the von Neumann algebra $\Ma$, then, since $\Ma=\Ma''$, it has these properties also for its commutant $\Ma'$.\\
Every nondegenerate representation of a $C^*$-algebra can be described as a direct sum of cyclic sub-representations. Therefore, the discussion of such representations can be restricted to the investigation of the cyclic ones.

\begin{defi}\label{definitionsforstates}
Let $\Aa$ be a unital C$\*$-algebra and $\o:\Aa\longrightarrow\C$ a linear functional on $\Aa$, then $\o$ is said to be

\begin{itemize}
\item[(i)] hermitian, if $\o(A^{*})=\overline{\o(A)}$ for all $A\in\Aa$;
\item[(ii)] positive, if $\o(A)\geq 0$ for all $A\in\Aa,\;A\geq 0$;
\item[(iii)] a weight, if $\o$ is positive;
\item[(iv)] a state, if $\o$ is positive and normalized, i.e. $\o(\1)=1$;
\item[(v)] a faithful state if $\o$ is a state and $\o(A)>0$ for all $A\in\Aa_+$, the set of positive elements of $\Aa$;
\item[(vi)] a trace, if $\o(AB)=\o(BA)$ for all $A,B\in\Aa$;
\item[(vii)] a vector state, if $\o(A)\equiv\o_\xi(A):=\big(\xi,\pi(A)\xi\big)$ for a non-degenerate representation $(\Ha,\pi)$ of $\Aa$ and some vector $\xi\in\Ha$ with $\|\xi\|=1$.
\end{itemize} 
\end{defi}
In the case of an abelian $C^*$-algebra $\Aa$, $\o$ is a pure state if and only if $\o(AB)=\o(A)\o(B)$ holds for all $A,B\in\Aa$. If $\Aa$ does not possess a unit $\1$, then the norm property $(iv)$ is replaced by the condition 
\begin{gather*}
\|\o\|:=\sup\big\{|\o(A)|\,|\;A\in\Aa\text{ and }\|A\|=1\big\}=1. 
\end{gather*}

\begin{defi}\label{definitionnormal}
If $\Ma$ is a von Neumann algebra, $\o$ a positive linear functional on $\Ma$ and $\o\big(\emph{l.u.b.}_n(A_{\alpha})\big)=\emph{l.u.b.}_n\big(\omega(A_{\alpha})\big)$ holds for all increasing bounded nets $(A_n)$ in $\Ma_{+}$, where `$\emph{l.u.b.}$' denotes the least upper bound, then $\o$ is said to be normal.
\end{defi} 

\begin{defi}\label{definitionsemifinite}
A trace $\o$ on a von Neumann algebra $\Ma$ is said to be semifinite, if the set
$$
\Ma_+^\o:=\big\{A\in\Ma_+|\;\o(A)<\infty\big\}
$$
is dense in $\Ma$. A von Neumann algebra $\Ma$ is called semifinite, if there exists a faithful, normal and semifinite weight on $\Ma$.
 
\end{defi}
It can be shown that the commutant $\Ma'$ of a semifinite von Neumann algebra $\Ma$ on a separable Hilbert space is also semifinite. \\
For the classification of von Neumann algebras it is useful to introduce some notations. Let us consider an involutive Banach algebra $\Aa$, then $P\in\Aa$ is called a projection if $P^{2}=P$ and $P^{*}=P$. Two projections $P,Q\in\Ma$, where $\Ma$ is a von Neumann algebra, are said to be equivalent, abbreviated by $P\sim Q$, if there exists a $W\in\Ma$, such that $W^{*}W=P$ and $WW^{*}=Q$. The projection $P$ is said to be finite if $Q\leq P$ and $Q\sim P$ imply $Q=P$, otherwise it is called infinite. If there is no nonzero finite projection $Q\leq P$, $Q\in\Ma$, then $P$ is called purely infinite. If $ZP\neq 0$ is infinite for every central projection $Z\in\Ma$, i.e. for every projection on $\Ma\cap\Ma'$, then $P$ is called properly infinite. $P$ is said to be $abelian$, if $P\Ma P$ is a commutative subalgebra of $\Ma$. Two projections $P_1$ and $P_2$ are said to be centrally orthogonal, if the smallest central projections $Z_{P_1}$ and $Z_{P_2}$ majorizing $P_1$ and $P_2$, respectively, are orthogonal. \\
Every projection $P\in\Ma$ can be uniquely described as the sum of two centrally orthogonal projections $P_1,P_2\in\Ma$, such that $P_1$ is finite and $P_2$ is properly infinite. Since the set spanned by the projections is dense in the von Neumann algebra $\Ma$, the property of the projections can be used to characterise their algebras.

\begin{defi}
A von Neumann algebra $\Ma$ is called finite, infinite, properly infinite or purely infinite if the identity projection $\1$ possesses these properties. The von Neumann algebra is said to be of\\
\\
\begin{tabular}{l l}

Type $I$,  & if every nonzero central projection in $\Ma$ majorises a nonzero\\
         & abelian projection in $\Ma$;\\
Type $II$, & if $\Ma$ has no nonzero abelian projection and every nonzero central\\
         & projection in $\Ma$ majorises a nonzero finite projection in $\Ma$;\\
Type $II_{1}$, & if $\Ma$ is of type II and finite;\\
Type $II_{\infty}$, & if $\Ma$ is of type II and has no nonzero central finite projection;\\
Type $III$, & if $\Ma$ is purely infinite.\\

\end{tabular}
\\
\end{defi}
If $\Ma$ is of type $I,II$ or $III$, then so is its commutant $\Ma'$ and contrariwise. For von Neumann algebras of type $II_{1}$ and $II_{\infty}$ this statement is in general not valid. The aforementioned characterisation allows one to decompose all von Neumann algebras completely in terms of these different types, namely every von Neumann algebra $\Ma$ is uniquely decomposable into the direct sum
$$
\Ma=\Ma_I\oplus\Ma_{II_1}\oplus\Ma_{II_\infty}\oplus\Ma_{III}.
$$
A von Neumann algebra is semifinite if and only if it contains no type $III$ component.

\begin{defi}\label{factor}
A von Neumann algebra $\Ma$ is called a factor, if it possesses a trivial center, i.e. $\Za(\Ma)=\Ma\cap\Ma'=\C\1$. A state $\o$ on a C$\*$-algebra $\Aa$ is said to be a factor state or primary state, if $\pi_\o(\Aa)''$ is a factor, where $\pi_\o$ is the corresponding cyclic representation.
\end{defi}
A factor is either of type $I,II_{1},II_{\infty}$ or $III$.

\begin{defi}
A state $\o$ is said to be pure, if the only positive linear functionals majorised by $\o$ are of the form $\lambda\o$ with $0\leq\lambda\leq 1$.
\end{defi}
Pure states are extremal points of the set of states $\EA$ on $\Aa$, which means that a pure state $\o$ is not describable as a convex linear combination 
$$
\o=\lambda\o_1+(1-\lambda)\o_2,\quad 0<\lambda<1,
$$
of different states $\o_1$ and $\o_2$.\\
To each arbitrary nondegenerate representation of a $C^*$-algebra and a vector $\xi\in\Ha$ with $\|\xi\|=1$ one can assign a state, the so-called vector state. The construction in the opposite direction is ensured by the next

\begin{theo}
For an arbitrary state $\o$ on a C$\*$-algebra $\Aa$ there exists a (up to unitary equivalence) unique cyclic representation $(\Ha_\o,\pi_\o,\xi_\o)$ of $\Aa$, the so-called canonical cyclic representation of $\Aa$ with respect to $\o$, such that 
\begin{equation*}
\o(A)=\big(\xi_{\o},\pi_{\o}(A)\xi_{\o}\big)
\end{equation*}
holds for all $A\in\Aa$. 
\end{theo}

\proof
First, we want to construct the Hilbert space $H_\o$. We define equivalence classes 
\begin{gather*}
\Psi_A:=\big\{\widehat{A}|\;\widehat{A}=A+I,\;I\in\mathfrak{I}\big\}
\end{gather*}
where 
\begin{gather*}
\mathfrak{I}:=\big\{A\in\Aa|\;\o(A*A)=0\big\}
\end{gather*}
is a left ideal of $\Aa$. It is easy to show, that these equivalence classes form a complex vector space and we choose its completion provided with the scalar product $\langle A,B\rangle:=\o(A^*B)$ as our representation space $\Ha_\o$. The representatives of elements $A\in\Aa$ and the vector $\xi_\o$ are defined as
\begin{gather*}
\pi_\o(A)\Psi_B:=\Psi_{AB}\quad\text{and}\quad \xi_\o:=\Psi_\1.
\end{gather*}
In order to prove the uniqueness of this triple let us assume the existence of another representation $(\Ha_\o',\pi_\o',\xi_\o')$ with
\begin{equation*}
\o(A)=\big(\xi_{\o}',\pi_{\o}'(A)\xi_{\o}'\big)
\end{equation*}
for all $A\in\Aa$.  Then, the algebras $\pi_\o(\Aa)$ and $\pi_o'(\Aa)$ should be unitary equivalent, i.e. there should exist a unitary operator $U$ such that
\begin{equation*}
U\pi_\o(A)U^{-1}=\pi_\o'(A)
\end{equation*}
for all $A\in\Aa$.
\begin{equation*}
U\xi_\o=\xi_\o'.
\end{equation*}
These requirements are satisfied if we set
\begin{gather*}
U\pi_\o(A)\xi_\o=\pi_\o'(A)\xi_\o',
\end{gather*}
leading directly to:
\begin{align*}
\big(U\pi_\o(A)\xi_\o,U\pi_\o(B)\xi_\o\big)&=\big(\pi_\o'(A)\xi_\o',\pi_\o'(B)\xi_\o'\big)\\
&=\o(A^*B)\\
&=\big(\pi_\o(A)\xi_\o,\pi_\o(B)\xi_\o\big)
\end{align*}

\qed

\begin{defi}\label{ergodic}
Let $\Gr$ be a locally compact group and $\Aa$ a $C^*$-algebra, then a continuous homomorphism 
\begin{align*}
\alpha:\Gr&\longrightarrow\Aut(\Aa)\\
g&\mapsto(\alpha^g:\Aa\longrightarrow\Aa)
\end{align*}
is called an action of $\Gr$ into $\Aut(\Aa)$ and the triple $(\Aa,\Gr,\alpha)$ a covariant system over $\Gr$. The covariant system is said to be ergodic if $\Ma^\alpha=\C\1$ and centrally ergodic if $\Za_\Ma=\C\1$.
\end{defi}
From now on we focus on one-parameter groups of $\*$-automorphisms of $C\*$-algebras or von Neumann algebras and use the algebraic setting for this purpose. The fundamental algebraic tool for the investigation of infinitesimal generators is the symmetric derivation. The defining characteristics of derivations are naturally motivated by the main algebraic properties of the groups:
\begin{gather*}
\alpha^t(A)^*=\alpha^t(A^*)\quad\text{and}\\
\alpha^t(AB)=\alpha^t(A)\alpha^t(B).
\end{gather*}  

\begin{defi}
A symmetric derivation $\delta$ of a $C\*$-algebra $\Aa$ is a linear operator from a $\*$-subalgebra $\mathbf{D}(\delta)$, the domain of $\delta$, into $\Aa$ satisfying for all $A,B\in\mathbf{D}(\delta)$ the conditions:
\begin{itemize}
\item[(ii)] $\delta(AB)=\delta(A)B+A\delta(B)$;
\item[(ii)] $\delta(A)^*=\delta(A^*)$.
\end{itemize}
\end{defi}
The terminology 'symmetric' is referred to the second property.
\begin{lem}
The set of derivations $\Da(\Aa)$ on a $C\*$-algebra $\Aa$ form with the usual addition and scalar multiplication of operators and the product
$$
[\delta_1,\delta_2]:=\delta_1\delta_2-\delta_2\delta_1
$$ 
a Lie algebra.
\end{lem}

\proof
The recuired properties for a Lie algebra can be shown straightforward by definition.
\qed
The set of derivations is also called derivation algebra.

\begin{prop}\label{bounded_derivation}
Each anywhere defined derivation $\delta$ from a $C\*$-algebra into a larger $C\*$-algebra is automatically bounded.
\end{prop}

\proof
In the case where the $C\*$-algebra does not contain the identity, we first unitilise it. Let us consider a positive element $A\in\mathbf{D}(\delta)$, then its square root lies also in $\mathbf{D}(\delta)$ and one has:
\begin{gather}\label{prop-bounded}
\delta(A^*A)\geq\delta(A^*)A+A^*\delta(A).
\end{gather}
Then, $\delta$ is norm-closeable which can be seen as follows. Let us assume $A_n\in\mathbf{D}(\delta)$ with $A_n\rightarrow 0$ and $\delta A_n\rightarrow A$, then we have to show $A=0$. We are allowed to prove equivalently this implication for its square $A_n^*A_n$. But this is a direct consequence of the estimation \eqref{prop-bounded}, and therefore our statement is proved, since an everywhere defined norm-closed operator is always bounded by the graph norm.

\qed
For the discussion of automorphism groups, the notion of spatial or inner derivations is of great importance, due to the fact that they occur as their infinitesimal generators and, therefore, are irreplaceable in their analysis.

\begin{defi}
A symmetric derivation $\delta$ of a $C\*$-algebra $\Aa$ of bounded operators on a Hilbert space $\Ha$ is called inner if there exists a symmetric operator $H\in\Aa$ with the properties
\begin{itemize}
\item[(i)] $\delta(A)=i[H,A]$ for all $A\in\mathbf{D}(\delta)$, and
\item[(ii)] $\mathbf{D}(\delta)\mathbf{D}(H)\subseteq\mathbf{D}(H)$.
\end{itemize}
If $H\in\LH$ then we call $\delta$ spatial. $H$ is said to implement $\delta$.
\end{defi}
Our investigation in chapter five will be structured corresponding to this defintion.
\begin{lem}
The set of inner derivations $\Da_i(\Aa)$ on a $C\*$-algebra $\Aa$ form an ideal in the derivation Lie algebra $\Da(\Aa)$.
\end{lem}

\proof
Obviously, $\Da_i(\Aa)$ is a linear space. Let $\delta\in\Da(\Aa)$ and $\delta_i\in\Da_i(\Aa)$, i.e. there exists a selfadjoint element $H$ with $\delta_i(A)=i[H,A]$, then one has:
\begin{align*}
[\delta_i,\delta](A)&=i\big[H,\delta(A)\big]-i\delta[H,A]\\
&=-i\delta(H)A+iA\delta(H)\\
&=-i\big[\delta(H),A\big].
\end{align*}
The proof is completed by showing the selfadjointness via the symmetrical property of the derivation:
\begin{gather*}
\big(\delta(H)\big)^*=\delta(H^*)=\delta(H).
\end{gather*}
\qed

\begin{theo}[Borchers-Averson]\label{Borchers-Averson}
For a ${}\,\sigma$-weakly continuous one-parame-ter group of ${}^*$-automorphisms $\big(\alpha^t\big)_{t\in\R}$ of a von Neumann algebra $\Ma$ the following conditions are equivalent: 
\begin{itemize}
\item[(i)] There exists a strongly continuous one-parameter group of unitaries $U_t$, $t\in\R$, in $\LH$ with non-negative spectrum such that
\begin{gather*}
\alpha^t(A)=U_tAU_t^*
\end{gather*}
for all $A\in\Ma$ and $t\in\R$.
\item[(ii)]  There exists a strongly continuous one-parameter group of unitaries $U_t$, $t\in\R$, in $\Ma$ with non-negative spectrum such that
\begin{gather*}
\alpha^t(A)=U_tAU_t^*
\end{gather*}
for all $A\in\Ma$ and $t\in\R$.
\item[(iii)] One has for the spectral subspaces:
\begin{gather*}
\bigcap_{t\in\R}\big[\Ma^\alpha[t,\infty\rangle\Ha\big]=\{0\}.
\end{gather*}
\end{itemize}
If these statements are given, then $U$ may be written as
\begin{gather*}
U_t=\int_\Re^{-itp}dP(p),
\end{gather*}
where $P$ is the unique projection valued measure on $\R$ such that
\begin{gather*}P[t,\infty\rangle\Ha=\bigcap_{s<t}\big[\Ma^\alpha[s,\infty\rangle\Ha\big].
\end{gather*}
\end{theo}

\proof
See  \cite[Theorem 3.2.46]{Bratteli:1979tw}.
\qed

\begin{cor}\label{Neumann-derivation}
If $\delta$ is an anywhere defined derivation, i.e. bounded derivation, of a von Neumann algebra $\Ma$, then there exists an element $H=H^*\in\Ma$, $\|H\|\leq\frac{\delta}{2}$, with
\begin{gather*}
\delta(A)=i[H,A]
\end{gather*}
for all $A\in\Ma$. In case of a $C\*$-algebra $\Aa$, for all representation $\pi$ of $\Aa$ there is always an $H=H^*\in\pi(\Aa)$, $\|H\|\leq\frac{\delta}{2}$, such that
\begin{gather*}
\pi\big(\delta(A)\big)=i\big[H,\pi(A)\big].
\end{gather*}

\end{cor}

\proof
In the case of von Neumann algebras $\Ma$, our derivation $\delta$ is due to Proposition \ref{bounded_derivation} bounded and allows on $\Ma$ the introduction of the norm-continuous one-parameter ${}^*$-automorphism group
\begin{gather*}
\alpha^t(A):=e^{t\delta}A.
\end{gather*}
This means for the spectrum that
\begin{gather*}
\sigma(\alpha)=i\sigma(\delta),
\end{gather*}
and 
\begin{gather*}
\Ma^\alpha[t,\infty\rangle=\{0\},\qquad\text{for}\quad t>\|\delta\|.
\end{gather*}
Finally, Theorem \ref{Borchers-Averson} ensures one the existence of a strongly continuous group of unitaries $(U_t)_{t\in\R}$ with positive selfadjoint infinitesimal generator $H_0\in\Ma$ which implements $\delta$. The last statement can be shown by setting $H:=H_0-\frac{\|\delta\|}{2}$.\\
In the case of a $C\*$-algebra $\Aa$, let us consider the kernel of the representation $\mathfrak{I}:=\ker\pi$. If $A\in\mathfrak{I}_+$, i.e. $A=B^2$ for some $B\in\mathfrak{I}_+$, then one obtains $\delta(A)=\delta(BB)=\delta(B)B+B\delta(B)\in\mathfrak{I}$ because $\mathfrak{I}$ is an ideal. Thus $\delta(\mathfrak{I})\subseteq\mathfrak{I}$ and we may introduce a derivation of $\pi(\Aa)$ as
\begin{align*}
\tilde{\delta}:\pi(\Aa)&\longrightarrow\pi(\Aa)\\
\pi(A)&\mapsto\tilde{\delta}\big(\pi(A)\big):=\pi\big(\delta(A)\big)
\end{align*}
This derivation has a unique $\sigma$-weakly closed extension onto the $\sigma$-weakly closure of $\pi(\Aa)$ and we may apply the aforementioned result for von Neumann algebras. 
\qed

\begin{lem}
Let $\Aa$ be a $C^*$-algebra and $\delta:\Aa\longrightarrow\Aa$ a derivation, i.e. it is everywhere defined, then $\delta\equiv 0$.
\end{lem}

\proof
This statement is an immediate consequence of Corollary \ref{Neumann-derivation}.
\qed


\section{One-Parameter Automorphism Groups}

In this section we give a brief insight into the general theory of one-parameter automorphism groups. These objects are used in many theories for the description of its dynamical nature, e.g. the time evolution in physics and stochastic processes. For more details we refer the reader to \cite{Engel:2000} and \cite{Jacob:2001}.
\begin{cor}
For a linear operator $\delta$ acting on a $C\*$-algebra $\Aa$ the following statements are equivalent:
\begin{itemize}
\item[(i)] $\delta$ is an everywhere defined symmetric derivation of $\Aa$.
\item[(ii)] $\delta$ generates a norm-continuous one-parameter group of ${}^*$-automorphisms $\big(\alpha^t\big)_{t\in\R}$ of $\Aa$. 
\end{itemize}
If these conditions are valid, then for any representation $\pi$ of $\Aa$ there is always an selfadjoint operator $H\in\pi(\Aa)$with
\begin{gather*}
\pi\big(\alpha^t(A)\big)=e^{itH}\pi(A)e^{-itH}, \qquad t\in\R,
\end{gather*}
for all elements $A\in\Aa$.
\end{cor}

\proof
The implication $(i)\Rightarrow(ii)$ is clear since, due to Proposition \ref{bounded_derivation}, $\delta$ is a bounded operator and therefore generates a norm-continuous ${}^*$-automorphism group.\\
For the direction $(ii)\Rightarrow(i)$ let us assume that $\alpha^t$ is norm-continuous at the origin. Since for sufficiently small $t\in\R$ one has
\begin{gather*}
\Big\|\frac{1}{t}\int_0^t \alpha^sds-\1\Big\|<1,
\end{gather*}
the integral 
\begin{gather*}
X_t=\frac{1}{t}\int_0^t \alpha^sds
\end{gather*}
has a bounded inverse. Moreover we may write
\begin{align*}
\frac{1}{h}(\alpha^h-\1)X_t&=\frac{1}{th}\int_0^t(\alpha^{s+h}-\alpha^s)ds\\
&=\frac{1}{th}\int_t^{t+h}\alpha^sds-\frac{1}{th}\int_0^h\alpha^sds\\
&=\frac{1}{t}(\alpha^t-\1)X_h\\
&\overset{h\rightarrow0}{\longrightarrow}\frac{1}{t}(\alpha^t-\1),
\end{align*}
where the convergence occurs with respect to the norm-topology. The group $\big(\alpha^t\big)_{t\in\R}$ is then, due to
\begin{gather*}
\lim_{h\rightarrow 0}\Big\|\frac{1}{h}\big(\alpha^h-\1\big)-\frac{1}{t}\big(\alpha^t-\1\big)X_t^{-1}\Big\|=0,
\end{gather*}
uniformly differentiable and we identify our derivation as
\begin{gather*}
\delta=\frac{1}{t}\big(\alpha^t-\1\big)X_t^{-1}.
\end{gather*}
This can be justified by the relation
\begin{gather*}
\alpha^t=\delta\int_0^t\alpha^sds+1,
\end{gather*}
which finally leads via iteration to  
\begin{gather*}
\alpha^t=e^{t\delta}=\sum_{n=0}^\infty\frac{t^n}{n!}\delta^n.
\end{gather*}
 The last statement, namely the existence of the selfadjoint operator $H$, is verified via Proposition \ref{Neumann-derivation}.
\qed

\begin{theo}
Let $\delta$ be a norm-densely defined and norm-closed element of the unital $C\*$-algebra $\Aa$ with domain $D(\delta)$. Then $\delta$ is the infinitesimal generator of a strongly continuous one-parameter ${}^*$-automorphism group $\big(\alpha^t\big)_{t\in\R}$ of $\Aa$ if and only if the following conditions are satisfied:
\begin{itemize}
\item[(i)] $\delta$ is a derivation and $D(\delta)$ is a ${}^*$-algebra.
\item[(ii)] $\delta$ has a dense set of analytic elements.
\item[(iii)] $\|(\1+\alpha\delta)(A)\|\geq\|A\|$ for all $A\in D(\delta)$ and $\alpha\in\R$.
\end{itemize}
 
\end{theo}

\proof
If we assume the group to be strongly continuous, then its infinitesimal generator is obviously a symmetric derivation with its ${}^*$-algebra as a domain. If $A$ is an analytic element for $\alpha^t$ in the strip $\{z\in\C|\;|\Im z|<\lambda\}$, then, via Cauchy estimates and due to the fact that each ${}^*$-automorphism of a $C\*$-algebra is norm-preserving, one obtains for some constant $M$:
\begin{gather*}
\Big\|\frac{d^n}{dt^n}\alpha^t(A)\Big\|=\|\alpha^t\delta^nA\|=\|\delta^nA\|\leq\frac{n!M}{\lambda^n}.
\end{gather*}
But this means analyticity of $A$ for $\delta$, i.e.
\begin{gather*}
\sum_{n=0}^\infty\frac{|z|^n}{n!}\|\delta^nA\|\leq M\sum_{n=0}^\infty\left(\frac{|z|}{\lambda}\right)^n<\infty,
\end{gather*}
and, since $\alpha^t$ has a dense set of analytic elements, one arrives at the statement (ii). \\
Contrariwise, the conditions (ii) and (iii) force $\delta$ to be the generator of a strongly continuous group of isommetries $\alpha^t$ satisfying, due to (iii),
\begin{gather*}
\alpha^t(\1)=\1+t\delta(\1)+o(t).
\end{gather*}
Because our derivation is symmetric, i.e. $\delta(\1)$ is selfadjoint, and $\|\alpha^t(\1)\|=1$, one obtains for the spectrum in the case of small positive and small negative $t$
\begin{gather*}
\sigma\big(\delta(\1)\big)\in\langle-\infty,0]\quad\text{and}\quad\sigma\big(\delta(\1)\big)\in[0,\infty\rangle,
\end{gather*}
respectively. Consequently, $\sigma\big(\delta(\1)\big)=0$ and, since $\delta(\1)$ is selfadjoint, it follows that $\delta(\1)=0$. Thus we conclude that $\alpha^t(\1)=\1$ for all $t\in\R$ and every $\alpha^t$ is a ${}^*$-automorphism.
\qed

\begin{prop}\label{Neumann-aut}
Every ${}^*$-homomorphism $\tau:\Ma\longrightarrow\Na$ between two von Neumann algebras $\Ma$ and $\Na$ is $\sigma$-weakly and $\sigma$-strongly continuous.
\end{prop}

\proof
Let us consider the increasing net $(A_n)_{n\in\N}$ in $\Ma_+$ and its least upper bound 
$$A:=\text{l.u.b.}_{n\rightarrow\infty}A_n=\sigma\text{-weak}\lim_{n\rightarrow\infty}A_n,
$$
for which, due to the fact that $\tau$ is surjective and obviously preserves positivity, one has:
$$
\tau(A):=\text{l.u.b.}_{n\rightarrow\infty}\tau(A_n)=\sigma\text{-weak}\lim_{n\rightarrow\infty}\tau(A_n).
$$
Thus $\o\circ\tau$ is normal if $\o$ is so. But since a state on a von Neumann algebra is normal if and only if it is $\sigma$-weakly continuous, and every $\sigma$-weakly continuous functional can be written as a linear combination of $\sigma$-weakly continuous states, a $\sigma$-weakly continuous functional passes on this property to the functional $\o\circ\tau$, i.e. $\tau$ is a $\sigma$-weakly continuous homomorphism. \\
Now, if $(A_n)_{n\in\N}$ converges strongly to 0, then $(A_n^*A_n)_{n\in\N}$ converges $\sigma$-weakly to 0. But this means that 
$$
\sigma\text{-weak}\lim_{n\rightarrow\infty}\tau(A_n)^*\tau(A_n)=\sigma\text{-weak}\lim_{n\rightarrow\infty}\tau(A_n^*A_n)=0,
$$
and, hence, the net $\big(\tau(A_n)\big)_{n\in\N}$ converges $\sigma$-strongly to 0.

\qed

\begin{theo}
Let $\delta$ be a $\sigma(\Ma,\Ma_*)$-densely defined and $\sigma(\Ma,\Ma_*)-\sigma(\Ma,\Ma_*)$-closed element of the unital von Neumann algebra $\Ma$ with domain $D(\delta)$ which contains the identity. Then $\delta$ is the infinitesimal generator of a $\sigma$-weakly continuous one-parameter ${}^*$-automorphism group of $\Ma$ if and only if the following conditions are satisfied:
\begin{itemize}
\item[(i)] $\delta$ is a derivation and $D(\delta)$ is a ${}^*$-algebra.
\item[(ii)] $\delta$ has a dense set of analytic elements.
\item[(iii)] $\|(\1+\alpha\delta)(A)\|\geq\|A\|$ for all $A\in D(\delta)$ and $\alpha\in\R$.
\end{itemize}
 
\end{theo}

\sproof
Because of Proposition \ref{Neumann-aut} and since each automorphism on a $C^*$-algebra is an isometry, each automorphism group on a von Neumann algebra forms a $\sigma$-weakly continuous group of $\mathcal{C}_0^*$-isometries on which the proof is based on. The proof is then carried out in a similar way to the case of $C^*$-algebras.

\qed
The analytic elements for a derivation on a von Neumann algebra form a ${}^*$-algebra.

\section{Modular Theory}

In this section we give a short and straightforward introduction to the Tomita-Takesaki theory, also called modular theory \cite{Takesaki:1970}. The outstanding result is the fact that only few assumptions are needed to formulate this theory, namely an underlying von Neumann algebra $\Ma$ and a cyclic and separating vector in the representing Hilbert space, but, surprisingly, this already ensures the existence of an one-parameter automorphism group of $\Ma$. Thus, one gets the dynamics in $\Ma$ for free.\\
Although one can formulate it in a more general setting, namely in terms of left or right Hilbert algebras which are identical for isometrical involutions, we will restrict ourselves to what is absolutely necessary for our aims. We will follow the standard literature, e.g. \cite{Takesaki2:2001}, \cite{Stratila:1981}, \cite{Bratteli:1979tw} and \cite{Kadison:1983}. Modular theory has also been investigated in the framework of O$\*$-algebras, i.e. $\*$-algebras of closable operators, see e.g. \cite{Inoue:1998}.\\
 Our starting point is a von Neumann algebra $\mathfrak{M}$ acting on a Hilbert space $\Ha$  and a vector $\xi\in\Ha$ which cyclic, i.e. $\mathfrak{M}\xi$ is dense in $\Ha$, and separating, i.e. $A\xi=0$ implies $A=0$ for $A\in\mathfrak{M}$. Because due to Corollary \ref{cyclicseperating}, cyclicity of $\xi$ for $\Ma$ is equivalent to it being separable for the commutant $\mathfrak{M}'$, the vector $\xi$ transports these two properties from the algebra onto its commutant. Thus the following two anti-linear operators are well defined:
\begin{align*}
S_0:D(S_0)=\mathfrak{M}\xi\subset\Ha&\longrightarrow\Ha\\
A\xi&\mapsto S_0 A\xi:=A^*\xi,\\
F_0:D(F_0)=\mathfrak{M}'\xi\subset\Ha&\longrightarrow\Ha\\
A\i&\mapsto F_0 A\xi:=A^*\xi.
\end{align*}
Both operators are closable, and one defines $F:=\overline{F_0}=S_0^*$ and $S:=\overline{S_0}=F_0^*$. Therefore, the Tomita operator $S$ allows for a unique polar decomposition into the positive, selfadjoint operator $\Delta$ and the anti-unitary operator $J$, the so-called modular operator and modular conjugation with respect to the pair $(\mathfrak{M},\xi)$, respectively:
\begin{equation*}
S=J\Delta^{1/2}.
\end{equation*}
The first often used properties of the modular objects are collected in the following proposition.
\begin{prop}
For the modular objects the following relations hold:
\begin{align*}
&(i)\;\Delta=FS;                     &&(ii)\;S^{-1}=S;\\
&(iii)\;J=J^2;                       &&(iv)\;J=J^*;\\
&(v)\;\Delta^{-1/2}=J\Delta^{1/2}J;  &&(vi)\;F=J\Delta^{-1/2};\\
&(vii)\;SF=\Delta^{-1};              &&(viii)\;S=\Delta^{-1/2}J.
\end{align*}
\end{prop}

\proof
\begin{itemize}
\item[(i)] $\Delta=S^*S=FS$.
\item[(ii)] This statement follows as $S$ is the closure of $S_0=S_0^{-1}$.
\item[(iii)] Thanks to $(ii)$ we conclude $J\Delta^{1/2}=S=S^{-1}=\Delta^{-1/2}J^*$ and therefore 
\begin{gather}\label{modularobjects}
J^2\Delta^{1/2}=J\Delta^{-1/2}J^*.
\end{gather}
The uniqueness of the polar decomposition ends the proof.
\item[(iv)] The assertion follows directly from Equation \ref{modularobjects}.
\item[(v)] The claim is a direct consequence of Equation \ref{modularobjects}.
\item[(vi)] $F=S^*=(\Delta^{-1/2})^*=J\Delta^{-1/2}$.
\item[(vii)] $SF=J\Delta^{1/2}J\Delta^{-1/2}=\Delta^{-1}$.
\item[(vii)] $S=J\Delta^{1/2}=J\Delta^{1/2}JJ=\Delta^{-1/2}J$.
\end{itemize}
\qed
We are now already in the position to formulate the core of modular theory.
\begin{theo}\label{TomitaTakesaki}
For the von Neumann algebra $\mathfrak{M}$ and the associated modular operator and modular conjugation the relations
\begin{gather*}
J\Ma J=\Ma',\quad J\Ma'J=\Ma\\
\text{and}\quad\Delta^{it}\Ma\Delta^{-it}=\mathfrak{M}
\end{gather*}
hold for all $t\in\R$.
\end{theo}

\proof
Let  $\lambda>0$, $\varphi,\psi\in D\big(\Delta^{1/2}\big)\cap D\big(\Delta^{-1/2}\big)$ and consider the functions
\begin{align}\label{Ilambda}
I_\lambda:\LH&\longrightarrow\LH\notag\\
A&\mapsto I_\lambda(A):=\lambda^{-1/2}\int_{-\infty}^\infty\frac{\lambda^{it}}{e^{\pi t}+e^{-\pi t}}\Delta^{it}A\Delta^{-it}dt.
\end{align}
and
\begin{align*}
f(\lambda):&=\big\langle\psi,\Delta^{-1/2}I_\lambda(A)\Delta^{1/2}+\lambda\Delta^{1/2}I_\lambda(A)\Delta^{-1/2}\varphi\big\rangle\\
&=\big\langle\Delta^{-1/2}\psi,I_\lambda(A)\Delta^{1/2}\varphi\big\rangle+\lambda\big\langle\Delta^{1/2}\psi,I_\lambda(A)\Delta^{-1/2}\varphi\big\rangle\\
&=\int_{-\infty}^\infty\frac{\lambda^{it}}{e^{\pi t}+e^{-\pi t}}\Big[\lambda^{-1/2}\big\langle\Delta^{-1/2-it}\psi,A\Delta^{1/2-it}\varphi\big\rangle\\
&\hspace{5cm}+\lambda^{1/2}\big\langle\Delta^{1/2-it}\psi,A\Delta^{-1/2-it}\varphi\big\rangle\Big]dt.\\
&=\int_{-\infty}^\infty\frac{\lambda^{it}}{e^{\pi t}+e^{-\pi t}}\int\left(\frac{\mu}{\rho}\right)^{it}\left[\left(\frac{\rho}{\mu\lambda}\right)^{1/2}+\left(\frac{\mu\lambda}{\rho}\right)^{1/2}\right]\\
&\hspace{7cm}d^2\big\langle E_\Delta(\mu)\psi,AE_\Delta(\rho)\varphi\big\rangle dt\\
&=\int\left(\frac{\mu}{\rho}\right)^{it}\left[\left(\frac{\rho}{\mu\lambda}\right)^{1/2}+\left(\frac{\mu\lambda}{\rho}\right)^{1/2}\right]\int_{-\infty}^\infty\frac{dt}{e^{\pi t}+e^{-\pi t}}\left(\frac{\lambda\mu}{\rho}\right)^{it}\\
&\hspace{7cm}d^2\big\langle E_\Delta(\mu)\psi,AE_\Delta(\rho)\varphi\big\rangle\\
&=\int d^2\big\langle E_\Delta(\mu)\psi,AE_\Delta(\rho)\varphi\big\rangle\\
&=\langle\psi,A\varphi\rangle,
\end{align*}
where we have made use of spectral decomposition of the modular operator
\begin{gather*}
\Delta=\int\mu\, dE_\Delta(\mu),
\end{gather*}
and the identity
\begin{gather*}
\int_{-\infty}^\infty\frac{e^{ipt}}{e^{\pi t}+e^{-\pi t}}dt=\frac{1}{e^{p/2}+e^{-p/2}}.
\end{gather*}
Thus we have shown 
\begin{gather*}
A=\Delta^{-1/2}I_\lambda(A)\Delta^{1/2}+\lambda\Delta^{1/2}I_\lambda(A)\Delta^{-1/2}.
\end{gather*}
Because $I_\lambda$ and $(D^{1/2}+\lambda D^{-1/2})$, $D^{1/2}(A):=\Delta^{1/2}A\Delta^{-1/2}$, commute on $\LH$ for all $\A\in\Ma$ and $\lambda\in\R$, one has
\begin{gather*}
I_\lambda=(D^{-1/2}+\lambda D^{1/2})^{-1}.
\end{gather*}
For $A'\in\Ma'$ define
\begin{gather*}
A_\lambda^*\xi:=(\Delta+\lambda\1)^{-1}A'\xi
\end{gather*}
which can be written in terms of $I_\lambda$ as:
\begin{align*}
A_\lambda&=(\Delta+\lambda\1)^{-1}(\Delta+\lambda\1)A_\lambda\\
&=(\Delta+\lambda\1)^{-1}\big(\Delta^{-1/2}A_\lambda\Delta^{1/2}+\Delta^{1/2}A_\lambda\Delta^{-1/2}\big)\\
&=I_\lambda(JA'J),
\end{align*}
i.e. $A_\lambda$ is an element of $\Ma$. For $A'\in\Ma'$ one has
\begin{gather*}
\big\langle\psi,[A',I_\lambda(JA'J)]\varphi\big\rangle=0
\end{gather*}
and therefore via Definition \ref{Ilambda}
\begin{gather*}
\int_{-\infty}^\infty\frac{e^{ipt}}{e^{\pi t}+e^{-\pi t}}\big\langle\psi,[A',\Delta^{it}JA'J\Delta^{-it}]\varphi\big\rangle dt=0,
\end{gather*}
where we have set $\lambda=e^p$ with arbitrary $p\in\R$. The application of the Fourier transformation leads to 
\begin{gather}\label{aut}
\Delta^{it}JA'J\Delta^{-it}\in\Ma'',
\end{gather}
but $\Ma$ is a von Neumann algebra and therefore we get for $t=0$
\begin{gather*}
J\Ma' J\subseteq\Ma.
\end{gather*}
Due to Proposition \ref{modularobjects}, the pairs $(\Ma,\xi)$ and $(\Ma',\xi)$ share the same modular conjugation and therefore we have also shown
\begin{gather*}
J\Ma J\subseteq\Ma'.
\end{gather*}
By Proposition \ref{modularobjects}, we finally conclude
\begin{gather*}
\Ma=J^2\Ma J^2\subseteq J\Ma'J\subseteq\Ma,
\end{gather*}
i.e.
\begin{gather*}
J\Ma J=\Ma'\quad\text{and}\quad J\Ma'J=\Ma.
\end{gather*}
We have seen that each element of $\Ma$ can be written as $A=JA'J$ for some $A'\in\Ma'$ and with Equation \ref{aut} the last statement of the theorem follows:\begin{gather*}
\Delta^{it}\Ma\Delta^{-it}=\Ma,\qquad t\in\R.
\end{gather*}

\qed

\begin{defi}
A von Neumann algebra is said to be $\sigma$-finite if it contains (at most) countably many pairwise orthogonal projections. 
\end{defi}
In statistical quantum mechanics and quantum field theory only $\sigma$-finite  von Neumann algebras appear, which therefore can be represented in a separable Hilbert space, while von Neumann algebras, which can be represented in a separable Hilbert space, need not be $\sigma$-finite in general.

\begin{lem}
For the von Neumann algebra $\Ma$ acting on the Hilbert space $\Ha$ the following statements are equivalent:
\begin{itemize}
\item[(i)] $\Ma$ is $\sigma$-finite.
\item[(ii)] There exists a countable subset of $\Ha$, which is separating for $\Ma$.
\item[(iii)] There exists a faithful and normal weight on $\Ma$.
\item[(iv)] $\Ma$ is isomorphic to a von Neumann algebra $\pi(\Ma)$ which admits a cyclic and separating vector.
\end{itemize}
\end{lem}

\proof
Let us assume $(i)$ and consider the existence of an maximal family $\{\xi_n\}$ of vectors $\xi_i\in\Ha$ such that $[\Ma'\xi_i]$ are pairwise orthogonal. $[\Ma'\xi_i]$ constitutes a projection in $\Ma$, i.e. $\{\xi_n\}$ must be countable and thanks to maximality one obtains
\begin{gather*}
\sum_i[\Ma'\xi_i]=\1,
\end{gather*}
i.e. $\{\xi_n\}$ is cyclic for $\Ma'$ and, due to Corollary \ref{cyclicseperating}, it has to be separating for $\Ma$.\\
Now let the sequence $(\xi_n)$ be separating for $\Ma$ with the property $\sum_i\|\xi \|^2=1$. The state $\o(A):=\sum_i\langle\xi_i,A\xi_i\rangle$ on $\Ma$ is $\sigma$-weakly continuous, thus normal. Moreover, the assumption $\o(A^*A)=0$ leads via
\begin{align*}
\|A\xi \|^2=\langle\xi_i,A^*A\xi_i\rangle=0,\qquad\forall\xi_i,
\end{align*}
to $A=0$, namely to faithfulness of $\o$.\\
We suppose the existence of a faithful and normal sate on the von Neumann algebra $\Ma$ and the corresponding cyclic representation $(\Ha,\pi,\xi)$. By Lemma \ref{homomorphism}$(v)$, $\pi(\Ma)$ is a $C^*$-algebra and due to 
\begin{gather*}
\pi(\Ma)=\pi(\Ma'')=\pi(\Ma)''
\end{gather*}
also a von Neumann algebra. The following implications for all $A\in\Ma$ show faithfulness for $\pi$ and the separating property for $\xi$ with respect to $\pi(\Ma)$ and therefore the validity of $(iv)$:
\begin{gather*}
\pi(A)=0\quad\Longrightarrow\quad\o(A^*A)=\|\pi(A)\xi\|^2\quad\Longrightarrow\quad A^*A=A=0.
\end{gather*}
Finally, we suppose $(iv)$ and consider the separating and cyclic vector $\xi\in\Ha$ for $\pi(\Ma)$ and a family $\{P_n\}$ of pairwise orthogonal projections in $\Ma$. With the abbreviation $P:=\sum_iP_i$ we obtain
\begin{align*}
\|\pi(P)\xi\|^2&=\big\langle\pi(P)\xi,\pi(P)\xi\big\rangle\\
&=\sum_{i,j}\big\langle\pi(P_i)\xi,\pi(P_j)\xi\big\rangle\\
&=\sum_{i,j}\|\pi(P_i)\xi\|^2<\infty,
\end{align*}
i.e. only a countable number of $\pi(P_i)\xi$ and thus of the projections $P_i$ can be nonzero, the statement of $(i)$.
\qed
Given a faithful and normal weight $\o$, one can derive the associated cyclic representation $(\Ha_{\o},\pi_{\o},\xi_{\o})$ through the GNS-construction and the modular operator for the pair $\big(\pi_{\o}(\mathfrak{M}),\xi_{\o}\big)$. The above theorem ensures the existence of a  $\sigma$-weakly continuous one-parameter group of $\*$-automorphisms, $\big(\sigma_{\o}^t\big)_{t\in\R}$,
\begin{align*}
\sigma_{\o}^t:\mathfrak{M}&\longrightarrow\mathfrak{M}\\
A&\mapsto\sigma_{\o}^t(A):=\pi_{\o}^{-1}\big(\Delta^{it}\pi_{\o}(A)\Delta^{-it}\big),
\end{align*}
the so-called modular automorphism group associated with $\big(\pi_{\o}(\mathfrak{M}),\xi_{\o}\big)$. Since we are concerned with von Neumann algebras, due to Proposition \ref{Neumann-aut}, the modular group is continuous with respect to the $\sigma$-strong topology, too. The modular group is a powerful and constructive tool for the investigation of von Neumann algebras and has made possible many applications in mathematics and theoretical physics. The main linkage between modular theory and physics is the following property, 
\begin{align*}
\big\langle\Delta^{1/2}\pi_\o(A)\xi_\o,\Delta^{1/2}\pi_\o(B)\xi_\o\big\rangle&=\big\langle J\pi_\o(A^*)\xi_\o,J\pi_\o(B^*)\xi_\o\big\rangle\\
 &=\big\langle\pi_\o(B^*)\xi_\o,\pi_\o(A^*)\xi_\o\big\rangle,
\end{align*}
the so-called modular condition or KMS-condition, which can equivalently be described by means of the modular group itself,
\begin{gather}\label{kmscondition}
\o\big(\sigma_\o^{i\beta/2}(A)\sigma_\o^{-i\beta/2}(B)\big)=\o(BA)
\end{gather}
for all $A,B\in\Ma$. A state $\o$ satisfying condition \ref{kmscondition} is called $(\sigma,\beta)$-KMS state, where the parameter $\beta=\frac{1}{T}$ is interpreted in statistical physics as the inverse of the temperature $T$.\\
A spatial derivation $\delta=i[H,A]$ of a von Neumann algebra $\Ma$ implemented by a self-adjoint operator $H$ can be extended such that it generates a $\*$-weak-continuous one-parameter group of $\*$-automorphisms
\begin{align*}
\alpha:\R&\longrightarrow\Aut(\Ma)\\
t&\mapsto\alpha^t,\quad \alpha^t(A):=e^{itH}Ae^{-itH}.
\end{align*} 
Let $\Delta_\xi$ be the modular operator with respect to the cyclic and separating vector $\xi\in\Ha$, and, with the help of the self-adjoint operator $H\in\Ma$, $H\xi=0$, define the domain
$$
\mathbf{D}(\delta):=\big\{A\in\Ma|\;i[H,A]\in\Ma\big\}.
$$ 
Then, the following two conditions are equivalent:
\begin{itemize}
\item[(i)] $e^{itH}\Ma e^{-itH}=\Ma$ for all $t\in\R$. 
\item[(ii)] $\mathbf{D}(\delta)\xi$ is a core for $H$, and $H$ and $\Delta$ commute strongly, i.e.
\begin{gather*}
\Delta^{it}H\Delta^{-it}=H,\qquad\forall t\in\R.
\end{gather*}
\end{itemize}
If there exists more than one faithful, normal and semifinite state on $\Ma$, it can be shown that the modular automorphism is unique up to unitaries. Let $\o_1$ and $\o_2$ be two faithful, normal and semifinite states on $\Ma$ and $\sigma^t_{\o_1}$ and $\sigma^t_{\o_2}$ the corresponding modular groups. Let us consider a faithful and normal weight defined as
\begin{equation*}
\rho\left[\begin{pmatrix}A_{11} & A_{12}\\A_{21} & A_{22}\end{pmatrix}\right]:=\frac{1}{2}\big(\o_1(A_{11})+\o_2(A_{22})\big)
\end{equation*}
on $\Ma\otimes M_2$ and the corresponding modular group $\sigma^t_{\rho}\,$, then the unitaries $\Gamma_t$, defined as 
\begin{equation*}
\begin{pmatrix}0 & \Gamma_{t}\\ 0 & 0\end{pmatrix}:=\sigma_{\rho}^t\left[\begin{pmatrix} 0 & \1\\ 0 & 0\end{pmatrix}\right],
\end{equation*}
connect the two original modular groups. The next theorem is dealing exactly with this situation, but first we have to introduce the following notions. We denote by $\A(\mathbf{D})$ the set of bounded and holomorphic functions on the domain $\mathbf{D}$ and define
$$
\mathrm{n}_\o:=\big\{A\in\Ma|\;\o(A^*A)<\infty\big\}.
$$

\begin{theo}[Connes' Cocycle Derivative]\label{Connes-theorem}
If  $\o_1$, $\o_2$ are two faithful, normal and semifinite states on $\mathfrak{M}$, then there exists a $\sigma$-strongly continuous one-parameter family of unitaries $\big(\Gamma_{t}\big)_{t\in\R}$ in $\mathfrak{M}$ with the following properties:

\begin{itemize}
\item [(i)]$\sigma_{\o_2}^t(A)=\Gamma_{t}\sigma_{\o_1}^t(A)\Gamma_{t}^*,\quad A\in\mathfrak{M},\;t\in\R$;
\item [(ii)] $\Gamma_{s+t}=\Gamma_{s}\sigma_{\o_1}^s(\Gamma_{t}),\quad s,t\in\R$.
\end{itemize} 
\end{theo}

\proof
Consider a faithful and normal state defined $\rho$ on the tensor product $\Ma\otimes M(2\times2, \C)$ defined by 
\begin{gather*}
\rho\left[\left(\begin{array}{cc}A_{11} & A_{12}\\ A_{21} & A_{22}\end{array}\right)\right]:=\frac{1}{2}[\o_1(A_{11})+\o_2(A_{22})]
\end{gather*}
for which one has:
\begin{gather*}
\rho\left[\left(\begin{array}{cc} 1 & 0\\ 0 & 0\end{array}\right)\left(\begin{array}{cc}A_{11} & A_{12}\\ A_{21} & A_{22}\end{array}\right)\right]=\frac{1}{2}\o_1(a_{11})=\rho\left[\left(\begin{array}{cc}A_{11} & A_{12}\\ A_{21} & A_{22}\end{array}\right)\left(\begin{array}{cc} 1 & 0\\ 0 & 0\end{array}\right)\right],\\
\rho\left[\left(\begin{array}{cc} 0 & 0\\ 0 & 1\end{array}\right)\left(\begin{array}{cc}A_{11} & A_{12}\\ A_{21} & A_{22}\end{array}\right)\right]=\frac{1}{2}\o_1(A_{22})=\rho\left[\left(\begin{array}{cc}A_{11} & A_{12}\\ A_{21} & A_{22}\end{array}\right)\left(\begin{array}{cc} 0 & 0\\ 0 & 1\end{array}\right)\right].
\end{gather*}
Since the centralizer of $\o$ is contained in the fixpoint algebra of its modular group $\sigma_\rho$ we conclude $\sigma_\rho^t(\Ma\oplus\Ma)\subseteq\Ma\oplus\Ma$ with
\begin{gather*}
\Ma\oplus\Ma=\left\{\left(\begin{array}{cc}A_{11} & 0\\ 0 & A_{22}\end{array}\right)| A_{ii}\in\Ma\right\}
\end{gather*}
The state $\rho$ can be shown to be a KMS-state with respect to the automorphism group $\sigma_\rho^t$ and due to uniqueness of such states one deduces
\begin{align*}
\left(\begin{array}{cc}A_{11} & 0\\ 0 & A_{22}\end{array}\right)\mapsto\sigma_\rho^t\left[\left(\begin{array}{cc}A_{11} & 0\\ 0 & A_{22}\end{array}\right)\right]=\left(\begin{array}{cc}\sigma_{\o_1}^t(A_{11}) & 0\\ 0 & \sigma_{\o_2}^t(A_{22})\end{array}\right)
\end{align*}
The following implication
\begin{gather*}
\left(\begin{array}{cc} 0 & 0\\ 1 & 0\end{array}\right)\left(\begin{array}{cc} 0 & 0\\ 0 & 1\end{array}\right)=0=\left(\begin{array}{cc} 1 & 0\\ 0 & 0\end{array}\right)\left(\begin{array}{cc} 0 & 0\\ 1 & 0\end{array}\right)\\
\Longrightarrow\quad\sigma_\rho^t\left[\left(\begin{array}{cc} 0 & 0\\ 1 & 0\end{array}\right)\right]\left(\begin{array}{cc} 0 & 0\\ 0 & 1\end{array}\right)=0=\left(\begin{array}{cc} 1 & 0\\ 0 & 0\end{array}\right)\sigma_\rho^t\left[\left(\begin{array}{cc} 0 & 0\\ 1 & 0\end{array}\right)\right]
\end{gather*}
ensures the existence of a $\sigma$-weakly continuous one-parameter family $\Gamma_t\in\Ma$ defined as
\begin{gather*}
\left(\begin{array}{cc} 0 & 0\\ \Gamma_t & 0\end{array}\right):=\sigma_\rho^t\left[\left(\begin{array}{cc} 0 & 0\\ 1 & 0\end{array}\right)\right].
\end{gather*}
Unitarity for all $\Gamma_t$, i.e. $\Gamma_t\Gamma_t^*=\Gamma_t^*\Gamma_t=\1$ for all $t\in\R$, follows directly from the application of $\sigma_\rho^t$ to 
\begin{gather*}
\left(\begin{array}{cc} 0 & 0\\ 1 & 0\end{array}\right)\left(\begin{array}{cc} 0 & 0\\ 1 & 0\end{array}\right)^*+\left(\begin{array}{cc} 0 & 0\\ 1 & 0\end{array}\right)^*\left(\begin{array}{cc} 0 & 0\\ 1 & 0\end{array}\right)=\1.
\end{gather*}
If we apply $\sigma_\rho^t$ to
\begin{gather*}
\left(\begin{array}{cc} 0 & 0\\ 0 & A\end{array}\right)=\left(\begin{array}{cc} 0 & 0\\ 1 & 0\end{array}\right)\left(\begin{array}{cc} A & 0\\ 0 & 0\end{array}\right)\left(\begin{array}{cc} 0 & 1\\ 0 & 0\end{array}\right)
\end{gather*}
we finally arrive at
\begin{gather*}
\sigma_{\o_2}^t(A)=\Gamma_{t}\sigma_{\o_1}^t(A)\Gamma_{t}^*,\quad A\in\mathfrak{M},\;t\in\R.
\end{gather*}
The last statement is justified by:
\begin{align*}
\left(\begin{array}{cc} 0 & 0\\ \Gamma_{s+t} & 0\end{array}\right)&=\sigma_\rho^{s+t}\left[\left(\begin{array}{cc} 0 & 0\\ 1 & 0\end{array}\right)\right]=\sigma_\rho^t\left[\left(\begin{array}{cc} A & 0\\ \Gamma_s & 0\end{array}\right)\right]\\
&=\sigma_\rho^t\left[\left(\begin{array}{cc} 0 & 0\\ 1 & 0\end{array}\right)\right]\sigma_\rho^t\left[\left(\begin{array}{cc} \Gamma_s & 0\\ 0 & 0\end{array}\right)\right]\\
&=\left(\begin{array}{cc} A & 0\\ \Gamma_t\sigma_{\o_1}^t(\Gamma_s) & 0\end{array}\right).
\end{align*}
\qed 
The main statement of this theorem is that two arbitrary modular automorphism groups are equivalent up to inner automorphisms. The specific problem with its statement is that it assures the existence of the unitaries only, but fails to give a method for their construction.\\
The above theorem still holds true with few adjustments if the faithfulness is given only for one of the weights. The family of unitaries is called cocycle derivative of $\o_1$ with respect to $\o_2$ and denoted by $(D_{\o_2}:D_{\o_1})_t:=\Gamma_t$. For two faithful, normal and semifinite states we have 
$$
(D_{\o_2}:D_{\o_1})_t=(D_{\o_2}:D_{\o_2})_t^{-1},\quad t\in\R.
$$
If an additional third faithful, normal and semifinite state on $\Ma$ is given, then one can establish the chain rule,
\begin{equation*}
(D_{\o_1}:D_{\o_2})_t=(D_{\o_1}:D_{\o_3})_t(D_{\o_3}:D_{\o_2})_t, \quad t\in\R.
\end{equation*}
Consequently, the equivalence of two cocycle derivatives of $\o_1$ with respect to $\o_2$ and $\o_3$ uniquely determines the identity of  $\o_2$ and $\o_3$:
$$
(D_{\o_1}:D_{\o_2})_t=(D_{\o_1}:D_{\o_3})_t\,,\quad\forall t\in\R\quad\Longleftrightarrow\quad\o_2=\o_3.
$$
One of the most important mathematical applications of the Tomita-Takesaki theory is the classification of factors. If for a fixed $t_0\in\R$ and a particular faithful, normal and semifinite state $\o$ its modular group $\sigma^{t_0}_{\o}$ is inner, then, because of the cocycle theorem, $\sigma_{\o'}^{t_0}$ is inner for any faithful, normal and semifinite weight $\o'$. Thus the modular period group 
\begin{equation*}
T(\mathfrak{M}):=\big\{t\in\R|\;\sigma^{t}_{\o} \text{ is inner}\big\}
\end{equation*}
characterises the von Neumann algebra. Let $\o$ be an arbitrary faithful and semifinite state, then $T(\Ma)$ is related to the so-called modular spectrum of $\Ma$, 
\begin{equation}\label{modularspectrum}
S(\mathfrak{M}):=\underset{\o}{\bigcap}\;{\bf Spec}\Delta_{\o},
\end{equation}
via the inclusion
$$
\ln\big(S(\Ma)\backslash\{0\}\big)\subset\big\{s\in\R|\;e^{ist}=1\quad\forall t\in T(\Ma)\big\}.
$$
This modular spectrum can be used for the classification of factors. For this purpose let $\Ma$ be a factor, then the following statements hold:
\begin{itemize}
\item[(i)] $\mathfrak{M}$ is of type $I$ or type $II$, if $S(\mathfrak{M})=\{1\}$;
\item[(ii)] $\mathfrak{M}$ is of type $III_0$, if $S(\mathfrak{M})=\{0,1\}$;
\item[(iii)] $\mathfrak{M}$ is of type $III_\lambda$, if $S(\mathfrak{M})=\{0\}\cup\{\lambda^{n}|\;0<\lambda<1,n\in \mathbb{Z}\}$;
\item[(iv)] $\mathfrak{M}$ is of type $III_1$, if $S(\mathfrak{M})=\R_+$.
\end{itemize}

\section{Non-Commutative Integration and \\Crossed Products}

The main target of this section is the introduction to crossed products, in particular the statement of Theorem \ref{generalTakesaki}, since it enables one always with a spatial automorphism group on the corssed product $\Ma\otimes_\alpha\Gr$. Therefore, we may transfer in case of general (non-spatial) automorphism groups our investigation from the von Neumann algebra $\Ma$ onto the crossed product and make use of results already obtained for the spatial case. \\
First of all, we define a new faithful, semifinite and normal state on $\Na_\rho:=\Ma\otimes M(2\times2,\C)$ with reference to the von Neumann algebra $\Ma$, the so-called balanced weight of the weights $\o_1$ and $\o_2$ on $\Ma$, by
\begin{equation}\label{balancedweight}
\rho\begin{pmatrix}A_{11} & A_{12}\\ A_{21} & A_{22}\end{pmatrix}:=(\o_1\oplus\o_2)\begin{pmatrix}A_{11} & A_{12}\\ A_{21} & A_{22}\end{pmatrix}:=\o_1(A_{11})+\o_2(A_{22}).
\end{equation}
We denote the so-called semi-cyclic representation of $\Na$ by $(\Ha_\rho,\pi_\rho,\eta_\rho)$ consisting of a representation $(\Ha_\rho,\pi_\rho)$ and the linear map 
\begin{align*}
\eta_\rho:\mathfrak{n}&\longrightarrow\Ha_\rho\\
A&\mapsto\eta_\rho(A)
\end{align*} 
from a left ideal $\mathfrak{n}$ of $\Ma$ into $\Ha_\rho$ such that
\begin{gather*}
\pi_\rho(B)\eta_\rho(A)=\eta_\rho(BA),\quad\forall A\in\mathfrak{n},B\in\Ma.
\end{gather*} 
Because of the properties of $\rho$ mentioned above, we may introduce the corresponding modular objects. One finds that the representation space can be expressed as a direct sum of mutually orthogonal subspaces,
$$
\Ha_\rho=\Ha_{\o_1,\o_1}\oplus\Ha_{\o_1,\o_2}\oplus\Ha_{\o_2,\o_1}\oplus\Ha_{\o_2,\o_2},
$$
with
\begin{align*}
\Ha_{\o_1,\o_1}&:=\Big[\big((\mathfrak{n}_{\o_1}\cap\mathfrak{n}_{\o_1}^*)\otimes e_{11}\big)+\mathbf{N}_\rho\Big],\\
\Ha_{\o_1,\o_2}&:=\Big[\big((\mathfrak{n}_{\o_1}^*\cap\mathfrak{n}_{\o_2})\otimes e_{12}\big)+\mathbf{N}_\rho\Big],\\
\Ha_{\o_2,\o_1}&:=\Big[\big((\mathfrak{n}_{\o_1}\cap\mathfrak{n}_{\o_2}^*)\otimes e_{21}\big)+\mathbf{N}_\rho\Big],\\
\Ha_{\o_2,\o_2}&:=\Big[\big((\mathfrak{n}_{\o_2}\cap\mathfrak{n}_{\o_2}^*)\otimes e_{22}\big)+\mathbf{N}_\rho\Big],
\end{align*}
where $\mathbf{N}_\rho:=\big\{A\in\Na_\rho|\;\rho(A^*A)=0\big\}$ is a left ideal of $\Na_\rho$, $\mathfrak{n}_{\o_1}:=\big\{A\in\Ma|\;\o_1(A^*A)<\infty\big\}$, $\mathfrak{n}_{\o_2}:=\big\{A\in\Ma|\;\o_2(A^*A)<\infty\big\}$
and 
\begin{align*}
e_{11}:=\begin{pmatrix}1 & 0\\ 0 & 0\end{pmatrix},\;e_{12}:=\begin{pmatrix}0 & 1\\ 0 & 0\end{pmatrix},\;e_{21}:=\begin{pmatrix}0 & 0\\ 1 & 0\end{pmatrix},\;e_{22}:=\begin{pmatrix}0 & 0\\ 0 & 1\end{pmatrix}
\end{align*}
are the matrix units. The brackets `$[\;\cdot\;]$' denote the closure in the Hilbert space $\Ha_\rho$. The von Neumann algebra $\Na_\rho$ inherits its involution from $\Ma$. Since we are concerned with states, not weights, the sets simplify to $\mathfrak{n}_{\o_1}=\mathfrak{n}_{\o_2}=\Ma$. Furthermore, because  $\rho$ is faithful, we obtain $\mathbf{N}_\rho=\{0\}$. For the domain of the $\*$-operation on has,
$$
\D^*=\big(\D^*\cap\Ha_{\o_1,\o_1}\big)\oplus\big(\D^*\cap\Ha_{\o_1,\o_2}\big)\oplus\big(\D^*\cap\Ha_{\o_2,\o_1}\big)\oplus\big(\D^*\cap\Ha_{\o_2,\o_1}\big),
$$
which determines the Tomita operator as
\begin{align*}
S\big(\D^*\cap\Ha_{\o_1,\o_1}\big)&=\D^*\cap\Ha_{\o_1,\o_1},\\
S\big(\D^*\cap\Ha_{\o_1,\o_2}\big)&=\D^*\cap\Ha_{\o_2,\o_1},\\
S\big(\D^*\cap\Ha_{\o_2,\o_1}\big)&=\D^*\cap\Ha_{\o_1,\o_2},\\
S\big(\D^*\cap\Ha_{\o_2,\o_2}\big)&=\D^*\cap\Ha_{\o_2,\o_2}.
\end{align*}
Written in a more compact form, we obtain for the Tomita operator the following matrix:
$$
S=
\begin{pmatrix}
S_{11} & 0 & 0 & 0\\ 0 & 0 & S_{23} & 0\\ 0 & S_{32} & 0 & 0\\ 0 & 0 & 0 & S_{44}
\end{pmatrix}.
$$
Now, we are able to identify the Tomita operators as
\begin{equation}\label{Tomitaoperators}
\begin{split}
S_{\o_1}&=U_1^*S_{11}U_1,\quad  S_{\o_1,\o_2}=U_2^*S_{23}U_3,\\
S_{\o_2,\o_1}&=U_3^*S_{32}U_2,\quad \hspace{3,8mm} S_{\o_2}=U_4^*S_{44}U_4,
\end{split}
\end{equation}
where we have set
\begin{equation}\label{operatorU}
\begin{split}
U_1\eta_{\o_1}(A)&:=\eta_\rho(A\otimes e_{11}),\quad\forall A\in\mathfrak{n}_{\o_1},\\
U_2\eta_{\o_1}(A)&:=\eta_\rho(A\otimes e_{21}),\quad\forall A\in\mathfrak{n}_{\o_1},\\
U_3\eta_{\o_2}(A)&:=\eta_\rho(A\otimes e_{12}),\quad\forall A\in\mathfrak{n}_{\o_2},\\
U_4\eta_{\o_2}(A)&:=\eta_\rho(A\otimes e_{22}),\quad\forall A\in\mathfrak{n}_{\o_2}.
\end{split}
\end{equation}
The adjoint of $S$ is then given by
$$
F=
\begin{pmatrix}
F_{11} & 0 & 0 & 0\\ 0 & 0 & F_{23} & 0\\ 0 & F_{32} & 0 & 0\\ 0 & 0 & 0 & F_{44}
\end{pmatrix},
$$
where the components are fixed as
$$
F_{11}:=S_{11} ^*,\quad F_{23}:=S_{32}^*,\quad F_{32}:=S_{23}^*\quad\text{and}\quad F_{44}:=S_{44}.
$$
Finally, via polar decomposition
$$
S=J\Delta^{1/2}=\Delta^{-1/2}J,
$$
which is equivalent to
\begin{align*}
S_{11}&=J_{11}\Delta_{11}^{1/2}=\Delta_{11}^{-1/2}J_{11},\\
S_{23}&=J_{23}\Delta_{33}^{1/2}=\Delta_{22}^{-1/2}J_{23},\\
S_{32}&=J_{32}\Delta_{22}^{1/2}=\Delta_{33}^{-1/2}J_{32},\\
S_{44}&=J_{44}\Delta_{44}^{1/2}=\Delta_{44}^{-1/2}J_{44},
\end{align*}
we arrive at the modular conjugation and, most important for us, the modular operator:
\begin{gather*}
J=
\begin{pmatrix}
J_{11} & 0 & 0 & 0\\ 0 & 0 & J_{23} & 0\\ 0 & J_{32} & 0 & 0\\ 0 & 0 & 0 & J_{44}
\end{pmatrix}
\quad\text{and}\quad
\Delta=
\begin{pmatrix}
\Delta_{11} & 0 & 0 & 0\\ 0 & \Delta_{22} & 0 & 0\\ 0 & 0 & \Delta_{33} & 0\\ 0 & 0 & 0 & \Delta_{44}
\end{pmatrix}.
\end{gather*}
Let $\Ma$ act on the Hilbert space $\Ha$, then one can find linear maps 
\begin{gather*}
\xi:\mathfrak{n}_{\o_1}\ni A\mapsto\xi(A)\in\Ha\quad\text{and}\\
\eta:\mathfrak{n}_{\o_2}\ni A\mapsto\eta(A)\in\Ha
\end{gather*}
such that for all $A\in\Ma$, $\;A_1,B_1\in\mathfrak{n}_{\o_1}$ and $\;A_2,B_2\in\mathfrak{n}_{\o_2}$
\begin{align*}
A\xi(A_1)&=\xi(AA_1),&\quad\o_1(A_1^*B_1)&=\big(\xi(A_1),\xi(B_1)\big),\\
A\eta(A_2)&=\eta(AA_2),&\quad\o_2(A_2^*B_2)&=\big(\eta(A_2),\eta(B_2)\big],
\end{align*}
\begin{align*}
\text{and}\quad\Ha=\big[\xi(\mathfrak{n}_{\o_1})\big]=\big[\eta(\mathfrak{n}_{\o_2})\big].
\end{align*}
The form of the representation space with respect to $\rho$ is reduced to
\begin{gather*}
\Ha_\rho=\Ha\oplus\Ha\oplus\Ha\oplus\Ha,
\end{gather*}
and the representation itself to:
\begin{gather*}
\pi_\rho\begin{pmatrix} A_{11} & A_{12}\\ A_{21} & A_{22}\end{pmatrix}=\begin{pmatrix}A_{11} & A_{12} & 0 & 0\\ A_{21} & A_{22} & 0 & 0\\ 0 & 0 & A_{11} & A_{12}\\ 0 & 0 & A_{21} & A_{22}\end{pmatrix}.
\end{gather*}
Since the modular automorphism group on $\mathfrak{N}_\rho$ should leave the representation space invariant, i.e. 
\begin{gather*}
\sigma_\rho^t\big(\pi_\rho(\mathfrak{N}_\rho)\big)=\Delta^{it}\pi_\rho(\mathfrak{N}_\rho)\Delta^{-it}=\pi_\rho(\mathfrak{N}_\rho),
\end{gather*}
we conclude
\begin{align*}
\Delta_{11}^{it}A_{11}\Delta_{11}^{-it}&=\Delta_{33}^{it}A_{11}\Delta_{33}^{-it},&\quad\Delta_{11}^{it}A_{12}\Delta_{22}^{-it}&=\Delta_{33}^{it}A_{12}\Delta_{44}^{-it},\\
\Delta_{22}^{it}A_{21}\Delta_{11}^{-it}&=\Delta_{44}^{it}A_{21}\Delta_{33}^{-it},&\quad\Delta_{22}^{it}A_{22}\Delta_{22}^{-it}&=\Delta_{44}^{it}A_{22}\Delta_{44}^{-it}.
\end{align*}
The identification of $\Delta_{11}$ and $\Delta_{44}$ with $\Delta_{\o_1}$ and $\Delta_{\o_2}$, respectively, yields the following expressions for the modular automorphism groups,
\begin{align*}
\sigma_{\o_1}^t(A)&=\Delta_{11}^{it}A\Delta_{11}^{-it}\in\Ma,&\quad\sigma_{\o_1,\o_2}^t(A)&=\Delta_{11}^{it}A\Delta_{22}^{-it}\in\Ma,\\
\sigma_{\o_2,\o_1}^t(A)&=\Delta_{22}^{it}A\Delta_{11}^{-it}\in\Ma,&\quad\sigma_{\o_2}^t(A)&=\Delta_{22}^{it}A\Delta_{22}^{-it}\in\Ma,
\end{align*}
for all $A\in\Ma$. Finally, we can write down the form of the modular group with respect to $\rho$ explicitly:
\begin{equation}
\sigma_\rho^t\left[\begin{pmatrix} A_{11} & A_{12}\\ A_{21} & A_{22}\end{pmatrix}\right]=\begin{pmatrix} \sigma_{\o_1}^t(A_{11})& \sigma_{\o_1,\o_2}^t(A_{12})\\ \sigma_{\o_2,\o_1}^t(A_{21}) & \sigma_{\o_2}^t(A_{22})\end{pmatrix}.
\end{equation}
As already mentioned, the unitary cocycle appearing in Connes' Theorem is defined by
\begin{equation*}
\begin{pmatrix}0 & \Gamma_{t}\\ 0 & 0\end{pmatrix}:=\sigma_{\rho}^t\left[\begin{pmatrix} 0 & \1\\ 0 & 0\end{pmatrix}\right],
\end{equation*}
and therefore one has its explicit structure
$$
\Gamma_{t}=\sigma_{\o_1,\o_2}^t(\1)
$$
for all $t\in\R$.\\
Thus, the modular automorphism group with respect to the state $\o_2$ on $\Ma$ is determined up to a cocycle, i.e. up to a perturbation term.

\begin{defi}
The dual cone $\Pa^\circ$ of a convex cone $\Pa$ in a Hilbert space is defined as the set of all vectors $\eta\in\Ha$ with $\langle\xi,\eta\rangle\geq0$ and
\begin{gather*}
\Pa^\circ:=\big\{\eta\in\Ha|\; \langle\xi,\eta\rangle\geq0,\;\forall\xi\in\Pa\big\}.
\end{gather*} 
$\Pa$ is said to be self-dual if $\Pa=\Pa^\circ$. Let $\Aa$ be a left Hilbert algebra associated to the faithful, semi-finite normal weight $\o$ and $\Aa_0$ the corresponding Tomita algebra, i.e. the algebra equipped with an one-parameter group of automorphisms, then we set
\begin{gather*}
\Pa_\o:=\big\{\xi\xi^*|\;\xi\in\Aa_0\big\}^-,\\
\Pa_\xi:=\big\{(\R_+\xi-\Pa_\o)\cap\Pa_\o|\;\xi\in\Pa_\o\big\}^-.
\end{gather*} 
\end{defi}

\begin{lem}\label{operatorU2}
Let $\o_1$ and $\o_2$ be two faithful, semi-finite normal weights on the von Neumann algebra $\Ma$ and $(\pi_{\o_1},\Ha_{\o_1},\eta_{\o_1})$ and $(\pi_{\o_2},\Ha_{\o_2},\eta_{\o_2})$ their cyclic representation, then there exists uniquely a unitary operator $U_{\o_1,\o_2}:\Ha_{\o_2}\longrightarrow\Ha_{\o_1}$ such that:
\begin{itemize}
\item[(i)] $U_{\o_1,\o_2}\pi_{\o_2}(A)U_{\o_1,\o_2}^*=\pi_{\o_1}(A)$ for all $A\in\Ma$;
\item[(ii)] $U_{\o_1,\o_2}\Pa_{\o_2}=\Pa_{\o_1}$.
\end{itemize}
\end{lem}

\sproof
Let us consider the balanced weight $\rho=\o_1\oplus\o_2$ on the product $\Na_\rho:=\Ma\otimes M(2\times2,\C)$ as defined in \ref{balancedweight} as well as $S_{\o_1}$, $S_{\o_1,\o_2}$, $S_{\o_2,\o_1}$ and $S_{\o_2}$ as defined in Equation \ref{Tomitaoperators} with their polar decompositions:
\begin{align*}
S_{\o_1}&=J_{\o_1}\Delta_{\o_1}^{1/2},&\quad  S_{\o_1,\o_2}&=J_{\o_1,\o_2}\Delta_{\o_1,\o_2}^{1/2},\\
S_{\o_2,\o_1}&=J_{\o_2,\o_1}\Delta_{\o_2,\o_1}^{1/2},&\quad  S_{\o_2}&=J_{\o_2}\Delta_{\o_2}^{1/2}.
\end{align*}
The unitary equivalence between $\pi_{\o_1}$ and $\pi_{\o_2}$ is then established by the operator defined in Equation \ref{operatorU} 
\begin{gather*}
U_{\o_1,\o_2}=J_{\o_1}J_{\o_1,\o_2}=J_{\o_1,\o_2}J_{\o_2}
\end{gather*}
and, moreover, one has:
\begin{gather*}
\Pa_{\o_1}=\big\{\pi_{\o_1}(A)J_{\o_1}\eta_{\o_1}(A)|\;A\in\mathfrak{n}_{\o_1}\cap\mathfrak{n}_{\o_1}^*\big\}^-,\\
\Pa_{\o_2}=\big\{\pi_{\o_2}(A)J_{\o_2}\eta_{\o_2}(A)|\;A\in\mathfrak{n}_{\o_2}\cap\mathfrak{n}_{\o_2}^*\big\}^-.
\end{gather*} 
For all $A\in\mathfrak{n}_{\o_1}\cap\mathfrak{n}_{\o_1}^*$ and $B\in\mathfrak{n}_{\o_2}\cap\mathfrak{n}_{\o_2}^*$ selfduality of $\Pa_{\o_1}$ can be shown to lead to the estimation
\begin{gather*}
\big\langle\pi_{\o_1}(A)J_{\o_1}\eta_{\o_1}(A),U_{\o_1,\o_2}\pi_{\o_2}(B)J_{\o_2}\eta_{\o_2}(B)\big\rangle\geq0,
\end{gather*}
i.e. one obtains the inclusion
\begin{gather*}
U_{\o_1,\o_2}\Pa_{\o_2}\subset\Pa_{\o_1}.
\end{gather*}
Via $U_{\o_1,\o_2}^*=U_{\o_2,\o_1}$ one derives the opposite inclusion
\begin{gather*}
U_{\o_1,\o_2}\Pa_{\o_2}\supset\Pa_{\o_1}. 
\end{gather*}
\qed
For a normal weight $\o$ one may introduce the projections $e$ and $f$ of $\Ma$ by
\begin{gather*}
\Ma e:=\mathfrak{n}_\o^-\quad\text{and}\\
\Ma f:=N_\o=\{A\in\Ma|\;\o(A^*A)=0\}.
\end{gather*}
Since $\o$ is semifinite on $e\Ma e$ and faithful on $(1-f)\Ma(1-f)$, one defines the support of $\o$, denoted by $s(\o)$, as the difference $e-f$.
\begin{lem}\label{orthogonalxi}
For $\xi_1,\xi_2\in\Pa_\o$ the following conditions are equivalent:
\begin{itemize}
\item[(i)] $\xi_1\perp\xi_2$;
\item[(ii)]$\Pa_{\xi_1}\perp\Pa_{\xi_2}$;
\item[(iii)] $s(\o_{\xi_1})\perp s(\o_{\xi_2})$.
\end{itemize}
\end{lem}

\proof
Let us assume $(i)$ then thanks to Lemma \ref{operatorU2} and the definition of $\Pa_\xi$ we only have to prove
\begin{gather*}
(\R_+\xi_1-\Pa_\o)\cap\Pa_\o\perp(\R_+\xi_2-\Pa_\o)\cap\Pa_\o.
\end{gather*}
But for $\xi_1,\eta_2,\,\lambda\xi_1-\eta_1,\,\lambda\xi_2-\eta_2\in\Pa_\o$, where $\lambda,\mu>0$, we obtain
\begin{gather*}
0\leq\langle\eta_1,\eta_2\rangle\leq\lambda\langle\xi_1,\eta_2\rangle\leq\lambda\mu\langle\xi_1,\xi_2\rangle=0,
\end{gather*}
that is $\eta_1\perp\eta_2$. $[\Pa_\xi]=[\Ma\xi]\cap[\Ma'\xi]$ for all $\xi\in\Pa_\o$ proves the implication $(ii)\Longrightarrow(iii)$. The direction $(iii)\Longrightarrow(ii)\Longrightarrow(i)$ follows directly by definition.
\qed

\begin{defi}
Let $(\Ma,\Ha)$ be a von Neumann algebra, then the quadruple $(\Ma,\Ha,J,\Pa)$, where $J$ is the corresponding modular conjugation and $\Pa$ a self-dual cone in $\Ha$, is said to be a standard form, if the following conditions are satisfied:
\begin{itemize}
\item[(i)] $JAJ=A^*,\quad A\in\Za_\Ma$\,;
\item[(ii)] $J\xi=\xi,\quad \forall\xi\in\Pa$;
\item[(iii)] $AJA\Pa\subset\Pa,\quad A\in\Ma$.
\end{itemize} 
\end{defi}

\begin{prop}\label{standardisomorphism}
Let $(\Ma_1,\Ha_1,J_1,\Pa_1)$ and $(\Ma_2,\Ha_2,J_2,\Pa_2)$ be two standard forms and $\pi:\Ma_1\longrightarrow\Ma_2$ an isomorphism, then there exists uniquely a unitary operator $U:\Ha_1\longrightarrow\Ha_2$ such that:
 \begin{itemize}
\item[(i)] $\pi(A)=UAU^*,\quad A\in\Ma_1$;
\item[(ii)] $J_2=UJ_1U^*$;
\item[(iii)] $\Pa_2=U\Pa_1$.
\end{itemize} 
\end{prop}

\proof
The existence of the unitary $U$ is ensured if we prove the isomorphism 
\begin{gather*}
(\Ma,\Ha,J,\Pa)\cong\big(\pi_\o(\Ma),\Ha_\o,J_\o,\Pa_\o\big)
\end{gather*}
for a standard form $(\Ma,\Ha,J,\Pa)$ and a faithful weight $\o$ on $\Ma$. For $\xi\in\Pa$ define $\Ha(\xi):=[\Ma\xi]\cap[\Ma'\xi]$ and $e:=s(\o_\xi)$ where $\o_\xi$ is the semifinite and faithful weight obtained by the GNS-construction. For the commonly defined Tomita operator $S_\xi$ and modular conjugation $J_\xi$ one identifies $\Delta_\xi^{1/2}=JS_\xi$ and $J_\xi=J$. With $\Pa_\xi\{AJA\xi|\;A\in\Ma_e\}$ one obtains the inclusion $\Pa_\xi\subset\Pa$ and moreover, due to self-duality of $\Pa_\xi$ in $\Ha(\xi)$, $\Pa_\xi=\Pa\cap\Ha(\xi)$. \\
Let us assume first $\sigma$-finiteness of $\Ma$ and consider a maximal orthogonal family $\{\xi_i\in\Pa\}_{i\in I}$ then, thanks to Lemma \ref{orthogonalxi}, $s(\o_{\xi_i})\perp s(\o_{\xi_j})$ for $i\neq j$ and thus $I$ is countable. By Lemma \ref{orthogonalxi} maximality of the family leads to faithfulness of $\o_\xi$ and consequently to $\Ha=\Ha(\xi)$ and $\Pa=\Pa_\xi$.\\
If $\sigma$-finiteness of $\Ma$ is not given, then set for each $\sigma$-finite $e\in\text{Proj}(\Ma)$ define $\Ha(e):=e\Ha\cap Je\Ha$. For such an projection one may always find a $\xi\in\Pa$ such that $\Pa_\xi=\Pa\cap\Ha(e)$ and a unique unitary operator 
\begin{gather*}
U_e:\Ha(e)\longrightarrow\pi_\o(e)J_\o\pi_\o(e)\Ha_\o
\end{gather*}
with 
\begin{gather*}
U_eA=\pi_\o(A)U_e\quad\text{and}\\
U_e\Pa_\xi=\Pa_\o\cap U_e\Ha(e),\quad\forall A\in\Ma_e.
\end{gather*}
The family of $\sigma$-finite projections is upward directed with supremum $1$ and thus all operators $U_e$ share the same extension $U$ with the desired properties due to $\Pa=\bigcup\Pa_\xi$. Finally, the uniqueness of $U$ is guaranteed by Lemma \ref{operatorU2}.
\qed
\begin{theo}\label{standardimplementation}
If $(\Ma,\Ha,J,\Pa)$ is a standard form, then there exists an isomorphism between the group $\U$ of all unitary operators fulfilling
\begin{gather*}
U\Ma U^*=\Ma,\quad UJU^*=J,\quad U_g\Pa=\Pa, 
\end{gather*}
and the group of all automorphisms of $\Ma$
\begin{align*}
\pi_\U:\U(\Ma)&\longrightarrow\Aut(\Ma)\\
U&\mapsto\pi_\U^U,\quad\pi_\U^U(A):=UAU^*
\end{align*}
$\pi_\U$ is a homeomorphism of $\U(\Ma)$ equipped with the strong operator topology onto $\Aut(\Ma)$. 
\end{theo}

\proof \cite[IX.1.15.]{Takesaki2:2001}
Obviously, $\pi_\U$ is a homomorphism and due to Proposition \ref{standardisomorphism} also surjective as well injective. The unique map
\begin{gather*}
\Pa\ni\xi\mapsto\o_\xi\in\Ma_*^+
\end{gather*}
is because of 
\begin{gather*}
\|\xi-\eta\|\leq\|\o_\xi-\o_\eta\|\leq\|\xi+\eta\|\|\xi-\eta\|
\end{gather*}
homeomorphic and so is $\pi_\U$.
\qed
The inverse map $\pi_\U^{-1}:\Aut(\Ma)\longrightarrow\U(\Ma)$ is called the standard implementation.


\begin{defi}
Let $\Gr$ be a locally compact group acting on the von Neumann algebra $\Ma$, $\alpha:\Gr\longrightarrow\autMa$, $(\pi,\Ha)$ a normal representation of $\Ma$ and $U:\Gr\longrightarrow\LH$ a unitary representation of $\Gr$ on $\Ha$. Then $\pi$ and $U$ are said to be covariant, if the following relation is satisfied:
\begin{gather}
\pi\circ\alpha^g(A)=U_g\pi(A)U_g^*, \quad\forall A\in\Ma, g\in\Gr.
\end{gather}
\end{defi}

\begin{theo}\label{generalTakesaki}
Let $(\Ma,\Ha,J,\Pa)$ be a standard form, $\Gr$ a locally compact group and $\alpha:\Gr\longrightarrow\autMa$. Define the representations $\pi_\alpha:\Ma\longrightarrow \L\big(L^2(\Gr,\Ha\big))$ and $U:\Gr\longrightarrow\L\big(L^2(\Gr,\Ha)\big)$ of $\Ma$ and $\Gr$, respectively, as follows:
\begin{align}
[\pi_\alpha(A)\xi](g)&:=(\alpha^g)^{-1}(A)\xi(g),\quad\forall A\in\Ma, g\in\Gr\notag,\\
[U_{g_1}\xi](g_2)&:=\xi(g_1^{-1}g_2),\quad\forall \xi\in L^2(\Gr,\Ha), g_1,g_2\in\Gr.
\end{align}
Then $\pi_\alpha$ and $U$ are covariant, i.e. 
\begin{gather}
\pi_\alpha\circ\alpha^g(A)=U_g\pi_\alpha(A)U_g^*, \quad\forall A\in\Ma, g\in\Gr.
\end{gather}
\end{theo}

\proof \cite[X.1.7.]{Takesaki2:2001}
Unitarity holds for $U$ as it is operating as a translation.
\qed

\begin{defi}
The von Neumann algebra generated by $\pi_\alpha(\Ma)$ and $\U(\Gr)$ is the so-called crossed product of $\Ma$ by $\alpha$ and denoted by $\Ma\rtimes_\alpha\Gr$.
\end{defi}
The algebraic structure of crossed product can be shown to be independent of the underlying representation space $\Ha$ and of cocycle perturbations, to be more precise one has for a cocycle $U$ with $(U\xi)(g)=u_{g^{-1}}\xi(g)$:
\begin{gather*}
U(\Ma\rtimes_\alpha\Gr)U^*=\Ma\rtimes_{u^\alpha}\Gr.
\end{gather*}
Although the concept of the crossed product seems to be somehow artificial, there are good reasons for its introduction. First of all, the algebra of a dynamical system $(\Ma,\alpha,\Gr)$ often turns out to be the crossed product of a fixed point subalgebra of $\Ma$ and the dual group of $\Gr$ whereas the action of $\alpha$ on $\Ma$ can be identified with the dual action on $\Ma\rtimes_\alpha\Gr$. Then, crossed products play a crucial role in the classification of factors; the special case $\Gr\equiv\R$ is for this purpose sufficient. These products enables one to construct interesting examples contributing to a better understanding of dynamical systems. For example, if the von Neumann algebra is abelian and the action free and ergodic then the corresponding crossed product is a factor. This result remains also valid for discrete groups. \\
Last but not least, Theorem \ref{generalTakesaki} ensures the existence of an one-to-one correspondence between covariant systems and crossed products. This is the main motivation for the application of crossed products in this thesis, since in our general case where $\alpha$ does not act spatially on $\Ma$, it endows us a spatial action on the corresponding crossed product. Thus we can transfer the analysis from $\Ma$ onto $\Ma\rtimes_\alpha\Gr$, apply the results we already have obtained for the spatial action on $\Ma$  and project them back onto $\Ma$ itself.

\section{Decomposition Theory}

\begin{cor}\label{linearfunctional}
For a non-vanishing element $X$ of a normed space $\X$ there always exists a bounded linear functional $\rho$ on $\X$ such that $\|\rho\|=1$ and $\rho(X)=\|X\|$.
\end{cor}

\proof On the subspace $\X_0$ of $\X$, generated by $X_0\in\X$ only, one can define the bounded linear functional 
\begin{gather*}
\rho_0(cX_0):=c\|X\|,\quad X_0\in\X_0,
\end{gather*}
which satisfies the demanded conditions $\|\rho_0\|=1$ and $\rho(X_0)=\|X_0\|$. Thanks to the Hahn-Banach theorem, this functional can be extended to a still normed functional $\rho$ on the whole algebra $\X$. 

\qed

\begin{theo}
Let $\X$ be a normed space and $\X^*$ its dual space, then the map
\begin{align*}
\X&\longrightarrow\Cs(\X^*,\C)\\
X&\mapsto \widehat{X}\\
\widehat{X}(\rho)&:=\rho(X)
\end{align*}
is an isometrical isomorphism from $\X$ onto the subspace $\widehat{\X}:=\Cs(\X^*,\C)$ of $\X^{**}$.
\end{theo}

\proof Obviously $\widehat{X}$ is a linear functional on the dual space $\X^*$ and the map $\X\ni X\longrightarrow\widehat{X}\in\Cs(\X^*,\C)$ is a linear operator from $\X$ into the double dual space $\X^{**}$. For $\rho\in\X^*$ and $X\in\X$ with
\begin{gather*}
|\widehat{X}(\rho)|=|\rho(X)|\leq\|\rho\|\|X\|
\end{gather*}
we can choose $\rho$ as in the Corollary \ref{linearfunctional}, namely 
\begin{gather*}
|\widehat{X}(\rho)|=\|X\|=\|\rho\|\|X\|.
\end{gather*}
Consequently $\|\widehat{X}\|=\|X\|$ holds and the isomorphism is isometrical.

\qed

\begin{defi}
An normed space $\X$ is called reflexive if $\widehat{\X}=\X^{**}$ holds.
\end{defi}

\begin{cor}\label{reflexive1}
A closed subspace of reflexive space $\X$ is reflexive. 
\end{cor}

\proof
If $\Y$ is a closed subspace of the reflexive space $\X$ and $Y^{**}\in\Y^{**}$, then the mapping 
\begin{gather*}
\X^{**}\ni X^{**}\mapsto Y^{**}(X|_\Y)
\end{gather*}
lies in $\X^{**}$; this conclusion is justified by
\begin{gather*}
|Y^{**}(X^*|_\Y)|\leq\|(X|_\Y)\|\|(X^*|_\Y)\|\leq\|(X|_\Y)\|\|X^*\|.
\end{gather*}
Since $\X$ is a normed, reflexive space one can find an element $X\in\X$ such that
\begin{equation}\label{reflexive2}
X^{*}(X)=Y^{**}(X^*|_{\Y^{**}})
\end{equation}
holds for all $X^{*}\in\X^*$. The assumption $X\notin\Y$ would lead, as $\Y$ is a closed subspace, to the existence of a functional $X^*\in\X^*$ with $X^*(X)=1$ and $X^*|_\Y=0$ and finally to $Y^{**}(X^*|_\Y)=0$, in contradiction to Equation \eqref{reflexive2}. Thus we are allowed to assume $X\in\Y$ and only the relation 
\begin{gather*}
Y^{**}(Y^*)=Y^*(Y)
\end{gather*}
remains to be proved for all $Y^*\in\Y^*$. But let herefore an arbitrary $Y^*\in\Y^*$ be given, then one can extend it, due to the Hahn-Banach theorem, to a functional onto the whole space $\X^*$ such that:
\begin{gather*}
Y^{**}(Y^*)=Y^{**}(X^*|_\Y)\overset{\eqref{reflexive2}}{=}X^*(Y)=Y^*(Y).
\end{gather*}
Thus one has $\Y^{**}=\widehat{\Y}$ and $\Y$ is reflexive.

\qed

\begin{cor}
For a normed space $\X$ the following statements are equivalent:
\begin{itemize}
\item[(i)] $\X$ is reflexive. 
\item[(ii)] $\X^*$ is reflexive.
\item[(iii)] The unit ball $\Aa_1$ of $\X$ is (weakly) $\sigma(\X,\X^*)$-compact.
\end{itemize}
\end{cor}

\proof (i)$\Longrightarrow$(ii) It has to be shown that the mapping
\begin{gather*}
i_{\X^*}:\X^*\longrightarrow\X^{***}
\end{gather*}
is surjective. For a given element $\X^{***}\in\X^{***}$ the mapping 
\begin{align*}
X^*:\X&\longrightarrow\C\\
X&\mapsto X^{***}\big(i_\X(X)\big)
\end{align*}
is linear and continuous. Because $\X$ is reflexive, each $X^{**}\in\X^{**}$ can be written as $X^{**}=i_\X(X)$ and one obtains 
\begin{gather*}
X^{***}(X^{**})=X^{***}\big(i_\X(X)\big)=X^*(X)=\big(i_\X(X)\big)(X^*)=X^{**}(X^*).
\end{gather*}
Consequently, each $X^{***}$ has a pre-image
\begin{gather*}
X^{***}=i_{\X^*}(X^*)
\end{gather*}
and $\X^*$ is reflexive.\\
(ii)$\Longrightarrow$(i) Repeating the argumentation of the first part $\X^{**}$ has to be reflexive and therefore, due to Corollary \ref{reflexive1}, the closed subspaces $i_\X(X)$ and $\X$, too.
(i)$\Longrightarrow$(iii) Thanks to the theorem of Alaoglu the unit ball of the dual space of a normed space is $\sigma(\X^*,\X)$-compact. This means that we are allowed to assume the unit ball $\X_1^{**}$ of the double dual space $\X^{**}$ being $\sigma(\X^*,\X)$-compact. But since we may identify $\X^{**}$ with $\X$ and the mapping
\begin{gather*}
\X\longrightarrow\X^{**}
\end{gather*}
is a homeomorphism, the unit ball $\X_1$ of $\X$ has to be $\sigma(\X,\X^*)$-compact.\\
(iii)$\Longrightarrow$(i) If $\X_1$ is $\sigma(\X,\X^*)$-compact then its image $i(\X_1)$ in $\X^{**}$ must be $\sigma(\X^*,\X)$-compact, this means in particular that it is closed. But this image is always $\sigma(\X^*,\X)$-dense in the unit ball of the double dual and we conclude that $i(\X_1)\equiv\X_1^{**}$, i.e. $\X$ is reflexive.

\qed

\begin{cor}\label{reflexive3}
The double dual space $\X^{**}$ is identical to $\X$ under the $\sigma(\X^*,\X)$-topology.
\end{cor}

\proof This statement is a direct consequence of the last two corollaries.

\qed





\begin{cor}
The states over a C$\*$-algebra $\Aa$ form a convex set denoted by $\Sa_\Aa$ in the dual space $\Aa^*$ of $\Aa$.
\end{cor}

\proof Let $\o_1$ and $\o_2$ be two states over $\Aa$, then a linear combination $\o:=\lambda\o_1+(1-\lambda)\o_2$ with $\lambda\in[0,1]$ is still positive and furthermore one has
\begin{align*}
\|\o_1+\o_2\|&=\lim_\alpha\big(\o_1(E_\alpha^2)+\o_2(E_\alpha^2)\big)\\
&=\lim_\alpha\o_1(E_\alpha^2)+\lim_\alpha\o_2(E_\alpha^2)\\
&=\|\o_1\|+\|\o_2\|.
\end{align*}
Therefore we obtain
\begin{gather*}
\|\o\|=\lambda\|\o_1\|+(1-\lambda)\|\o_2\|=1
\end{gather*} 
and $\o$ is shown to be a state.

\qed

\begin{cor}\label{weaklystarcompact}
Let $\Ba_\Aa$ be the convex set of states $\o$ with $\|\o\|\leq 1$ over the C$\*$-algebra $\Aa$ then $\Ba_\Aa$ is a weakly$\*$-compact subset of $\Aa^*$ whose extremal points are $0$ and the pure states and  $\Ba_\Aa$ is the weak$\*$-closure of the convex envelope of its extremal points.
\end{cor}

\proof The first statement follows from the fact that $\Ba_\Aa$ is a weakly$\*$-closed subset of the weakly$\*$-compact unit ball $\Aa_1^*:=\{\o\in\Aa^*|\;\|\o\|\leq 1\}$ of $\Aa^*$.\\
Concerning the extremal points, the assumption $\o,-\o\in\Ba_\Aa$ leads to $\o(A^*A)=0$ and therefore $\o(A)=0$ for all $A\in\Aa$. Consequently $\o=0$ and $0$ has to be an extremal point. Let us suppose next that the pure state $\o\in\Pa_\Aa$ could be described as a convex combination 
\begin{gather*}
\o=\lambda\o_1+(1-\lambda)\o_2,\qquad\lambda\in(0,1),
\end{gather*} 
of the states $\o_1,\o_2\in\Ba_\Aa$. This would mean $\o\geq\lambda\o_1$ and thus due to purity $\lambda\o_1=\mu\o$ for an appropriate $\mu\in[0,1]$. Since all three states are normed, i.e. 
\begin{gather*}
\|\o\|=\lambda\|\o_1\|+(1-\lambda)\|\o_2\|=1\quad\text{and}\|\o_1\|=\|\o_2\|=1,
\end{gather*} 
it follows that $\lambda=\mu$ and therefore $\o=\o_1$. With the same arguments one concludes $\o=\o_2$ and consequently the pure state $\o$ must be an extremal point.\\ 
Let us now assume the existence of an extremal but not pure state $\o\in\Ba_\Aa$. Accordingly there exists a state $\o_1\neq\o$ and a number $\lambda\in(0,1)$ such that $\o\geq\lambda\o_1$, and we can define a new state $\o_2:=(1-\lambda)^{-1}(\o-\lambda\o_1)$. But this is contradictory to our assumption because $\o$ can be described as a convex combination of $\o_1$ and $\o_2$ and therefore can not be an extremal point.\\
The last assertion follows from the Klein-Milman theorem.
 
\qed

\begin{cor}
For a C$\*$-algebra $\Aa$ the following statements are equivalent:
\begin{itemize}
\item[(i)] The set of states $\Sa_\Aa$ is $\sigma(\Aa^*,\Aa)$-compact;
\item[(ii)] $\Aa$ is unital.  
\end{itemize}
If these conditions are fulfilled then the extremal points of $\Sa_\Aa$ are the pure states $\Pa_\Aa$ and $\Sa_\Aa$ is the weak$\*$-closure of the convex envelope of $\Pa_\Aa$.
\end{cor}

\proof (i)$\Longrightarrow$(ii): We want to show that if $\Aa$ does not contain the unity $\1$ then  the set of states $\Sa_\Aa$ is not weakly$\*$-compact. For this purpose it is sufficient to prove that every arbitrary neighborhood of $0$ includes a state. Since each element of $\Aa$ can be described as a linear combination of four positive elements we can restrict ourselves to the neighborhoods indexed by $A_1,...,A_n\in\Aa_+$. Let us consider an arbitrary element $A=A_1+\cdots+A_n\in\Aa_+$ and a faithful nondegenerate representation $(\Ha,\pi)$ of $\Aa$. Because $A$ has no inverse in $\Aa$ and thus no one in the unitalization $\tilde{\Aa}:=\Aa+\C\1$ of $\Aa$, $\pi(A)$ cannot be invertible in $\LH$. Thus there has to exist a unit vector $\psi\in\Ha$ with $\o_\psi(A)=(\psi,A,\psi)<\epsilon$ which holds the demanded property.\\ 
(ii)$\Longrightarrow$(i): For a unital C$\*$-algebra $\Aa$ the set of all states $\Sa_\Aa$ can be written as the intersection of $\Ba_\Aa$ with the hyperplane $\o(\1)=1$ and we can apply Corollary \ref{weaklystarcompact}.\\

\qed

\begin{theo}[Central decomposition]
Let $\rho$ be a state on a separable C$\*$-algebra $\Aa$ and $\mu_\rho$ its central measure on $\widehat{\Aa}$, then there exists an essentially unique Borel function $\widehat{\Aa}\ni\hat{A}\longrightarrow\rho_{\hat{A}}\in\Fa_\Aa$ with its relative Borel structure such that 
\begin{gather*} 
\rho(A)=\int\rho_{\hat{A}}(A)d\mu_\rho(\hat{A})
\end{gather*} 
for all $A$ in the enveloping Borel$\,\*$-algebra $\Ba_\Aa$ of $\Aa$. 
\end{theo}

\proof
Confer \cite[4.8.7.]{Pedersen:1970}.

\qed


\begin{lem}
Let $\Ebf$ be a locally compact topological vector space and $\Kbf$ a compact convex subset. Then for each measure $\mu\in\Mbf_1^+(\Kbf)$ there exists a unique point $y$, the so-called barycenter of $\mu$,  such that for $\varphi\in\Ebf^*$ one has
\begin{gather*}
\varphi(y)=\int_\Kbf\varphi(x)d\mu(x).
\end{gather*}
\end{lem}

\proof First, for a measure $\mu$ with finite support,
\begin{gather*}
\mu=\sum_{i=1}^\infty\lambda_i\delta_{x_i},
\end{gather*}
where $\lambda_i\geq 0$ and $\sum_{i=1}^\infty\lambda_i=1$, the barycenter is given by
\begin{gather*}
y=\sum_{i=1}^\infty\lambda_i\delta_{x_i}
\end{gather*}
and is included in $\Kbf$ as $\Kbf$ is a convex set. In general, a measure $\mu\in\Mbf_1^+(\Kbf)$ can be written as a $\sigma(\Ebf^*,\Ebf)$-limit of a sequence of measures $(\mu_n)_{n\in\N}$, $\mu_n\in\Mbf_1^+(\Kbf)$, having finite support. Since the barycenter for each $\mu_n$ lies in $\Kbf$ and $\Kbf$ is compact, there exists a convergent sequence of barycenters $(y_n)_{n\in\N}$ such that
\begin{align*}
\varphi(y):&=\lim_{n\rightarrow\infty}\varphi(y_n)\\
&=\lim_{n\rightarrow\infty}\int_\Kbf\varphi(x)d\mu_n(x)\\
&=\int_\Kbf\varphi(x)d\mu(x),
\end{align*}
for all $\varphi\in\Ebf^*$. The uniqueness of $y$ is assured as $\varphi$ separates the points of $\Kbf$ by the Hahn-Banach Theorem.

\qed

\begin{defi}
Let $\Kbf$ is a compact convex set, then for each real-valued continuous function $f\in\Cs_\R(\Kbf)$ one introduces the boundary set 
\begin{gather*}
B_f:=\big\{x\in\Kbf|\;\hat{f}(x)=f(x)\big\},
\end{gather*}
where $\hat{f}$ is the upper envelope of $f$. A measure on $\Kbf$ is said to be a boundary measure if the absolute value of $\mu$ vanishes for all $f\in\Cs_\R(K)$ on the complement of $B_f$, i.e. $|\mu|(B_f^c)=0$.
\end{defi}

\begin{lem}
Let $\Kbf$ is a compact convex set, then one has 
\begin{gather*}
\partial_e(\Kbf)=\bigcap_{n=1}^\infty\big\{B_f|\;f\in\Cs_\R(\Kbf\big\}.
\end{gather*}
\end{lem}

\proof Let $x$ be an extremal point of $\Kbf$, then the set of all measures representing $x$ consists only of the point mass measure, $\Mbf_x^+(\Kbf)=\{\delta_x\}$. Thus we obtain $\hat{f}(x)=f(x)$ and with this $x\in\bigcap_{n=1}^\infty\big\{B_f|\;f\in\Cs_\R(\Kbf\big\}$.\\
Let us assume now $x\in\bigcap_{n=1}^\infty\big\{B_f|\;f\in\Cs_\R(\Kbf\big\}$, then one gets, due to the equation $-\check{f}=(-f)\hat{}$, the identities $f(x)=\hat{f}(x)=\check{f}$, where $\check{f}$ denotes the lower limit of $f$. Because of $\delta_x(g)<\nu(g)$ for all $\nu\in\Mbf_x^+(\Kbf)$ and all elements $g$ of $\A(\Kbf)$, the set of real-valued continuous affine functions on $\Kbf$, the following chain of inequalities
\begin{gather*}
f(x)=\check{f}(x)\leq\nu(\check{f})\leq\nu(f)\leq\nu(\hat{f})\leq\hat{f}(x)=f(x)
\end{gather*}
holds, and we can conclude $\nu\equiv\delta_x$. Consequently, $x$ has to be an extreme point, since it is not describable as a non-trivial convex combination of other points in $\Kbf$.

\qed
Every state $\varphi\in\Sa_\Aa$ on a $C^*$-algebra $\Aa$ may be represented by a boundary measure $\mu$ via
\begin{gather*}
\varphi(A)=\int_{\Sa_\Aa}\psi(A)d\mu(\psi),
\end{gather*}
where $A\in\Aa$. In the case of $\Aa$ being a separable $C^*$-algebra, $\mu$ has its support on the extreme boundary of the state space, $\supp\mu\subseteq\partial_e(\Sa_\Aa)\equiv\Pa_\Aa$, so that we may write:
\begin{gather*}
\varphi(A)=\int_{\Pa_\Aa}\psi(A)d\mu(\psi),
\end{gather*}

\begin{prop}\label{boundary-extreme}
A measure $\mu\in\Mbf^+(\Kbf)$ is extreme, i.e. $\mu(f)\geq\nu(f)$ for all $\nu\in\Mbf^+$ and all $f\in\Cs_\R(\Kbf)$, if and only if it is a boundary measure.
\end{prop}

\proof For a maximal measure $\mu\in\Mbf^+(\Kbf)$ one has $\mu(f)=\mu(\hat{f})$ for all $f\in\Cs_\R(\Kbf)$, and therefore $\mu$ has to  be a boundary measure since $\mu\big(B_f^c\big)=0$.\\
Let us now assume that $\mu$ is a boundary measure. In this case $\mu(f)=\mu(\hat{f})$ and therefore $\mu(f)=\mu(\check{f})$ follow for all $f\in\Cs_\R(\Kbf)$. If we suppose the existence of a measure $\nu\in\Mbf^+(\Kbf)$ with $\nu(f)>\mu(f)$ for all $f\in\Cs_\R(\Kbf)$, then one obtains:
\begin{gather*}
\mu(f)=\mu(\check{f})\leq\nu(\check{f})=\nu(f)\leq\nu(\hat{f})\leq\mu(\hat{f})=\mu(f).
\end{gather*}
Consequently, $\mu$ and $\nu$ must be the same measure.

\qed

\begin{defi}
Let $Aa$ be a unital $C^*$-algebra, then a measure $\mu\in\Mbf_1^+(\Kbf)$ is said to be orthogonal if the functionals
\begin{gather*}
\varphi_E(A):=\int_E\psi(A)d\mu(\psi)\quad\text{and}\\
\varphi_{E^c}(A):=\int_{E^c}\psi(A)d\mu(\psi)
\end{gather*}
are orthogonal, i.e. $0\leq\psi'\leq\varphi_E$ and $0\leq\psi'\leq\varphi_{E^c}$ imply $\psi'=0$.
\end{defi}

\begin{prop}\label{orthogonal-boundary}
If the measure $\mu\in\Mbf_1^+(\Sa_\Aa)$ with barycenter $\varphi$ is orthogonal then it is a boundary measure of $\Mbf_\varphi^+(\Sa_\Aa)$.
\end{prop}

\proof Let $\mu$ be an orthogonal measure, then the map
\begin{align}\label{Thetamu}
\Theta_\mu:L^\infty(\Sa_\Aa,\mu)&\longrightarrow\pi_\varphi(\Aa)'\\
f&\mapsto\Theta_\mu(f)\notag\\
\big(\Theta_\mu(f)\pi_\varphi(A)\xi_\varphi,\xi_\varphi\big)&:=\int_{\Sa_\Aa}f(\psi)\psi(A)d\mu(\psi),\notag\\
\Theta(\1)&:=1,\notag
\end{align}
where $A$ is an element of $Aa$, can be shown to be injective. If $f$ is orthogonal to the set $\A_\C(\Sa_\Aa):=\{\tilde{A}|\;A\in\Aa\}$, where $\tilde{A}(\psi):=\psi(A)$, then we obtain for all elements $A,B\in\Aa$:
\begin{gather*}
\big(\Theta_\mu(f)\pi_\varphi(A)\xi_\varphi,\pi_\varphi(B)\xi_\varphi\big):=\int_{\Sa_\Aa}f(\psi)\psi(B^*A)d\mu(\psi)=0.
\end{gather*}
This leads to $\Theta_\mu(f)=0$ and therefore to $f=0$. Because of the duality relation $L^1(\Sa_\Aa,\mu)^*=L^\infty(\Sa_\Aa,\mu)$ the set $\A_\C(\Sa_\Aa)$ has to be dense in $L^1(\Sa_\Aa,\mu)$. Let us assume that our measure is not extreme, i.e. it can be written as a convex combination $\mu=\frac{1}{2}(\mu_1+\mu_2)$ with $\mu_1,\mu_2\in\Mbf_\varphi^+(\Sa_\Aa)$. Obviously $0\leq\mu_i\leq 2\varphi$, $i=1,2$, and thus there exist functions $h_1,h_2\in L^\infty(\Sa_\Aa,\mu)$ with $0\leq h_i\leq 2$, $i=1,2$, and satisfying 
\begin{gather*}
\int_{\Sa_\Aa}f(\psi)d\mu_i(\psi)=\int_{\Sa_\Aa}f(\psi)h_i(\psi)d\mu(\psi)
\end{gather*}
for all $f\in L^\infty(\Sa_\Aa,\mu)$. Additionally, both $\mu_i$ are representing $\varphi$, $\mu_i\in\Mbf_\varphi^+(\Sa_\Aa)$, resulting in 
\begin{gather*}
\varphi(A)=\int_{\Sa_\Aa}\tilde{A}(\psi)d\mu_i(\psi)=\int_{\Sa_\Aa}\tilde{A}(\psi)h_i(\psi)d\mu(\psi)
\end{gather*}
for all $A\in\Aa$. Hence $1-h_i$ has to be orthogonal to $\A_\C$, and since $\A_\C$ is dense in $L^1(\Sa_\Aa,\mu)$, we obtain  the identity $h_i=1$, $i=1,2$, in $L^\infty(\Sa_\Aa,\mu)$. Finally, we conclude $\mu=\mu_1=\mu_2$ and $\mu$ must therefore be an extremal measure. Due to Proposition \ref{boundary-extreme} $\mu$ is also a boundary measure.

\qed
We will need in the following the range of the isomorphism $\Theta_\mu$, defined in \eqref{Thetamu}, which can be shown easily to be an abelian von Neumann subalgebra of $\pi_\varphi(\Aa)'$. We denote it by $\Ca_\mu$. The next proposition justifies its importance since all informations of the associated orthogonal measure $\mu$ is encoded in $\Ca_\mu$.

\begin{prop}\label{abelian-subalgebra}
For each abelian von Neumann subalgebra $\Ca$ of $\pi_\varphi(\Aa)'$, where $\varphi\in\Sa_\Aa$, there exists a unique orthogonal measure $\mu\in\Mbf_\varphi^+(\Sa_\Aa)$ such that $\Ca=\Ca_\mu$.
\end{prop}
 
\proof Because the vector $\xi_\varphi$ is cyclic for $\pi_\varphi(\Aa)$ and $\Ca$ lies in $\pi_\varphi(\Aa)'$, $\xi_\varphi$ has to be separating for $\Ca$ and the normal state $\bar{\varphi}(A):=(A\xi_\varphi,\xi_\varphi)$ must be faithful. If $E$ is the projection onto $[\Ca\xi_\varphi]$ then $\Ca_E$ is maximal abelian, because $\xi_\varphi$ is also cyclic with respect to $\Ca_E$ and the map 
\begin{gather*}
\Ca\ni A\mapsto A_E\in\Ca_E
\end{gather*}
is isomorphic. Because of the identity $(\Ca_E)'=(\Ca')_E$ and the inclusion $E\pi_\varphi(A)E\subset\Ca_E$ we may define the positive mapping
\begin{align*}    
\Theta:\Aa&\longrightarrow\Ca\\
A&\mapsto\Theta(A)\\
\Theta(A)\xi_\varphi&:=\pi_\varphi(A)\xi_\varphi,
\end{align*}
with $\quad\Theta(\1)=1$. One may establish with the help of the transposes of $\Theta$ a continuous map $\Theta^t:\xi\longrightarrow\Sa_\Aa$, where $\xi$ is the spectrum of $\Ca$, and also a map $\Theta^{tt}:\Ca(\Sa_\Aa)\longrightarrow\Ca$, an extension of the original map $\Theta$. We obtain now by defining $\mu:=\bar{\varphi}\Theta^{**}\in\Mbf_1^+(\Sa_\Aa)$ a representing measure for $\varphi$, since one has:
\begin{align*}
\varphi(A)&=\big(\pi_\varphi(A)\xi_\varphi,\xi_\varphi\big)=\big(\Theta(A)\xi_\varphi,\xi_\varphi\big)\\
&=\bar{\varphi}\circ\Theta(A)=\bar{\varphi}\circ\Theta(\tilde{A})\\
&=\int_{\Sa_\Aa}\tilde{A}(\psi)d\mu(\psi)\\
&=\int_{\Sa_\Aa}\psi(A)d\mu(\psi).
\end{align*}
Because $\Theta:\Ca(\Sa_\Aa)\longrightarrow\Ca$ is a homomorphism transporting $\bar{\varphi}$ into $\mu$, $\Theta^{**}$ is extended to a normal isomorphism 
\begin{gather*}
\Theta_\mu:L^\infty(\Sa_\Aa,\mu)\longrightarrow\Ca.
\end{gather*}
The algebra $\Ca_0$ generated by $E\pi_\varphi(\Aa)E$ on $E\Ha_\varphi=[\Ca\xi_\varphi]$ is a von Neumann subalgebra of $\Ca_E$ which itself is generated by $\Theta(\Aa)$. Because of 
\begin{gather*}
[\Ca_0\xi_\varphi]\supset[E\pi_\varphi(\Aa)E\xi_\varphi]=[E\pi_\varphi(\Aa)\xi_\varphi]=E\Ha_\varphi
\end{gather*}
the vector $\xi_\varphi$ must be cyclic for $\Ca_0$ and, consequently, $\Ca_0$ has to be a maximal abelian von Neumann subalgebra, which means $\Ca_0=\Ca$. Therefore $E\pi_\varphi(\Aa)E$ generates $\Ca_E$, $\Theta(\Aa)$ generates $\Ca$ and accordingly we obtain $\Theta^{**}\big(L^\infty(\mu)\big)=\Ca$. Now, we can identify $\Theta$ with $\Theta_\mu$ defined in \eqref{Thetamu} due to the following relation with $f\in L^\infty(\mu)$:
\begin{align*}
\big(\pi_\varphi(A)\Theta(f)\xi_\varphi,\xi_\varphi\big)&=\big(E\pi_\varphi(A)E\Theta(f)\xi_\varphi,\xi_\varphi\big)\\
&=\big(\Theta(A)\Theta(f)\xi_\varphi,\xi_\varphi\big)\\
&=\big(\Theta^{**}(\tilde{A}f)\xi_\varphi,\xi_\varphi\big)\\
&=\int_{\Sa_\Aa}\tilde{A}(\psi)f(\psi)d\mu(\psi)\\
&=\int_{\Sa_\Aa}\psi(A)f(\psi)d\mu(\psi).
\end{align*}
The last thing to show is the uniqueness of the orthogonal measure. Let us consider two orthogonal measures $\mu,\nu\in\Mbf_1^+(\Sa_\Aa)$ with their associated abelian von Neumann algebras $\Ca_\mu$ and $\Ca_\nu$. We may assume $\Ca_\mu=\Ca_\nu$ and $[\Ca_\mu\xi_\varphi]=[\Ca_\nu\xi_\varphi]$. As shown above, we have for all $A\in\Aa$
\begin{gather*}
\Theta_\mu(\tilde{A})\xi_\varphi=E\pi_\varphi(A)\xi_\varphi=\Theta_\nu(\tilde{A})\xi_\varphi,\\
\Theta_\mu(\tilde{A})E=E\pi_\varphi(A)\xi_\varphi E=\Theta_\nu(\tilde{A})E,
\end{gather*}
and herewith $\Theta_\mu(\tilde{A})=\Theta_\nu(\tilde{A})$ for all $\A\in\Aa$. Thus we have proved that $\Theta_\mu$ and $\Theta_\nu$ coincide on the $C^*$-algebra generated by the set $\big\{\tilde{A}|\;A\in\Aa\big\}=\A_\C(\Sa_\Aa)$, which is, due to the Stone-Weierstrass Theorem, equivalent to $\Cs(\Sa_\Aa)$. Finally, we conclude for all $f\in\Cs(\Sa_\Aa)$:
\begin{gather*}
\mu(f)=\big(\Theta_\mu(f)\xi_\varphi,\xi_\varphi\big)=\big(\Theta_\nu(f)\xi_\varphi,\xi_\varphi\big)=\nu(f).
\end{gather*}

\qed

\begin{defi}
A measure $\mu\in\Mbf_1^+(\Kbf)$ is said to be pseudo-concentrated on a subset $\mathbf{K'}$ if $\mu(\mathbf{L})=0$ for all Baire sets $\mathbf{L}\subset\Kbf$ disjoint from $\mathbf{K'}$.
\end{defi}

\begin{theo}\label{pseudo-concentrated}
If $\Aa$ is a unital $C^*$-algebra and $\mu$ is an orthogonal representing measure of $\varphi\in\Sa_\Aa$ whose associated abelian von Neumann algebra $\Ca_\mu$ is maximal abelian in $\pi_\varphi(\Aa)'$, then $\mu$ is pseudo-concentrated on the pure state space $\Pa_\Aa$. In the case of a separable $C^*$-algebra $\Aa$ the measure $\mu$ is concentrated on $\Pa_\Aa$.
\end{theo}

\proof  First of all, $\Theta_\mu\big(\Cs(\Sa_\Aa)\big)\subset\Ca_\mu$ commute with $\pi_\varphi(\Aa)$, and therefore there exists a representation of the injective tensor product $\Ca_\Aa:=\Cs(\Sa_\Aa)\otimes_\text{min}\Aa$ such that 
\begin{gather*}
\pi(f\otimes A)=\Theta_\mu(F)\pi(A),\qquad\forall f\in\Cs(\Sa_\Aa),
\end{gather*}
and $A\in\Aa$. Because the vector $\xi_\varphi$ is also cyclic for $\Ca_\Aa$, the representation $\pi$ of $\Ca_\Aa$ is cyclic for with respect to $\xi_\varphi$ and unitarily equivalent to $\pi_{\tilde{\varphi}}$, where 
\begin{gather*}
\tilde{\varphi}(X):=\big(\pi(X)\xi_\varphi,\xi_\varphi\big),\qquad X\in\Ca_\Aa,
\end{gather*}
 is a state on $\Ca_\Aa$. For this reason, we may identify the two representations $\pi$ and $\pi_{\tilde{\varphi}}$, and show the multiplicity freedom of $\pi_{\tilde{\varphi}}$,
\begin{align*}
\pi_{\tilde{\varphi}}(\Ca_\Aa)'&=\big(\pi_\varphi(\Aa)\cup\Theta_\mu(\Ca(\Sa_\Aa))\big)'\\
&=\pi_\varphi(\Aa)'\cap\Theta_\mu\big(\Ca(\Sa_\Aa)\big)'\\
&=\pi_\varphi(\Aa)'\cap\Ca_\mu'\\
&=\Ca_\mu,
\end{align*}
where we have used the maximal commutativity of $\Ca_\mu$ in $\pi_\varphi(\Aa)'$. Therefore there exists  for $\pi_{\tilde{\varphi}}$ the unique maximal representing measure $\tilde{\mu}$ which is pseudo-concentrated on the extremal points of the state space of $\Ca_\Aa$, i.e. the pure states $\Pa_{\Ca_\Aa}$. But we also have $\Pa_{\Ca_\Aa}=\Sa_\Aa\times\Pa_\Aa$.\\
Now we want to investigate the maximal measure $\tilde{\mu}$ and its uniqueness. Let us consider the the homomorphism
\begin{align*}
\Phi:\Sa_\Aa&\longrightarrow\Sa_{\Ca_\Aa}\\
\psi&\mapsto\Phi(\psi):=\psi\otimes\psi,
\end{align*}
and we will show $\tilde{\mu}=\Phi(\mu)$. If we set $g:=h\circ\Phi$, where $h\in L^\infty\big(\Sa_{\Ca_\Aa},\Phi(\mu)\big)$, then we obtain for every $f\in\Cs(\Sa_\Aa)$ and $A\in\Aa$:
\begin{align*}
\int_{\Sa_{\Ca_\Aa}}h(\tilde{\psi})(f\otimes A)\,\tilde{}\;(\tilde{\psi})d\Phi(\mu)_{\tilde{\psi}}&=\int_{\Sa_\Aa}h(\psi\otimes\psi)(f\otimes A)\,\tilde{}\;(\psi\otimes\psi)d\nu_{\psi}\\
&=\int_{\Sa_\Aa}g(\psi)f(\psi)\tilde{A}(\psi)d\nu_{\psi}\\
&=\big(\Theta_\mu(gf)\pi_\varphi(A)\xi_\varphi,\xi_\varphi\big)\\
&=\big(\Theta_\mu(g)\pi_{\tilde{\varphi}}(f\otimes A)\xi_\varphi,\xi_\varphi\big).
\end{align*}
A linearisation leads to
\begin{gather*}
\int_{\Sa_{\Ca_\Aa}}h(\tilde{\psi})\tilde{X}(\tilde{\psi})d\Phi(\mu)_{\tilde{\psi}}=\big(\Theta_\mu(g)\pi_{\tilde{\varphi}}(X)\xi_\varphi,\xi_\varphi\big)
\end{gather*}
for $X\in\Ca_\Aa$, which shows that $\Phi(\mu)$ is a representing measure of $\tilde{\varphi}$ with $\Theta_{\phi(\mu)}=\Theta_\mu\circ\Phi^*$, where $\Phi^*(h):=h\circ\Phi$ for all $h\in L^\infty\big(\Sa_{\Ca_\Aa},\Phi(\mu)\big)$. Our assumption $\tilde{\mu}=\Phi(\mu)$ is valid, because of the uniqueness of $\tilde{\mu}$ and the isomorphisms 
\begin{gather*}
\Phi^*:L^\infty\big(\Phi(\mu)\big)\longrightarrow L^\infty(\mu)\quad\text{and}\\
\Theta_{\phi(\mu)}:L^\infty\big(\Phi(\mu)\big)\longrightarrow\Ca_\mu.
\end{gather*}
It remains to prove that $\mu$ is pseudo-concentrated on the pure states $\Pa_\Aa$. For this purpose we define the embedding $i:\Aa\longrightarrow\Ca_\Aa$ and with this map the continuous restriction $i^*:\Sa_{\Ca_\Aa}\longrightarrow\Sa_\Aa$ for which we have $i^*\circ\Phi(\psi)=\psi$. For a $G_\delta$-compact subset $E\subset\Sa_\Aa$ which is disjoint from $\Pa_\Aa$, the set $(i^*)^{-1}(E)$ is a $G_\delta$-compact subset of $\Pa_{\Ca_\Aa}$. As shown above, we obtain
\begin{gather*}
i^*\big(\Pa_{\Ca_\Aa}\big)=i^*\big(\Sa_\Aa\times\Pa_\Aa\big)=\Pa_\Aa,
\end{gather*}
leading to $(i^*)^{-1}\big(E\cap\Pa_{\Ca_\Aa}\big)=\emptyset$. But $\tilde{\mu}$ is a boundary measure on $\Pa_{\Ca_\Aa}$ so that we get $\mu\big((i^*)^{-1}(E)\big)=0$ and herewith:
\begin{align*}
\mu(E)&=\mu\big(i^*\circ\Phi)^{-1}(E)\big)\\
&=\mu\big(\Phi^{-1}\big((i^*)^{-1}(E)\big)\big)\\
&=\tilde{\mu}\big((i^*)^{-1}(E)\big)\\
&=0.
\end{align*}
Thus we have shown that $\mu$ is pseudo-concentrated on the pure states $\Pa_\Aa$ of $\Aa$.
 
\qed



\chapter{Non-Commutative Martingales}

\begin{flushright}
 \emph{Zwei, denen ich auf der Promenade begegnete, \\stellten sich mir gleich zweimal vor. Erst der Herr a \\und dann der Herr b und dann der Herr b und drauf\\ der Herr a, und dann fragten sie mit s\"uffisanter Miene,\\ ob das nicht gleich sei...}\\
\vspace{0,5cm}
Wilhelm Busch\\
Eduards Traum
\end{flushright}
\vspace{0,5cm} 
The purpose of this chapter is to give a straightforward introduction to the transition from classical probability theory to its non-abelian counterpart,  where we restrict ourselves to the very basic notions and concepts which are crucial for our investigation of non-commutative martingales. The first part on commutative probability theory is excerpted mainly from \cite{Bauer:1995}, \cite{Doob:1953} whereas our non-commutative formulation is based on \cite{Coja:2005}.\\ \\
In classical probability theory one considers the underlying triple $(\Omega,\SA,P)$ comprising a set $\Omega$, a $\sigma$-algebra $\SA$ of subsets of $\Omega$ and a probability measure $P:\SA\longrightarrow\R$. An $\SA-\SA'$-measurable mapping 
\begin{gather*}
X:\Omega\longrightarrow\Omega',
\end{gather*}
where $(\Omega',\SA')$ is a measurable space, is called a random variable with values in $\Omega'$ or a $(\Omega',\SA')$-random variable, where we suppose $P$ to be fixed.\\
A family of events $(A_i)_{i\in\TT}$, where $A_i\in\SA$ and $\TT$ is an arbitrary discrete or continuous index set, is said to be (stochastically) independent with respect to $P$ if for all non-void subsets $\{i_1,...,i_n\}$ of distinct elements of $\TT$ one has
\begin{gather*}
P\big(A_{i_1}\cap...\cap A_{i_n}\big)=P\big(A_{i_1}\big)\cdots P\big(A_{i_n}\big).
\end{gather*} 
A family $(\SA_i)_{i\in\TT}$ of subsets of the $\sigma$-algebra is called independent if all possible combinations of events $A_{i_k}\in\SA_{i_k}$, $k=1,...,n$, are independent.
\begin{defi}
A family $(X_t)_{t\in\TT}$ of random variables on a common probability space $(\Omega,\SA,P)$ is said to be independent if the corresponding generated $\sigma$-algebras form an independent family $\big(\sigma(X_i)\big)_{i\in\TT}$.
\end{defi}
A stochastic process is a quadruple $\big(\Omega,\SA,P,(X_t)_{t\in\TT}\big)$ consisting of a probability space and a family $(X_t)_{t\in\TT}$ of random variables on a common probability space. Let us suppose $\TT$ to be an ordered set, usually it is assumed to be $\N$ or $\R$, then a filtration $(\SF_t)_{t\in\TT}$ is an increasing family of $\sigma$-algebras of $\Omega$, i.e.
\begin{gather*}
s\leq t\quad\Longrightarrow\quad\SF_s\subset\SF_t,\qquad s,t\in\TT.
\end{gather*}
In case of $\SF_t$ being a $\sigma$-subalgebra of $\SA$ for all $t\in\TT$, $(\SF_t)_{t\in\TT}$ is called a filtration in $\SA$. The family $(X_t)_{t\in\TT}$ is said to be adapted to the filtration $(\SF_t)_{t\in\TT}$ if $X_t$ is $\SF_t$-measurable for all $t\in\TT$. For our purposes it suffices to consider real valued random variables, i.e. $(\R,\mathbb{B})$-random variables, which transform $P$ to a probability measure $P_X:=X(P)$ on $\R$.
\begin{defi}
Let $(\Omega,\SA,P)$ be a probability space and $(X_t)_{t\in\TT}$ a family of integrable real random variables adapted to the filtration $(\SF_t)_{t\in\TT}$ in $\SA$. Then $(X_t)_{t\in\TT}$ is said to be a supermartingale with respect to $(\SF_t)_{t\in\TT}$ if for all $s,t\in\TT$ with $s\leq t$ one of the following equivalent conditions holds:
\begin{itemize}
\item[(i)] $E(X_t|\SF_s)\leq X_s$ $P$-almost surely;
\item[(ii)] $\int_C X_tdP\leq\int_C X_sdP$ for all $C\in\SF_s$.
\end{itemize}
$(X_t)_{t\in\TT}$ is called a submartingale with respect to $(\SF_t)_{t\in\TT}$ if $(-X_t)_{t\in\TT}$ is a supermartingale. If $(X_t)_{t\in\TT}$ is both, then it is called a martingale.
\end{defi}
\begin{examp}\emph{
Let us consider a game between two players A and B described by an independent sequence $(X_t)_{t\in\N}$ of random variables with values $\pm1$ and expectation $E(X_t)=2p-1$ for all $t\in\N$ and $0\leq p\leq1$, where the outcomes $+1$ and $-1$ are interpreted as win and loss for the player A, respectively. Player A owns an initial capital $S_0$ and chooses before the start of the game a sequence of functions 
\begin{gather*}
b_t:\{-1,+1\}^t\longrightarrow\R_+
\end{gather*}
representing his bet in the $t$-th round, namely $b_t(X_1,...,X_t)$ with $0\leq b_0\leq S_0$, determined by his fortune in the last $t$ rounds, constitutes his wager for the $(t+1)$-th round. We suppose that contrary to the first player player B does never bet. Then, the cumulative gains of player player A 
\begin{gather*}
S_{t+1}=S_t+b_t(X_1,...,X_t)\cdot X_{t+1}
\end{gather*}
form a sequence of integrable random variable $(S_t)_{t\in\N}$. If we set 
\begin{gather*}
\SF_t:=\SA(X_1,...,X_t),
\end{gather*}
then $(\SF_t)_{t\in\N}$ is obviously a filtration in $\SA$. Moreover, $S_n$ is adapted to $(\SF_t)_{t\in\N}$, i.e. is $\SF_t$-measurable for all $n\in\N$, and one gets almost surely for all $t\in\N$:
\begin{align*}
E(S_{t+1}|X_1,...,X_t)&=S_t+b_t(X_1,...,X_t)\cdot E(X_{t+1}|X_1,...,X_t)\\
&=S_t+b_t(X_1,...,X_t)\cdot E(X_{t+1})\\
&=S_t+(2p-1)b_t(X_1,...,X_t).
\end{align*}
Thus, $(S_t)_{t\in\N}$ constitutes a supermartingale for $p\leq\frac{1}{2}$, a submartingale for 
$p\geq\frac{1}{2}$ and a martingale for $p=\frac{1}{2}$. One may refer only to the last case as a fair game.} 
\end{examp}
In this thesis, we want to exemplify the transfer of statements on convergence properties of martingales to non-commutative probability with the following theorem. Its non-commutative analogue will be dealt with in Chapter five.
\begin{theo}\label{convergencemartingale}
Let $(X_t)_{t\in\TT}$ be a martingale with respect to the filtration $(\SF_t)_{t\in\TT}$ satisfying
\begin{gather*}
E\big(|X_1|^2\big)\leq E\big(|X_2|^2\big)\leq E\big(|X_3|^2\big)\leq...
\end{gather*}
and
\begin{gather*}
\lim_{i\rightarrow\infty}E\big(|X_i|^2\big)=:c<\infty,
\end{gather*}
then
\begin{itemize}
\item[(i)] $X_\infty:=\lim_{i\rightarrow\infty}X_i$ exists with probability 1 and 
\item[(ii)] the family $(X_1,X_2,...,X_\infty)$ constitutes a martingale.
\end{itemize}
\end{theo}

\proof
\cite[Theorem VII.4.1.]{Doob:1953} 
\begin{itemize}
\item[(i)] Since the sequence of expectations $\Big(E\big(|X_t|^2\big)\Big)$ is increasing and bounded, the existence of its limit $\lim_{i\rightarrow\infty}E\big(|X_i|^2\big)$ and therefore of $\lim_{i\rightarrow\infty}E\big(|X_i|\big)$ is ensured. The existence of $X_\infty$ then follows from Doob's upcrossing inequality.
\item[(ii)] Moreover, boundedness implies also uniform integrability for all $X_t$, $t\in\TT$. The martingale property of $(X_t)_{t\in\TT}$ means for $s<t$ and $\SF\subset\SF_s$
\begin{gather*}
\int_\SF X_t\,dP=\int_\SF E(X_t|\SF_s)dP=\int_\SF X_sdP.
\end{gather*}
Finally, uniform integrability allows to take the limit $t\rightarrow\infty$ under the integral of and one obtains
\begin{gather*}
\int_\SF X_\infty\, dP=\int_\SF X_sdP,\qquad\forall s\in\TT.
\end{gather*}
i.e. the martingale property of the family $(X_1,X_2,...,X_\infty)$.
\end{itemize}
\qed
\vspace{1cm}
In the functional analysis framework, the space of probability measures on $\R$ may be identified with the space of all states on $\Cnn(\R)$, the space of all continuous real-valued functions vanishing at infinity, and therefore the $X$ may be rewritten as:
\begin{align}\label{randomvariable}
X:\Cnn(\R)&\longrightarrow L^\infty(\Omega,P)\\
f&\mapsto f\circ X.\notag
\end{align}
or
\begin{align*}
X^t:L^1(\Omega,P)=L^\infty(\Omega,P)_*&\longrightarrow\Cnn(\R)^*\\
P&\mapsto P^X,
\end{align*}
Hence, the random variable $X$ constitutes a $\*$-homomorphism from the abelian $C^*$-algebra $\Cnn(\R)$ to the abelian von Neumann algebra $L^\infty(\Omega,P)$.\\
The approach of \cite{Coja:2005} for the non-commutative formulation of probability theory is to substitute in two steps both algebras by non-abelian ones. They start by allowing the right hand-side of \ref{randomvariable} to be a general (non-abelian) von Neumann algebra $\Ma$ acting on a Hilbert space $\Ha$, i.e.
\begin{align*}
X:\Cnn(\R)&\longrightarrow(\Ma,\Ha)\equiv\Ma.
\end{align*}
The probability measure $P$ is represented by a normal state $\varphi:\Ma\longrightarrow\C$ which is transformed by $X^t$ to a state on $\Cnn(\R)$
\begin{align*}
X^t:(\Ma,\Ha)_*&\longrightarrow\Cnn(\R)^*\\
\varphi&\mapsto X^t(\varphi),
\end{align*}
where we have identified $\Ma_*$ with $\Ma^*$. Since $X^t(\varphi)$ determines uniquely a probability measure on $\Cnn(\R)$ and the spectrum of $\Cnn(\R)$ is $\R$ itself, classical probability theory can be applied.\\
In a second step, the left hand-side of \ref{randomvariable} may be generalised to an arbitrary (non-abelian) $C^*$-algebra $\Aa$,
\begin{align*}
X:\Aa&\longrightarrow\Ma,\\
X^t:\Sa(\Ma)&\longrightarrow\Sa(\Aa),
\end{align*}
and $X$ becomes a representation of $\Aa$ on the $\Ma$-valued Hilbert space $\Ha$. Thus, a probability space in the non-commutative sense is a pair $(\Ma,\varphi)$.

\begin{defi}\label{noncommcondexp}
Let $\varphi$ be a faithful (semi-finite normal) weight on a von Neumann algebra $\Ma$ and $\Na$ a von Neumann subalgebra of $\Ma$ such that the restriction $\varphi\big|_\Na$ is semi-finite. A linear mapping $\E$ is said to be the conditional expectation of $\Ma$ onto $\Na$ with respect to $\varphi$ if the following properties are satisfied:
\begin{itemize}
\item[(i)] $\|\E(A)\|\leq\|A\|$ for all $A\in\Ma$;
\item[(ii)] $\E(A)=A$ for all $A\in\Na$;
\item[(iii)] $\varphi=\varphi\circ\E$.
\end{itemize}
\end{defi}
Faithfulness of $\varphi$ is the non-commutative analogue of a non-vanishing probability measure $P(A)$ if $A$ is not a null set whereas semi-finiteness and normality reduce to $\sigma$-finiteness and monotone convergence, respectively. Moreover, if we let 
\begin{gather*}
\Ma:=L^\infty(\Omega,\SA,P), X\in L^1(\Omega,\SA,P)\\
\text{and}\quad Y\in \Na:=L^1(\Omega,\SB,P), 
\end{gather*}
where $\SB\subset\SA$, then the classical case can be also easily recovered for the defining conditions as following:
\begin{itemize}
\item[(i)] $E\big(|E(X|\SB)|\big)\leq E(|X|)$;
\item[(ii)] $E(Y|\SB)=Y$;
\item[(iii)] $P\big(E(X|\SB)\big)=P(X)$.
\end{itemize}
The existence of a conditional expectation $\E$ of $\Ma$ onto $\Na$ is also equivalent to the invariance of the subalgebra with respect to the corresponding modular automorphism group, i.e. $\sigma_\varphi^t(\Na)=\Na$ for all $t\in\R$. Moreover, $\E$ can be shown to have the following properties:
\begin{gather*}
\E(X^*X)\geq0,\qquad X\in\Ma;\\
\E(AXB)=A\E(X)B,\qquad X\in\Ma,\,A,B\in\Na;\\
\E(X)^*\E(X)\leq\E(X^*X),\qquad X\in\Ma.
\end{gather*}
\begin{theo}
A conditional expectation $\E$ of $\Ma$ onto $\Na$ with respect to $\varphi$ exists if and only if $\Na\subset\Ma$ is invariant under the modular automorphism group $\sigma_\varphi^t$ with respect to $\varphi$, i.e. $\sigma_\varphi^t(\Na)=\Na$ for all $t\in\R$.
\end{theo}

\proof
Confer \cite[Theorem 4.2.]{Takesaki2:2001}.
\qed
\begin{defi}
Let $\Aa$ be a $C^*$-algebra and $\Ma$ a von Neumann algebra. Then a non-commutative stochastic process is defined as a sequence $(X_t)_{t\in\TT}$ of integrable random variables  $X_t:\Aa\longrightarrow\Ma$.
\end{defi}

\begin{prop}\label{condexpminimalaction}
If $\Ma$ is a von Neumann algebra, $\Na$ a von Neumann subalgebra of $\Ma$ with $\Na\rc\equiv\Na'\cap\Ma=\Za_\Na$, and $\E$ a normal projection from $\Ma$ onto $\Na$, then $\E$ is faithful and unique. Moreover, $\E$ is the conditional expectation of $\Ma$ onto $\Na$ with respect to $\varphi:=\psi\circ\E$ for all faithful, semifinite and normal weights $\psi$ on $\Na$.
\end{prop}

\proof
See \cite[Proposition 4.3.]{Takesaki2:2001}.
\qed

\begin{defi}
Let $(\Ma,\varphi)$ be a non-commutative probability space, then a family $(\Ma_t)_{t\in\TT}$ of subalgebras of $\Ma$ is said to be independent if the following conditions hold:
\begin{itemize}
\item[(i)] $[A_{t_j},A_{t_k}]=0$ for all $A_{t_l}\in\Ma_{t_l}$ and $t_j\neq t_k$;
\item[(ii)] $\varphi(A_{t_1}\cdots A_{t_n})=\varphi(A_{t_k})\cdots\varphi(A_{t_n})$ for all combinations of $A_{t_k}\in\Ma_{t_k}$, where $\{t_1,...,t_n\}$ is a non-void subset of $\TT$.
\end{itemize}
\end{defi}
There exists a weaker version of independence formulated in terms of conditional expectations which we will denote by $\E$-independence.
\begin{defi}
Let $\Ma_t$ and $\Na$, $\Ma_t\supset\Na$ for all $t\in\TT$, be subalgebras of $\Ma$ and $\E_\Na:\Ma_t\longrightarrow\Na$ a conditional expectation, then the family $(\Ma_t)_{t\in\TT}$ is called $\E$-independent if 
\begin{gather*}
\E_\Na(AB)=\E_\Na(A)\E_\Na(B),
\end{gather*}
where $A\in\Ma_t$, $t\in\TT$, and $B$ is an element of the von Neumann algebra generated by $(\Ma_s)_{s\neq t}$.
\end{defi}

\begin{defi}
Let $(\Ma,\varphi)$ be a non-commutative probability space, then a family of random variables $(X_t)_{t\in\TT}$ is said to be independent if the family $\big(X_t(\Aa)\big)_{t\in\TT}$ of the corresponding subalgebras has this property.
\end{defi}
In classical probability theory, independence in terms of the probability measure $P$ implies the one via the expectation. This is directly verified by the multiplication theorem for independent and integrable random variables $X_1, X_2,...,X_n$
\begin{gather*}
E\left(\prod_{i=1}^n X_i\right)=\prod_{i=1}^n E(X_i).
\end{gather*}
We continue with the non-commutative counterpart.
\begin{lem}
For algebras, independence implies $\E$-independence.
\end{lem}

\proof
For all $A\in\Ma_t$ and all elements $B$ of the von Neumann algebra generated by $(\Ma_s)_{s\neq t}$ we obtain:
\begin{align*}
\varphi\big(\E_\Na(A)\big)\varphi(B)&=\varphi\big(\E_\Na(A)B\big)\\
&=\varphi\big(\E_\Na(AB)\big)\\
&=\varphi(AB)\\
&=\varphi(A)\varphi(B)\\
&=\varphi\big(\E_t(A)\big)\varphi(B).
\end{align*}
Therefore we conclude
\begin{align*}
\varphi\big(\E_\Na(A)\big)&=\varphi\big(\E_t(A)\big),\\
\Longleftrightarrow\qquad\varphi\big(\E_\Na(A)-\E_t(A)\big)&=0 \\
\Longrightarrow\qquad\E_\Na(A)-\E_t(A)&=0 ,
\end{align*}
where we have used for the last implication faithfulness of $\varphi$. The last statement is equivalent to
\begin{gather*}
\E_\Na(AB)=\E_\Na(A)\E_\Na(B),
\end{gather*}
confer \cite[p.128]{Bauer:1995}.
\qed

\begin{defi}
A filtration in $\Ma$ of the non-commutative probability space $(\Ma,\varphi)$ is an increasing sequence $(\Ma_t)_{t\in\TT}$ of von Neumann subalgebras of $\Ma$ such that $\bigcup_{t\in\TT}\Ma_t$ is $\sigma$-weakly dense in $\Ma$. A family $(X_t)_{t\in\TT}$ of random variables on $(\Ma,\varphi)$ is said to be adapted with $(\Ma_t)_{t\in\TT}$ if $X_t\in L^1(\Ma_t)$. For a random variable $X$ on $(\Ma,\varphi)$ we define its conditional expectation with respect to $\Ma_t$ as
\begin{gather*}
\E_s(X_t):=\E(X_t|\Ma_s).
\end{gather*}
\end{defi}
Obviously, one has for two subalgebras $\Ma_t$ and $\Ma_s$ 
\begin{gather*}
\E_t\circ\E_s=\E_s\circ\E_t=\E_{\min\{s,t\}}.
\end{gather*}
The following diagram shall illustrate the aforementioned notions.\\
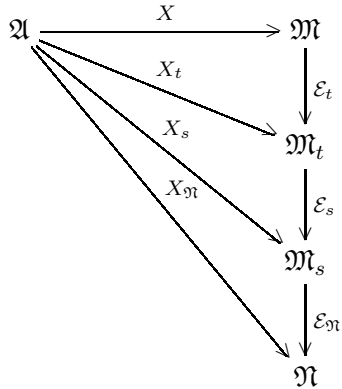
\begin{diagram}
\Aa  & \rTo^{\hspace{1,5cm} X}                                &&& \Ma\\
     &\rdTo(4,2)^{X_{t}}\rdTo(4,4)^{X_{s}}\rdTo(4,6)^{X_\Na}&&& \dTo_{\E_{t}}  \\
     &                                                        &&& \Ma_{t}\\
     &                                                        &&& \dTo_{\E_{s}}   \\ 
     &                                                        &&& \Ma_{s}\\
     &                                                        &&&  \dTo_{\E_\Na}\\
     &                                                        &&& \Na\\
\\
\end{diagram}
Now, we are in a position to define the main object of interest, namely a non-commutative martingale which will be investigated in Chapter five in the context of Galois correspondence and von Neumann subalgebras. 
\begin{defi}
A sequence $(X_t)_{t\in\TT}$ of random variables on $(\Ma,\varphi)$ is called a non-commutative martingale with respect to the filtration $(\Ma_t)_{t\in\TT}$ in $\Ma$ if for all $s,t\in\TT$ with $s<t$
\begin{gather*}
\E(X_{t}|\Ma_s)=X_s
\end{gather*}
holds.
\end{defi}



\chapter{Galois Correspondence and Non-Commutative Martingales}

\begin{flushright}
 \emph{What I tell you three times is true.}\\
\vspace{0,5cm}
Lewis Carroll\\
The Hunting of the Snark
\end{flushright}
\vspace{0,5cm}
The main results of this thesis are comprised in this chapter. We start with the discussion of invariant states with respect to automorphism groups and we show the existence of a non-trivial invariant von Neumann subalgebras.\\
The next section is concerned with Galois correspondence for compact automorphism groups and is opened by a summary of general notations and some well-known facts on the structure of von Neumann algebras, especially that of factors. The bottom-up analysis is divided in three consecutive cases, namely in inner, spatial and general automorphism groups. An short comparison with the approach of Izumi et al. \cite{Izumi:1998}  concludes this section.\\
Finally, we apply the whole apparatus and results to non-commutative probability theory as introduced in chapter five, i.e. non-commutative martingales are identified.






\section{Invariant Spaces}

According to the principle of Galois theory we assign to each subgroup $\Hr$ a subalgebra $\Ma^\Hr$ of $\Ma$ invariant under the action of $\Hr$, but first we should guarantee at least one of them not to be trivial. Here, this is done indirectly through the dual von Neumann algebra $\Ma^*$.

\begin{prop}
Let $\Gr$ be a compact group acting on a von Neumann algebra $\Ma$, $\alpha:\Gr\longrightarrow\autMa$, and define the annihilators
\begin{gather*}
\Ma_\perp^\Gr:=\Big\{\varphi\in\Ma^*|\;\varphi(A)=0\;\forall A\in\Ma^\Gr\Big\}\quad\text{and}\\
(\Ma^*)_\perp^\Gr:=\Big\{A\in\Ma|\;\varphi(A)=0\;\forall\varphi\in(\Ma^*)^\Gr\Big\}.
\end{gather*}
Then one has the following identities:
\begin{gather*}
\Ma=\Ma^\Gr+(\Ma^*)_\perp^\Gr,\\
\Ma^*=(\Ma^*)^\Gr+\Ma_\perp^\Gr.
\end{gather*}

\end{prop}

\proof Let us consider the following injections and their implication:
\begin{align*}
\Ma^\Gr\hookrightarrow\Ma\quad&\Longrightarrow\quad\Ma^*\longrightarrow(\Ma^\Gr)^*,\\
(\Ma^*)^\Gr\hookrightarrow\Ma^*\quad&\Longrightarrow\quad\Ma^{**}\longrightarrow\big[(\Ma^*)^\Gr\big]^*.
\end{align*}
The mappings on the right hand-side are due to the Hahn-Banach theorem both surjective. The second dualization  has been carried out with respect to the $\sigma(\Ma^*,\Ma)$-topology and therefore we may identify the double dual $\Ma^{**}$ with $\Ma$ itself. Thus we obtain, thanks to this identification, the next two injections
\begin{gather*}
i:\Ma^\Gr\longrightarrow\big[(\Ma^*)^\Gr\big]^*\quad\text{and}\\
i^*:(\Ma^*)^\Gr\longrightarrow(\Ma^\Gr)^*,
\end{gather*}
where $i^*$, obviously, is the dual mapping of $i$. Also here, we have used the closeness and therefore reflexivity of the subspace $(\Ma^*)^\Gr$, see Corollary \ref{reflexive1} . We want to investigate, firstly, the mapping $i$ and choose an arbitrary linear functional $\varphi$ in $\big[(\Ma^*)^\Gr\big]^*$, then we can find for it, due to the Hahn-Banach theorem, an extension $\Phi$ operating on the whole double dual space $\Ma^{**}\equiv\Ma$. Because our group $\Gr $ is compact, each $\Phi$ defines through the integral
\begin{gather*}
\int_\Gr\alpha^g(\Phi)d\mu_(g),
\end{gather*}
where $\mu(g)$ is the corresponding Haar measure, an invariant element of the algebra $\Ma$. Consequently, the mapping $i$ has to be surjective and therefore bijective:
\begin{gather*}
\Ma^\Gr\cong\big[(\Ma^*)^\Gr\big]^*.
\end{gather*}
With the same arguments we prove the isomorphism
\begin{gather}\label{MGS-MSG}
(\Ma^\Gr)^*\cong(\Ma^*)^\Gr.
\end{gather}
Each element $A\in\Ma$ defines through the equation
\begin{gather*}
A^\Gr:=\int_\Gr\alpha^g(A)d\mu_(g)
\end{gather*}
a $\Gr$-invariant element and their difference lies in $(\Ma^*)_\perp^\Gr$, as shown by the calculation:
\begin{align*}
\varphi\big(A-A^\Gr\big)&=\varphi(A)-\varphi\left(\int_\Gr\alpha^g(A)d\mu(g)\right)\\
&=\varphi(A)-\int_\Gr\varphi\big(\alpha^g(A)\big)d\mu(g)\\
&=\varphi(A)-\int_\Gr\varphi(A)d\mu(g)\\
&=\varphi(A)-\varphi(A)\int_\Gr d\mu(g)\\
&=0.
\end{align*} 
In the last equation we made use subsequently of the compactness of the group $\Gr$. Thus we may write
\begin{gather*}
\Ma=\Ma^\Gr+(\Ma^*)_\perp^\Gr.
\end{gather*}
In the same manner one may perform the proof of the second assertion.

\qed

\begin{theo}
To each $\sigma$-compact (locally compact) group $\Gr$ and von Neumann algebra $\Ma$ there always exists a non-trivial subspace $\Sa_\Aa^\Gr\subset\Sa_\Aa$ of $\Gr$-invariant states. 
\end{theo}

\proof Let $(U_n)_{n\in\N}$ be an increasing sequence of open neighborhoods of the unit element of $\Gr$ covering the group $\Gr$. Then there exists a 'locally left invariant' measure $\mu_{U_n}$ which is the restriction of the Haar measure $\mu_\Gr$ on $\Gr$, i.e. one has
\begin{gather}\label{localinvariance}
\int_{U_n^2}\alpha^{g'}\circ\alpha^g f(g)d\mu(U_n)=\int_{U_n^2}\alpha^g f(g)d\mu(U_n)
\end{gather}
for all $g\in\Gr$ and test function $f$ with $\supp f\subset U_n$.\\
Next choose a sequence of compact subsets, $(K_n)_{n\in\N}$, $K_n=\overline{U_n}$ with $\bigcup_{n=1}^\infty K_n=\Gr$ with
$$
\lim_{n\rightarrow\infty}\frac{\mu(K_n)}{\mu(K_{n+1})}=0.
$$
Since $\Sa_\Aa$ is $\sigma(\Aa^*,\Aa)$-compact one can always find a subsequence such that 
\begin{gather*}
\varphi:=\lim_{n\rightarrow\infty}\varphi_n:=\lim_{n\rightarrow\infty}\int_{K_n^2}\alpha^g(\phi)f(g)\frac{1}{\mu(K_n)}d\mu(K_n),
\end{gather*}
where $f$ is a test function with support in $K_n$ and $\phi\in\Sa_\Aa$, exists in the state space. Consequently, the state $\varphi$ is $\Gr$-invariant because due to \eqref{localinvariance} one derives
\begin{align*}
\alpha^{g'}(\varphi_n)&=\int_{K_n^2}\alpha^{g'}\circ\alpha^g(\phi)f(g)\frac{1}{\mu\big(\alpha^{g'}(K_n)\big)}d\mu(K_n)\\
&=\int_{K_n^2}\alpha^g(\phi)f(g)\frac{1}{\mu\big(\alpha^{g'}(K_n)\big)}d\mu(K_n)\\
&=\varphi_n
\end{align*}
for all $n\in\mathbb{N}$. Because each state is of norm one the limit cannot be zero.

\qed

\rem The construction of an $\Gr$-invariant element in $\Sa_\Aa$ would in general not apply to the algebra $\Aa$ itself because the limit of the integral could converge against zero. The subset $\Aa_\1$ would not be an alternative as it is not compact.

\begin{prop}
If $\Hr$ is a normal subgroup of the locally compact $\Gr$, then the subset $(\Ma^*)^\Hr$ of the dual space $\Ma^*$ is $\Gr$-invariant. 
\end{prop}

\proof Since $\Hr$ is a normal subgroup, for $h\in\Hr$ and $g\in\Gr$ the element $h':=ghg^{-1}$ is again contained in $\Hr$. But then one obtains for a functional $\varphi\in(\Ma^*)^\Hr$ with
\begin{gather*}
\varphi:=\lim_{n\rightarrow\infty}\varphi_n:=\lim_{n\rightarrow\infty}\int_{K_n}\alpha^g(\phi)f(g)\frac{1}{\mu(K_n)}d\mu(K_n),
\end{gather*}
where $K_n$
\begin{gather*}
\alpha^{g}(\varphi_n)=\int_{K_n^2}\alpha^{g}\circ\alpha^g(\phi)f(g)\frac{1}{\mu\big(\alpha^{g'}(K_n)\big)}d\mu(K_n).
\end{gather*}

\qed

\begin{theo}
Let $\Gr$ be a locally compact group acting on a $C^*$-algebra $\Aa$. Then the $\Gr$-ergodic states on $\Aa$, i.e. the pure states in $\Sa_\Aa^\Gr$, are contained in $\Pa_\Aa$, 
\begin{gather*}
\partial_e\Big(\Sa_\Aa^\Gr\Big)\subseteq\Pa_\Aa.
\end{gather*}
\end{theo}

\proof Let us consider the continuous unitary representation $U_\varphi$ of $\Gr$ on the Hilbert space $\Ha_\varphi$ defined by
\begin{gather*}
U_\varphi(g)\pi_\varphi(A)\xi_\varphi:=\pi_\varphi\circ\alpha^g(A)\xi_\varphi,
\end{gather*}
where $g\in\Gr$ and $A\in\Aa$. $U_\varphi$ is unique up to unitary equivalence and possesses the following properties:
\begin{gather*}
U_\varphi(g)\pi_\varphi(A)U_\varphi(g)^*=\pi_\varphi\circ\alpha^g(A)\\
\text{and}\quad U_\varphi(g)\xi_\varphi=\xi_\varphi.
\end{gather*}
Due to Proposition \ref{abelian-subalgebra} we can choose the maximal abelian von Neumann subalgebra $\Ca_\mu\subset\pi_\varphi(\Aa)'\cap U_\varphi(\Gr)'$ such that $\mu$ is the unique representing orthogonal measure of the $\Gr$-invariant state $\varphi$. With the choice of $\Aa^\Gr$ the conditions of Theorem \ref{pseudo-concentrated} are complete and, therefore, $\mu$ has to be pseudo-concentrated on $\Pa_{\Aa^\Gr}$, the extremal points of the state space on $\Aa^\Gr$. But due to the isomorphism 
\begin{gather*}
\big(\Ma^\Gr\big)^*\equiv(\Ma^*)^\Gr,
\end{gather*}
proven in \eqref{MGS-MSG} $\mu$ is pseudo-concentrated of $\partial_e\big(\Sa_\Aa^\Gr\big)$. The measure is due to Proposition \ref{abelian-subalgebra} unique and, as an orthogonal measure thanks to Proposition \ref{boundary-extreme} and Proposition \ref{orthogonal-boundary}, must also be pseudo-concentrated on the pure states of $\Aa$. Consequently, we obtain the inclusion $\partial_e\big(\Sa_\Aa^\Gr\big)\subseteq\Pa_\Aa$.
\qed

\section{Galois Correspondence for Compact Groups}

First of all, we want to clarify our notation, where we follow standard definitions, and give some vital well-known facts.\\
Let $\Gr$ be a compact group and $M(\Gr)$ the vectors pace of all bounded Radon measures on $\Gr$. On $\Cg(\Gr)$ one introduces a continuous linear form, the convolution of the Radon measures $\mu,\nu\in M(\Gr)$ by
\begin{gather*}
(\mu\star\nu)(f):=\int_\Gr\int_\Gr f(g_1g_2)d_\mu(g_1)d_\nu(g_2),\quad f\in\Cg(\Gr),g_1,g_2\in\Gr,
\end{gather*}
satisfying the following well-known properties for all $\mu,\nu,\nu'\in M(\Gr)$:
\begin{itemize}
\item[(i)] $(\alpha\mu)\star\nu=\mu\star(\alpha\nu)=\alpha(\mu\star\nu)$;
\item[(ii)] $\mu\star(\nu+\nu')=\mu\star\nu+\mu\star\nu'$;
\item[(iii)] $\mu\star(\nu\star\nu')=(\mu\star\nu)\star\nu'$;
\item[(iv)] $\|\mu\star\nu\|\leq\|\mu\|\|\nu\|$;
\item[(v)] $M(\Gr)$ is commutative if an only if $\Gr$ is so.
\end{itemize}
For $f,h\in L^2(\Gr)$ the convolution is defined as
\begin{gather*}
(f\star h)(g_1):=\int_\Gr f(g_2)h\big(g_2^{-1}g_1\big)dg_2,
\end{gather*}
and the convolution of $\mu\in M(\Gr)$ and $f\in L^1(\Gr)$ by
\begin{gather*}
(\mu\star f)(g_1):=\int_\Gr f\big(g_2^{-1}g_1\big)d_\mu(g_2),
\end{gather*}
which determines 
\begin{gather*}
(f\star\mu)(g_1):=\int_\Gr\Delta\big(g_2^{-1}\big) f\big(g_1g_2^{-1}\big)d_\mu(g_2),
\end{gather*}
where $\Delta$ is the modular function. For $f,h\in L^p(\Gr)$, $1\leq p\leq\infty$, $\frac{1}{p}+\frac{1}{q}=1$ and $\mu\in M(\Gr)$ the functions $f\star h$, $\mu\star f$ and $f\star\mu$ are also contained in $L^p(\Gr)$. Moreover, one has
\begin{equation}\label{measurenorm}
\begin{split}
\|f\star h\|_p&\leq\|f\|_1\|h\|_p,\\
\|\mu\star f\|_p&\leq\|\mu\|\|f\|_p,\\
\|f\star\mu\|_p&\leq\|f\|_p\int_\Gr\Delta\big(g^{-1}\big)^{1/q}d|\mu|(g).
\end{split}
\end{equation}
The adjoint Radon measure is defined by
\begin{gather*}
\mu^*(f):=\overline{\int_\Gr\overline{f\big(g^{-1}\big)}d\mu(g)}=\int_\Gr f\big(g^{-1}\big)\overline{d\mu(g)}=\overline{\mu(\overline{f^*})}
\end{gather*}
where we set $f^*(g):=f(g^{-1})$. The adjoint measure satisfies in addition to the defining conditions of an involution also 
\begin{equation}\label{betatildestar}
\begin{split}
|\mu^*|=|\mu|^*\quad\quad\|\mu^*\|=\|\mu\|,\notag\\
\int_\Gr f(g)d\mu^*(g)=\overline{\int_\Gr\overline{f\big(g^{-1}\big)}d\mu(g)},\quad\text{and}\notag\\
(\mu^*\star f)(g_1)=\overline{\int_\Gr\overline{f(g_2g_1)}d\mu(g_2)}.
\end{split}
\end{equation}
In particular $\Gr$ can be represented on the space of square integrable functions on $\Gr$:
\begin{align*}
\beta:\Gr&\longrightarrow\L\big(L^2(\Gr)\big)\\
g_1&\mapsto\beta(g_1),\\
 \big(\beta(g_1)f\big)(g_2)&:=f\big(g_1^{-1}g_2\big),\quad f\in L^2(\Gr).
\end{align*}
The mapping of $\Gr$ on $M(\Gr)$ is denoted by
\begin{align*}
\lambda:\Gr&\longrightarrow M(\Gr)\\
g&\mapsto\delta_g,
\end{align*}
which transforms the multiplication in $\Gr$ to the convolution of measures
\begin{gather*}
g_1\cdot g_2\mapsto\delta_{g_1}\star\delta_{g_2}. 
\end{gather*}
By the following mapping we may interprete the measures on $\Gr$ as operators on $L^2(\Gr)$:
\begin{align*}
\tilde{\beta}:M(\Gr)&\longrightarrow\L\big(L^2(\Gr)\big)\\
\delta_g&\mapsto\tilde{\beta}(\delta_{g_1})\\
\big(\tilde{\beta}(\delta_{g_1})f\big)(g_2)&:=f\big(g_1^{-1}g_2\big), g_1,g_2\in\Gr,f\in L^2(\Gr).
\end{align*}
where $g_1,g_2\in\Gr$ and $f\in L^2(\Gr)$. $\tilde{\beta}$ is then defined for an arbitrary $\mu$ through the one of one-point measures $\delta_g$,
\begin{align*}
\tilde{\beta}(\mu)(f):=\int_\Gr\tilde{\beta}(\delta_g)d\mu(g)=\mu\star f,
\end{align*}
and the norm in $\L\big(L^2(\Gr)\big)$ is given by
\begin{gather*}
\|\tilde{\beta}(\mu)\|:=\sup_{\substack{f\in L^2(\Gr)\\ \|f\|_2\leq1}}\|\mu\star f\|_2.
\end{gather*}
Since $\mu\star f$ is finite, see Equations \ref{measurenorm}, $\tilde{\beta}(\mu)$ constitutes for all $\mu\in M(\Gr)$ a bounded operator on $L^2(\Gr)$.
\begin{lem}
The mapping $\tilde{\beta}$ is a faithful representation of $M(\Gr)$ on the Hilbert space $L^2(\Gr)$.
\end{lem}

\proof
Linearity of $\tilde{\beta}$ follows from the definition of $\mu\star f$. Moreover, homomorphism property is also given:
\begin{align*}
\tilde{\beta}(\mu\star\nu)(f)&=(\mu\star\nu)(f)\\
&=\int_\Gr\int_\Gr f\big(g_2^{-1}g_1\big)d{\mu\star\nu}(g_2)\\
&=\int_\Gr\int_\Gr f\big(g_2^{-1}g_1\big)d{\nu}(g_2)d_\mu(g_2)\\
&=\mu\star(\nu\star f)\\
&=\tilde{\beta}(\mu)\big(\tilde{\beta}(\nu)(f)\big).
\end{align*}
$\tilde{\beta}$ is a $\*$-homomorphism because, due to Fubini and Equation \ref{betatildestar}, for $f,h\in\Cg(\Gr)$ one obtains:
\begin{align*}
\big\langle f,\tilde{\beta}(\mu^*)(h)\big\rangle&=\int_\Gr f(g_1)\int_\Gr \overline{h(g_2g_1)}d\mu(g_2)dg_1\\
&=\int_\Gr\int_\Gr f(g_1)\overline{h(g_2g_1)}dg_1d\mu(g_2)\\
&=\int_\Gr\int_\Gr f\big(g_2^{-1}g_1\big)\overline{h(g_1)}dg_1d\mu(g_2)\\
&=\int_\Gr\int_\Gr f\big(g_2^{-1}g_1\big)d\mu(g_2)\overline{h(g_1)}dg_1\\
&=\big\langle\tilde{\beta}(\mu^*)(f),h\big\rangle.
\end{align*}
Consider now a non-vanishing measure $\mu\neq0$, thus $\int_\Gr f\big(g^{-1}\big)d\mu\neq0$, which passes over its continuity to $\mu\star f$. But $\mu\star f(\1)=\int_\Gr f\big(g^{-1}\big)d\mu$ and therefore $\|\mu\star f\|^2=\int_\Gr|\mu\star f|^2d\lambda$ is positive, i.e. $\tilde{\beta}(\mu)\neq0$.
\qed

\begin{defi}
Let $\Sigma\equiv\Sigma(\Gr)$ be the system of equivalence classes of finite-dimensional, irreducible, continuous, unitary representations of $\Gr$ then we define on the Hilbert subspace $\Ha_\sigma$, $\sigma\in\Sigma$, the representation:
\begin{align*}
\tilde{\beta}_\sigma:M(\Gr)&\longrightarrow\L(\Ha_\sigma)\\
\mu&\mapsto\tilde{\beta}_\sigma(\mu)(f):=\int_\Gr\sigma^g(f)d\mu(g)
\end{align*}
\end{defi}

\begin{lem}
One has $\im\tilde{\beta}\subseteq(\im\beta)''$.
\end{lem}

\proof
The statement follows from the fact that $\im\tilde{\beta}$ lies in the weak closure of all linear combinations of $\delta_g$, $g\in\Gr$. 
\qed

\begin{defi}
Let $\big(\Ma_i,\pi_i,\Ha_i\big)$, $i\in I$, be a triple consisting of von Neumann algebras $\Ma_i$ acting on the Hilbert spaces $\Ha_i$, i.e. $\pi_i:\Ma_i\longrightarrow\L(\Ha_i)$ for all $i$. If we set for the direct sums $\Ha:=\sum^\oplus\Ha_i$ and $\xi:=\sum^\oplus\xi_i$, $\xi_i\in\Ha_i$, then we define the direct sum of the representation as
\begin{gather*}
\pi(A)\xi:=\sum_{i\in I}^\oplus\pi_i(A)\xi_i.
\end{gather*}
For each bounded sequence $(A_i)_{i\in I}$, $A_i\in\Ma_i$, we define an operator on $\Ha$ by
\begin{gather*}
A\,\sum_{i\in I}^\oplus\xi_i:=\sum_{i\in I}^\oplus A_i\xi_i
\end{gather*}
and call the von Neumann algebra generated by these operators the direct sum of $\Ma_i$, denoted as $\Ma:=\sum_{i\in I}^\oplus\Ma_i$. 
\end{defi}

\begin{defi}
If $\Aa$ and $\Ba$ are two $C^*$-algebras, then we denote by $\Aa\otimes\Ba$ the algebraic tensor product of them and follow the common definitions:
\begin{gather*}
(A_1\otimes B_1)(A_2\otimes B_2):=A_1A_2\otimes B_1B_2,\\
(A_1\otimes B_1)^*:=A_1^*\otimes B_1^*
\end{gather*}
for all $A_1,A_2\in\Aa$ and $B_1,B_2\in\Ba$. The injective $C^*$-tensor product $\Aa\otimes_{\text{min}}\Ba$ is defined as the completion of $\Aa\otimes\Ba$ with respect to the so-called injective $C^*$-crossnorm
\begin{gather*}
\|A\|_{\text{min}}:=\sup\|(\pi_1\otimes\pi_2)(A\otimes B)\|,\quad (A\otimes B)\in\Aa\otimes\Ba.
\end{gather*}
\end{defi}
Let us now consider two $W^*$-algebras $\Ma$ and $\Na$ with preduals $\Ma_*$ and $\Na_*$, respectively. First of all, the closure $\Ma_*\overline{\otimes}\,\Na_*$ of $\Ma_*\otimes\,\Na_*$ is seen as a closed subset in $(\Ma\otimes_{\text{min}}\Na)^*$. Since $\Ma_*\otimes\,\Na_*$ is invariant under the algebraic tensor product $\Ma\otimes\Na$, and therefore so is $\Ma_*\overline{\otimes}\,\Na_*$, there exists uniquely a central projection $Z$ in the universal enveloping $(\Ma\otimes_{\text{min}}\Na)^\sim$ of $\Ma\otimes_{\text{min}}\Na$ such that $\Ma_*\overline{\otimes}\,\Na_*=(\Ma\otimes_{\text{min}}\Na)^*Z$. The fact that $(\Ma\otimes_{\text{min}}\Na)$ is embedded in $(\Ma_*\otimes\Na_*)^*$ which in turn is isomorphic to the $W^*$-algebra $(\Ma\otimes_{\text{min}}\Na)^\sim Z$, justifies the following definition.

\begin{defi}
The $W^*$-tensor product $\Ma\overline{\otimes}\,\Na$ of two $W^*$-algebras $\Ma$ and $\Na$ is defined as the $W^*$-algebra $(\Ma\otimes_{\text{min}}\Na)^\sim Z$.
\end{defi}
The table given below illustrates the resulting type of the crossproduct $\Ma\overline{\otimes}\Na$ of two von Neumann algebras $\Ma$ and $\Na$. 
\vspace{1cm}

\begin{center}
\begin{tabular}{c|ccccc}
$\Ma\diagdown\Na$ & $I_n$ & $I_\infty$ & $II_1$ & $II_\infty$ & $III$\\
\hline
$I_m$ & $I_{mn}$ & $I_\infty$ & $II_1$ & $II_\infty$ & $III$\\
$I_\infty$ & $I_\infty$ & $I_\infty$ & $II_\infty$ & $II_\infty$ & $III$\\
$II_1$ & $II_1$ & $II_\infty$ & $II_1$ & $II_\infty$ & $III$\\
$II_\infty$ & $II_\infty$ & $II_\infty$ & $II_\infty$ & $II_\infty$ & $III$\\
$III$ & $III$ & $II$ & $III$ & $III$ & $III$
\end{tabular}
\end{center}

\subsection{Groups of Inner Automorphisms}

We begin our analysis with compact groups $\Gr$. We denote the system of equivalence classes of finite-dimensional, irreducible, continuous, unitary representations of $\Gr$ by $\Sigma\equiv\Sigma(\Gr)$. We will follow the standard definition of the following spaces:
\begin{gather*}
\La_{\Sigma,\Gr}^p:=l^p\left(\L\big(\Ha_{\sigma,\Gr}\big)_{\sigma\in\Sigma}\right):=\left\{(T_\sigma)_{\sigma\in\Sigma}\Big|\;\Big(\sum_\sigma\|T_\sigma\|^p_{\text{uniform}}\Big)^{1/p}\hspace{-3mm}<\infty,T_\sigma\in\L\big(\Ha_{\sigma,\Gr}\big)\right\},\\
\Lainf_{\Sigma,\Gr}:=l^\infty\left(\L\big(\Ha_{\sigma,\Gr}\big)_{\sigma\in\Sigma}\right):=\Big\{(T_\sigma)_{\sigma\in\Sigma}\big|\;\sup_\sigma\|T_\sigma\|_{\text{uniform}}<\infty,\;T_\sigma\in\L\big(\Ha_{\sigma,\Gr}\big)\Big\},\\
\La_{\Sigma,\Gr}^0:=l^0\left(\L\big(\Ha_{\sigma,\Gr}\big)_{\sigma\in\Sigma}\right):=\left\{(T_\sigma)_{\sigma\in\Sigma}\big|\;\lim_\sigma\|T_\sigma\|_{\text{uniform}}=0,\;T_\sigma\in\L\big(\Ha_{\sigma,\Gr}\big)\right\},\\
\Ha_{\Sigma,\Gr}^2:=l^2\left(\bigoplus_{n=1}^{n_\sigma}\Ha_{\sigma,\Gr}\right):=\left\{\big((\xi_\sigma^1,...,\xi_\sigma^{n_\sigma})\big)_{\sigma\in\Sigma}\big|\;\sum_{k=1}^{n_\sigma}\|\xi_\sigma^k\|^2<\infty\right\},
\end{gather*}
where $n_\sigma$ is the dimension of $\Ha_{\sigma,\Gr}$. $\Lainf_{\sigma,\Gr}$ is a Banach space of bounded functions on $\L(\Ha_\sigma)$. Because its elements can be seen as bounded operators on $\Ha_{\sigma,\Gr}$, it is together with the complex conjugation as involution also a von Neumann algebra in $\LH$.\\
Let $\Ma$ be a von Neumann algebra of type I acting on the Hilbert space $\Ha$ and $\{Z_i\}$ its central projections with $\sum_i Z_i=\1$ and $Z_i$ being the greatest $\alpha$-homogenous projection, i.e. the sum of $i$ orthogonal abelian projections having the central support $Z$. Since the commutant of $\Ma$ is automatically also of type I we may choose central projections $\{Z_j\}$ for $\Ma'$ with the same properties. Due to $\Za_{\Ma'}=\Ma'\cap\Ma''=\Ma'\cap\Ma=\Za_\Ma$ the projections lie in the same center $\{Z_i\},\{Z_i\}\in\Za_\Ma$ and therefore we may introduce a finer partition of the identity by
\begin{gather*}
Z_{i,j}:=Z_iZ_j,\quad \sum_{i,j}Z_{i,j}=\1.
\end{gather*}
In order to decompose the pair $(\Ma,\Ha)$, let $\Ha_{i,j}:=Z_{i,j}\Ha$, $\Ma_{i,j}:=Z_{i,j}\Ma$, $\Na_{i,j}$ an abelian von Neumann algebra represented on $\Ka_{i,j}$ as maximal abelian which is isomorphic to $\Za_{\Ma_{i,j}}$ and $\Ha_i$ and $\Ha_j$ $i$-dimensional and $j$-dimensional Hilbert spaces, respectively. Then one has the following isomorphism
\begin{gather*}
\Ma_{i,j}\cong\Na_{i,j}\overline{\otimes}\,\L(\Ha_i).
\end{gather*}
Let $\{E_k\}$ be an orthogonal family of projections in $\Na_{i,j}\otimes\C\overline{\otimes}\,\L(\Ma)$ such that for all $k$ on has $Z(E_k)=\1$, $\sum_kE_k=\1$ and
\begin{align*}
(\Ma_{i,j},\Ha_{i,j})\cong\Big(E\big(\Na_{i,j}\overline{\otimes}\,\L(\Ha_i)\otimes\C\big),E\big(\Na_{i,j}\otimes\Ha_i\otimes\Ma\big)\Big).\\
\end{align*}
If one chooses an orthogonal family of abelian projections $\{F_l:=\1\times P_l\}$, where $\{P_l\}$ is a set of minimal projections in $\L(\Ma)$ and set $F:=\sum_lF_l$ and $P:=\sum_lP_l$, then one has in $\Na_{i,j}\otimes\C\overline{\otimes}\,\L(\Ma)\big)$ the equivalence $E\sim F$ leading finally to
\begin{align*}
(\Ma_{i,j},\Ha_{i,j})&\cong\Big(E\big(\Na_{i,j}\overline{\otimes}\,\L(\Ha_i)\otimes\C\big),E\big(\Na_{i,j}\otimes\Ha_i\otimes\Ma\big)\Big)\\
&\cong\Big(F\big(\Na_{i,j}\overline{\otimes}\,\L(\Ha_i)\otimes\C\big),F\big(\Na_{i,j}\otimes\Ha_i\otimes\Ma\big)\Big)\\
&\cong\Big(\big(\Na_{i,j}\overline{\otimes}\,\L(\Ha_i)\otimes\C\big),\big(\Na_{i,j}\otimes\Ha_i\otimes\Ha_j\big)\Big).
\end{align*}
We formulate this result for type I von Neumann algebras in the next theorem.

\begin{theo}\label{TakesakiStruktur}
A von Neumann algebra $(\Ma,\Ha)$ of type I has the unique decomposition:
\begin{align*}
(\Ma,\Ha)\cong\sum_{i,j}^\oplus(\Na_{i,j},\Ka_{i,j})\overline{\otimes}\big(\L(\Ha_i),\Ha_i\big)\overline{\otimes}(\C,\Ha_j),
\end{align*}
where $(\Na_{i,j},\Ka_{i,j})$ is maximal abelian and $\Ha_i$ and $\Ha_j$ are $i$-dimensional and $j$-dimensional, respectively.
\end{theo}

\proof
Done, see \cite[V.1.31.]{Takesaki1:1970}.
\qed



\begin{defi}
A representation $\pi:\Aa\longrightarrow\LH$ of a $C^*$-algebra $\Aa$ on a Hilbert space $\Ha$ is said to be of type I if the von Neumann algebra generated by $\pi(\Aa)$ is of type I.
\end{defi}

\begin{prop}\label{typeI}
Let $\{\pi_i\}$, $i\in I$ for some index set $I$, be a family of pairwise disjoint representations of a von Neumann algebra $\Ma$ and set $\pi:=\sum_{i\in I}^\oplus\pi_i$, then the following statements are equivalent:
\begin{itemize}
\item[(i)] $\pi$ is of type I.
\item[(ii)] $\pi_i$ is of type I for all $i\in I$.
\end{itemize}
\end{prop}

\proof
The implication $(i)\Longrightarrow(ii)$ is a direct consequence of Theorem \ref{TakesakiStruktur}

\qed

\begin{lem}
All von Neumann subalgebras $\Aa$ of $\Lainf_\sigma$ are of type I and have the following structure:
\begin{gather*}
(\Aa,\Ha_\sigma)\cong\sum_{i,j}^\oplus\sum_k(\C_{i_k,j_k},\Ka_{i_k,j_k})\overline{\otimes}\big(\L(\Ha_i),\Ha_i\big)\overline{\otimes}(\C,\Ha_j).
\end{gather*}
\end{lem}

\proof
Since $\Lainf_{\Sigma,\Gr}$ is a von Neumann algebra of type one we may apply Theorem \ref{TakesakiStruktur}. Let $\{E_k\}$, $k\in I$ for some index set $I$, be a family of projections in $\Aa'$ with $\sum_kE_k=\1$, the we obtain
\begin{align*}
\big(\Aa,\Ha_{\sigma,\Gr}\big)&\cong\big(E_k\Lainf_{\sigma,\Gr},E_k\Ha_{\sigma,\Gr}\big)\\
&\cong\sum_{i,j}^\oplus(E_k\Na_{i,j},E_k\Ka_{i,j})\overline{\otimes}\big(E_k\L(\Ha_i),E_k\Ha_i\big)\overline{\otimes}(E_k\C,E_k\Ha_j)\\
&=\sum_{i,j}^\oplus\sum_k(\C_{i_k,j_k},\Ka_{i_k,j_k})\overline{\otimes}\big(\L(\Ha_i),\Ha_i\big)\overline{\otimes}(\C,\Ha_j).
\end{align*}
Each part of the decomposition is finite-dimensional and therefore of type I. Due to Proposition \ref{typeI}, the von Neumann subalgebra $\Aa$ as a direct sum of them has to be also of type I.
\qed

\begin{lem}\label{vonNeumannsubalg}
Let $\Ma$ be a von Neumann algebra and $\Gr$ a compact group, then one has the following identity:
\begin{gather*}
\big(\Ma^\Gr\big)\rc\rc=\Ma^\Gr.
\end{gather*}
\end{lem}

\proof 
One has obviously $\Ma^\Gr\subseteq\big(\Ma^\Gr\big)\rc\rc$. Let us assume that $\alpha^g$, $g\in\Gr$, is an ${}^*$-automorphism of $\Ma$ with its fix-algebra $\Ma^\Gr$. Its generator $\delta$ is a bounded inner ${}^*$-derivation and there exists a selfadjoint element $H\in\Ma$ such that $\delta(A)=i[H,A]$ for all $A\in\Ma$. Since an element $A\in\Ma$ belongs to the fix-algebra if and only if $\delta(A)=0$, we may write:
\begin{align*}
\Ma^\Gr=\big\{A\in\Ma|\;\delta(A)=i[H,A]=0\big\}.
\end{align*}
Thus $H$ lies in $\big(\Ma^\Gr\big)\rc$, the relative commutant of $\Ma^\Gr$. This leads to the following chain of implications:
\begin{align*}
A\in\big(\Ma^\Gr\big)\rc\rc\quad&\Longrightarrow\quad AB=BA,\quad\forall B\in\big(\Ma^\Gr\big)\rc\\
&\Longrightarrow\quad AH=HA\\
\quad&\Longrightarrow\quad\alpha^\Gr(A)=HAH^*=A.
\end{align*}
Thus one obtains also $\big(\Ma^\Gr\big)\rc\rc\subseteq\Ma^\Gr$.

\qed

\begin{lem}\label{mappingpi}
The image of the standard implementation 
\begin{align*}
\pi:=\pi_\U^{-1}:\autMa&\longrightarrow\U(\Ma)\\
\alpha&\mapsto U,
\end{align*}
defined in Theorem \ref{standardimplementation}, is contained in $\big(\Ma^\Gr\big)\rc$, the relative commutant of the fixed point algebra $\Ma^\Gr:=\{A\in\Ma|\;g(A)=A,\;\forall g\in\Gr\}$.
\end{lem}

\proof
If $A\in\Ma^\Gr$, i.e. $\alpha^t(A)=UAU^*=A$, then $A\in U(\Ma)\rc$. But due to Lemma \ref{vonNeumannsubalg} $\big(\Ma^\Gr\big)\rc\rc=\big(\Ma^\Gr\big)\subseteq U(\Ma)\rc$ leads to $\big(\Ma^\Gr\big)\rc\supseteq U(\Ma)$. 
\qed

\begin{lem}\label{universalmap}
If $(\pi,\Ha)$ is a representation of the $C^*$-algebra $\Aa$ and $\Ma:=\pi(\Aa)''$ the von Neumann algebra generated by $\Aa$, then there exists a unique linear mapping $\bar{\pi}$ of the second conjugate $\Aa^{**}$ of $\Aa$ onto $\Ma$ with the following properties:
\begin{itemize}
\item[(i)] $\bar{\pi}\circ i=\pi$ where $i$ denotes the canonical embedding of $\Aa$ in $\Aa^{**}$;
\item[(ii)] $\bar{\pi}$ is continuous with respect to the $\sigma$-weak topology.
\end{itemize} 
\end{lem}

\proof \cite[III.2.2.]{Takesaki1:1970}
Let us consider the Banach space $\Ma_*$ of all $\sigma$-weakly continuous linear functionals on $\Ma$ and let $\pi_*$ be the restriction of the transpose $\pi^T$ of the representation $\pi$, i.e. $\pi_*:=\pi^T\big|_{\Ma_*}$. We then define 
\begin{gather*}
\bar{\pi}:=(\pi_*)^T,
\end{gather*}
which constitutes, as the von Neumann algebra $\Ma$ is the conjugate space of $\Ma_*$, a map from $\Aa^{**}$ onto $\Ma$.  The second property is ensured by definition of $\bar{\pi}$ and statement $(i)$ follows directly by
\begin{gather*}
\big\langle\pi(A),\o\big\rangle=\big\langle i(A),\pi_*(\o)\big\rangle=\big\langle\bar{\pi}\circ i(A),\o\big\rangle,\quad\forall\o\in\Ma_*,\,\forall A\in\Ma.
\end{gather*}
\qed

\begin{defi}
We call $\pi(\Gr)$ proper, if $\int_\Gr\pi(g)d\mu(g)\neq0$ for $\mu\neq 0$.
\end{defi}

\begin{rem}
The center $\Za$ of $\Ma$ lies in $\Ma^\Gr\cap(\Ma^\Gr)\rc$ and $\pi(\Gr)\cap\Za=\S^1$.
\end{rem}
Due to Theorem \ref{standardimplementation} the following maps of $L^1(\Gr)$ and $M(\Gr)$ on $\left(\Ma^\Gr\right)\rc$ are well-defined:
\begin{align*}
\tilde{\pi}_{L^1}:L^1(\Gr)&\longrightarrow\left(\Ma^\Gr\right)\rc\\
f&\mapsto\tilde{\pi}_{L^1}(f):=\int_\Gr f(g)\pi(g)dg,\quad g\in\Gr,\\
\tilde{\pi}_M:M(\Gr)&\longrightarrow\left(\Ma^\Gr\right)\rc\\
\mu&\mapsto\tilde{\pi}_M(\mu):=\int_\Gr\pi(g)d\mu(g),\quad g\in\Gr.
\end{align*}

\begin{rem}
$\pi(\Gr)$ is proper if and only if $\tilde{\pi}_M$ is injective.\\
\end{rem}
Definition of $\rep(\Gr)$ as the equivalence classes of all representations, $\Sigma$ are irreducible ones. 

\begin{theo}\label{compactinner}
Let $\Gr$ be a compact group and $\Ma$ a von Neumann algebra, then there exists a $\Sigma':=\Sigma\backslash\Sigma^0\subseteq\Sigma$, $\Sigma^0:=\ker(\tilde{\pi}_M\circ\iota_1\circ\iota_2)$, such that the following diagram is commutative:
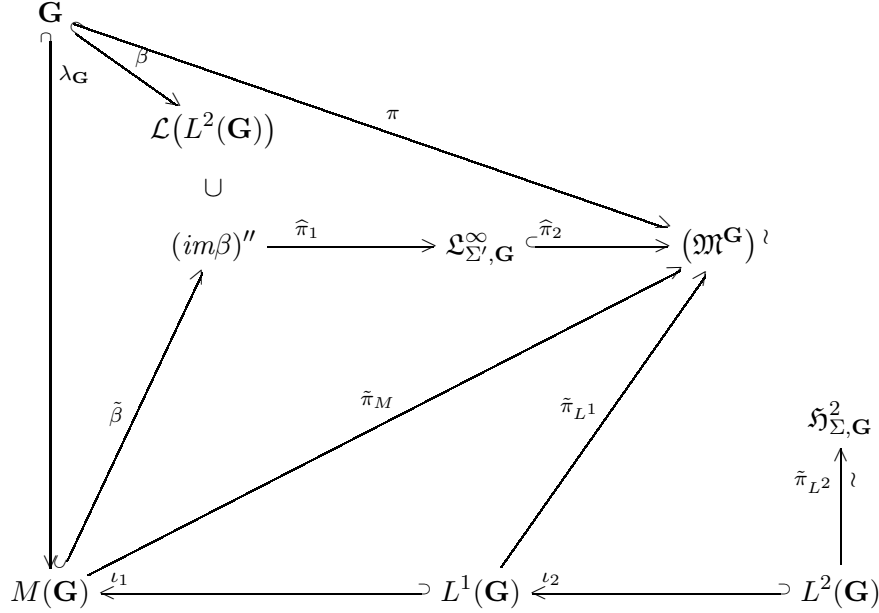
\begin{diagram}
\Gr       &     &      &     &&& &          &&&&&\\
\dInto(0,12)_{\lambda_\Gr}&\rdInto(2,2)^\beta\rdTo(10,4)^{\pi}&&&&&\\ 
         &          &\L\big(L^2(\Gr)\big) \\
         &          &\cup \\
         &          &(\im\beta)''    &\rTo^{\widehat{\pi}_1} &&&\Lainf_{\Sigma',\Gr}&\rInto^{\widehat{\pi}_2}      &&&\left(\Ma^\Gr\right)\rc&&\\
         &\ruInto(2,6)^{\tilde{\beta}}&      &     &&& &    &&\ruTo(10,6)^{\tilde{\pi}_M}\ruTo(4,6)^{\tilde{\pi}_{L^1}}&&&\\
\\ 
         &          &      &     &&& &          &&&&&\Ha_{\Sigma,\Gr}^2\\
         &          &      &     &&& &          &&&&&\uTo(0,4)_\wr^{\tilde{\pi}_{L^2}}\\
\\
M(\Gr)    &\lInto^{\iota_1}     &      &     &&&L^1(\Gr)&\lInto^{\iota_2}      &&&&&L^2(\Gr)
\end{diagram}
\end{theo} 

\begin{rem}
In the case of $\pi(\Gr)$ being proper $\tilde{\pi}_M$ is injective and hence $\Sigma'=\Sigma$.
\end{rem}

\hspace{-5,5mm}The proof of this theorem will be based on the following lemmas.

\begin{lem}\label{pihat}
The mapping $\widehat{\pi}:=\tilde{\pi}_M\circ\tilde{\beta}^{-1}$ is ${}^*$-preserving and can be extended to $(\im\beta)''$.
\end{lem}

\proof
As aforementioned, the norm $\|\tilde{\beta}(\mu)\|$ in $(\im\beta)''$ is given by 
\begin{align*}
\|\tilde{\beta}(\mu)\|=\sup_{\substack{f\in L^2(\Gr)\\ \|f\|_2\leq1}}\|\mu\star f\|_2.
\end{align*}
and the norm $\|\tilde{\pi}_M(\mu)\|$ in $(\Ma^\Gr)\rc$ is 
\begin{align*}
\|\tilde{\pi}_M(\mu)\|=\|\tilde{\pi}(\mu)\tilde{\pi}^*(\mu)\|^{1/2n}.
\end{align*}
Due to 
$$\|\mu\star\mu^*\|\geq\|\tilde{\pi}_M(\mu\star\mu^*)\|=\|\tilde{\pi}(\mu)\tilde{\pi}^*(\mu)\|
$$
the mapping $\tilde{\pi}_M\circ\tilde{\beta}^{-1}$ is continuous with respect to the uniform topology and can be continuously extended to $(\im\beta)''$.

\qed

\begin{lem}\label{propersigma}
$\pi(\Gr)$ is proper if and only if  one has $\Sigma=\Sigma'$, i.e. if $\widehat{\pi}=\widehat{\pi}_2\circ\widehat{\pi}_1$ is injective.
\end{lem}

\proof
Since $\tilde{\beta}$ is injective, faithfulness implies properness. 
\qed
\begin{lem}\label{imbeta}
One has the following identity:
\begin{align*}
(\im\beta)''&=\Lainf_{\Sigma,\Gr}.
\end{align*}
\end{lem}

\proof
One has the inclusion $\Lainf_{\Sigma,\Gr}\subset\im\tilde{\beta}$, since $\L(\Ha_{\sigma,\Gr})$ is contained in $\im\tilde{\beta}$ for all $\sigma\in\Sigma$ . Let $\chi_\sigma$ be the character with respect to $\sigma\in\Sigma$, then due to Lemma \ref{character} $e_\sigma:=\tilde{\beta}(\chi_\sigma)$ is a projector on $\L\big(\Ha_{\sigma,\Gr}\big)$. If $\L\big(\Ha_{\sigma,\Gr}\big)\subset\La_{\Sigma,\Gr}^0$, then one gets
\begin{gather*}
\L\big(\Ha_{\sigma,\Gr}\big)\cong\tilde{\beta}\big(M(g)\big)e_\sigma=\tilde{\beta}\big(M(g)\big)\tilde{\beta}(\chi_\sigma)=\tilde{\beta}\big(M(g)\star(\chi_\sigma)\big),
\end{gather*}
thus $\im\tilde{\beta}$ is a dense subset of $\La_\Sigma^0$.
 Moreover, by well-known theorems one obtains
\begin{align*}
\left(\La_{\Sigma,\Gr}^0\right)^*&\cong\La_{\Sigma_0,\Gr}^1:=l^1\left(\mathcal{N}\big(\Ha_{\sigma,\Gr}\big)_{\sigma\in\Sigma_0}\right)\\
&:=\left\{(T_\sigma)_{\sigma\in\Sigma_0}\big|\;\sum_\sigma\|T_\sigma\|_{trace}<\infty,\;T_\sigma\in\mathcal{N}\big(\Ha_{\sigma,\Gr}\big)\right\},
\end{align*}
where $\mathcal{N}\big(\Ha_{\sigma,\Gr}\big)$ is the set of nuclear operators on $\Ha_{\sigma,\Gr}$, and thus
\begin{gather*}
\Lainf_{\Sigma',\Gr}\cong\left(\La_{\Sigma_0,\Gr}^1\right)^*\cong\left(\La_{\Sigma,\Gr}^0\right)^{**}.
\end{gather*}
The statement of the lemma follows by Lemma \ref{universalmap}.\\

\qed

\begin{lem}\label{KernBildPi}
The kernel of $\widehat{\pi}$ is a two sided $\Gr$-invariant, weakly closed ideal of $(\im\beta)''$ and hence we obtain the following isomorphisms:
\begin{gather*}
\ker\widehat{\pi}\cong\Lainf_{\Sigma^0,\Gr},\\
\im\widehat{\pi}\cong\Lainf_{\Sigma',\Gr},
\end{gather*}
where $\Sigma':=\Sigma\backslash\Sigma^0$ for some subset of $\Sigma$.
\end{lem}

\proof
The kernel $\ker\widehat{\pi}$ is a two-sided ideal of $(\im\beta)''$. If $p_\sigma$ are projections of $(\im\beta)''$  onto $\L(\Ha_{\sigma,\Gr})$, then $p_\sigma\in\Za_{(\im\beta)''}$ and therefore $p_\sigma\cdot\ker\widehat{\pi}$ is an two-sided ideal in $\L(\Ha_{\sigma,\Gr})$. Since $\L(\Ha_{\sigma,\Gr})$ is, due to the finite-dimensionality of $\Ha_{\sigma,\Gr}$, simple, the ideal is either trivial or the whole algebra $\L(\Ha_{\sigma,\Gr})$. Thus we obtain:
\begin{gather*}
\sum_{\sigma\in\Sigma}p_\sigma\cdot\ker\widehat{\pi}=\ker\widehat{\pi}=\sum_{\sigma\in\Sigma_0}p_\sigma\cdot\ker\widehat{\pi},
\end{gather*}
where $\Sigma_0:=\{\sigma\in\Sigma|\;p_\sigma\cdot\ker\widehat{\pi}\neq 0\}$. Evidently, the image has the structure given in the theorem.

\qed

\proofof Theorem \ref{compactinner}: 
\begin{itemize}
\item $\beta(g)=\big(\tilde{\beta}\circ\lambda\big)(g)$, $g\in\Gr$, is valid by definition.
\item The involution is transformed properly:
\begin{align*}
\langle\tilde{\beta}(\mu),h\rangle&=\langle f,(\tilde{\beta}(\mu))(h)\rangle\\
&=\int\int f(g_1^{-1}g_2)d\mu(g_1)\cdot\overline{h(g_2)}g_2\\
&=\int\int f(g_3)d\mu(g_1)\overline{h(g_1g_3)}dg_3\\
&=\int\int f(g_3)\overline{h(g_1g_3)}d\mu(g_1)dg_3\\
&=\int f(g_3)\overline{\overline{\int\overline{h(g_1g_3)}d\mu(g_1)}}dg_3\\
&=\langle f,(\tilde{\beta}(\mu^*))(h)\rangle.
\end{align*}
\item $\pi(g)=\tilde{\pi}_M\circ\lambda(g)$, $g\in\Gr$, by definition.
\item $\pi(g)=\widehat{\pi}\circ\beta(g)=\widehat{\pi}_2\circ\widehat{\pi}_1\circ\beta(g)$, $g\in\Gr$, follows by Lemma \ref{pihat} and Lemma \ref{imbeta}.
\item $\tilde{\pi}_M=\tilde{\pi}_{L^1}\circ\iota_1$ is standard.
\item The homomorphism $\tilde{\pi}_{L^1}$ is ensured by the Peter-Weyl theorem \ref{PeterWeyl}.
\end{itemize}
\qed

\begin{theo}\label{StrukturInner}
Let $\Ma$, $\Gr$ and $\widehat{\pi}_2$ be defined as aforementioned and denote the center of $\Ma$ by $\Za\equiv\Za_\Ma$, then one has the following identity:
\begin{gather*}
\left(\Ma^\Gr\right)\rc\equiv(\im\widehat{\pi}_2\cup\Za)\rc\rc.
\end{gather*}

\end{theo}

\proof
We will show first the inclusion ``$\supseteq$''. Obviously, the center of $\Ma$ lies in $\big(\Ma^\Gr\big)\rc$. First we prove $\im\widehat{\pi}_2\subseteq\left(\Ma^\Gr\right)\rc$. Thanks to the commutative property of the diagram it suffices to show $\tilde{\pi}_M(\Gr)\subseteq\left(\Ma^\Gr\right)\rc$. For $A\in\Ma^\Gr$ and $B\in\tilde{\pi}_M(\Gr)$ we obtain 
\begin{align*}
BA&=\int\pi(g)d\mu(g)A\\
&=\int\pi(g)Ad\mu(g)\\
&=\int A\pi(g)d\mu(g)\\
&=AB,
\end{align*}
i.e. $B\in\left(\Ma^\Gr\right)\rc$ and, since due to Lemma \ref{vonNeumannsubalg} $\Ma^\Gr$ is closed, the completion of $\im\widehat{\pi}_2\cup\Za$ is contained in $\big(\Ma^\Gr)\rc$.  \\
On the other hand one has $\left[\tilde{\pi}_M\big(M(\Gr)\big)\right]'\cap\Ma\subseteq\Ma^\Gr$, because for $A\in\left[\tilde{\pi}_M\big(M(\Gr)\big)\right]'\cap\Ma$ and $B\in\tilde{\pi}_M\big(M(\Gr)\big)$ one has
\begin{align*}
A\int\pi(g)d\mu(g)=\int\pi(g)d\mu(g)A
\end{align*}
for all measures $\mu\in M(\Gr)$. We obtain in particular for one-point measures:
\begin{gather*}
A\pi(\delta_g)=\pi(\delta_g)A,\\
\Longleftrightarrow\quad\pi(\delta_g)A\pi(\delta_g)^{-1}=A,\quad\forall g\in\Gr,
\end{gather*}
i.e. $A\in\Ma^\Gr$. Finally, we derive:
\begin{align*}
&\hspace{-1,5cm}&\big[\tilde{\pi}_M\big(M(\Gr)\big)\big]'\cap\Ma&\subseteq\Ma^\Gr\\
&\Longleftrightarrow\hspace{-1,5cm}&\big[\tilde{\pi}_M\big(M(\Gr)\big)\big]'\cap\Za'\cap\Ma&\subseteq\Ma^\Gr\\
&\Longleftrightarrow\hspace{-1,5cm}&\big[\tilde{\pi}_M\big(M(\Gr)\big)\cup\Za\big]'\cap\Ma&\subseteq\Ma^\Gr\\
&\Longleftrightarrow\hspace{-1,5cm}&\big[\tilde{\pi}_M\big(M(\Gr)\big)\cup\Za\big]\rc&\subseteq\Ma^\Gr\\
&\Longleftrightarrow\hspace{-1,5cm}&\big[\tilde{\pi}_M\big(M(\Gr)\big)\cup\Za\big]\rc\rc&\supseteq\big(\Ma^\Gr\big)\rc.
\end{align*}
\qed

\begin{lem}\label{StrukturInnerFactor}
Let $\Gr$ and $\widehat{\pi}_2$ be defined as aforementioned and $\Ma$ a factor, then one has $\im\widehat{\pi}_2\cong\big(\Ma^\Gr\big)\rc$.
\end{lem}

\proof
We may apply Lemma \ref{StrukturInner} to the following identity:
\begin{align*}
(\im\widehat{\pi}_2)\rc\rc&=\big[(\im\widehat{\pi}_2)'\cap\Ma\big]'\cap\Ma\\
&=\big[(\im\widehat{\pi}_2)''\cup\Ma'\big]\cap\Ma\\
&=\big[\im\widehat{\pi}_2\cap\Ma\big]\cup\big[\Ma'\cap\Ma\big]\\
&=\im\widehat{\pi}_2.
\end{align*}
\qed
 
\begin{defi}
Let $\tilde{\Sigma}=(\Ha_{\sigma,\Gr}, \sigma)$ a subset of the set of all irreducible representations, then a compact group $\Gr$ is called full, if $\Gr\cong\prod_{\sigma\in\tilde{\Sigma}}U(\Ha_{\sigma,\Gr})$. A subgroup $\Hr$ of $\Gr$ is called full, if $\Hr\cong\prod_{\tau\in} U\big(\Ja_\tau\big)$, where $\Ja_\tau$ is a subspace of $\Ha_{\sigma,\Gr}$ for some $\sigma$ and the isomorphism is given by the restriction of the former isomorphism.
\end{defi}

\begin{lem}
Let 
\begin{gather*}
\Ka_\sigma:=\left\{\xi\in\Ha_{\sigma,\Gr}|\;\xi=\int_\Hr\sigma(h)\eta d_\Hr\text{ for some }\eta\in\Ha_{\sigma,\Gr}\right\},\\
\La_\sigma:=\left\{T\in\L\big(\Ha_{\sigma,\Gr}\big)|\;T=\int_\Hr\int_\Hr\sigma(h)S(\sigma(h'))\xi d_\Hr^2\text{ for some }S\in\L(\Ha_{\sigma,\Gr})\right\},
\end{gather*} 
then $\La_\sigma$ is isomorphic to $\L(\Ka_\sigma)$.
\end{lem}

\proof
The isomorphism is given by 
\begin{align*}
\varphi:\L_\sigma&\longrightarrow\L(\Ka_\sigma)\\
T&\mapsto\varphi(T)\\
\varphi(T)\xi&:=\int_\Hr\int_\Hr\sigma(h)T(\sigma(h'))\xi d_\Hr^2.
\end{align*}
\qed

\begin{defi}
Let the set of unitaries $\U(\Ha_{\sigma,\Gr}\big)$ be defined as above, then we call $\Ge:=\prod_{\sigma\in\Sigma}\U(\Ha_{\sigma,\Gr}\big)$ the enveloping of the group $\Gr$.  
\end{defi}

\begin{defi}
We call $\Ge_{\widehat{\pi}_2}:=\prod_{\sigma\in\Sigma'}\U\big(\widehat{\pi}_2(\Ha_{\sigma,\Gr})\big)\subset(\Ma^\Gr)\rc$ the enveloping of the group $\Gr$ with respect to $\widehat{\pi}_2$.  
\end{defi}

\begin{lem}
The enveloping group $\Ge$ of a compact group $\Gr$ is compact.
\end{lem}

\proof
$\Ge$ is as a product of compact groups by Tychonov's theorem compact, too. The continuity of the multiplication and the inverse operation in each components $\sigma\in\Sigma'$ implies that $\Ge$ is a compact group.
\qed

\begin{rems}
\begin{itemize}
\item[(i)] Locally compactness of $\Gr$ does not necessarily imply locally compactness for its enveloping group.
\item[(ii)] $\pi(\Gr)$ is a subgroup of $\Ge$.
\item[(iii)] $\widehat{\pi}_2\circ\sigma$, where $\sigma\in\Sigma'$ is representation, can be extended to $\Ge$ and is obviously irreducible.
\item[(iv)] $\Ge$ does not depend on the choice of $\Ma$. 
\end{itemize}
\end{rems}

\begin{defi}
Two subgroups $\Hr_1$ and $\Hr_2$ of $\Ge$ are said to be equivalent, if $\text{span}\sigma({\Hr_1})=\text{span}\sigma({\Hr_2})$ for all $\sigma\in\Sigma'$.
\end{defi}
Obviously, $\Gr$ and its enveloping $\Ge$ are equivalent. The commutative diagram of Theorem \ref{compactinner} has the same structure for both groups $\Gr$ and $\Ge$.\\
In what follows we will denote by $\Ka_\sigma$, $\sigma\in\Sigma'$, Hilbert subspaces of $\Ha_{\sigma,\Gr}$ satisfying for all $\sigma_1,\sigma_2\in\Sigma'$:
\begin{gather}\label{ksigma}
\begin{split}
\Ka_{\sigma_1}\oplus\Ka_{\sigma_2}&\subseteq\Ka_{\sigma_1\oplus\sigma_2},\\
\Ka_{\sigma_1}\otimes\Ka_{\sigma_2}&\subseteq\Ka_{\sigma_1\otimes\sigma_2},\\
\overline{\Ka_\sigma}&=\Ka_{\overline{\sigma}}.
\end{split}
\end{gather}
Furthermore, we set:
\begin{gather*}
\Ha_{\Sigma,\Hr}^2:=l^2\left(\bigoplus_{n=1}^{n_\sigma}\Ha_{\sigma,\Gr}/\Ka_\sigma\right):=\left\{\big((\xi_\sigma^1,...,\xi_\sigma^{n_\sigma})\big)_{\sigma\in\Sigma}\big|\;\sum_{k=1}^{n_\sigma}\|\xi_\sigma^k\|^2<\infty\right\}.
\end{gather*}

\begin{theo}\label{compactinnergal}
If $\pi(\Gr)$ is proper, then there exists an injective map from the set of all closed subgroups $\Hr$ of $\Gr$ and $\Gr$-invariant von Neumann subalgebras $\Aa=\Aa''$ of $\Ma$, $\Ma^\Gr\subset\Aa$, (which are constructed by $\Ka_\sigma$) according to the following diagram:
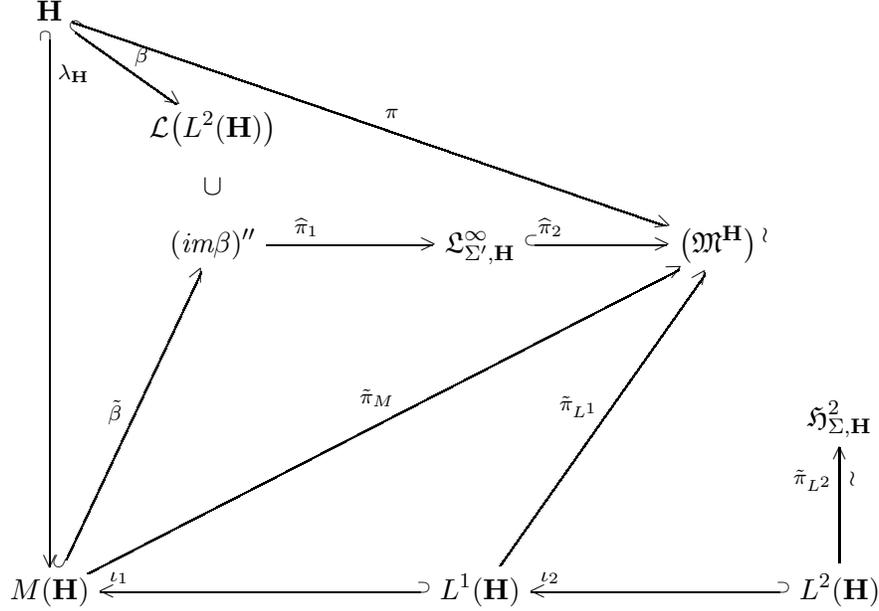
\begin{diagram}
\Hr       &     &      &     &&& &          &&&&&\\
\dInto(0,12)_{\lambda_\Hr}&\rdInto(2,2)^\beta\rdTo(10,4)^{\pi}&&&&&\\ 
         &          &\L\big(L^2(\Hr)\big) \\
         &          &\cup \\
         &          &(\im\beta)''    &\rTo^{\widehat{\pi}_1} &&&\Lainf_{\Sigma',\Hr}&\rInto^{\widehat{\pi}_2}      &&&\left(\Ma^\Hr\right)\rc&&\\
         &\ruInto(2,6)^{\tilde{\beta}}&      &     &&& &    &&\ruTo(10,6)^{\tilde{\pi}_M}\ruTo(4,6)^{\tilde{\pi}_{L^1}}&&&\\
\\ 
         &          &      &     &&& &          &&&&&\Ha_{\Sigma,\Hr}^2\\
         &          &      &     &&& &          &&&&&\uTo(0,4)_\wr^{\tilde{\pi}_{L^2}}\\
\\
M(\Hr)    &\lInto^{\iota_1}     &      &     &&&L^1(\Hr)&\lInto^{\iota_2}      &&&&&L^2(\Hr)
\end{diagram}
\end{theo}

\proof
Let $\Hr$ be a subgroup of $\Gr$, then there exists a family of $\Ka_\sigma$ satisfying the aforementioned conditions which are the fix-point subspaces of $\Ha_\sigma$ by $\sigma(\Hr)$. Hence, we obtain
\begin{gather*}
\Ha_{\sigma,\Hr}=\Ha_{\sigma,\Gr}/\Ka_\sigma=\Ka_\sigma^\perp
\end{gather*}
and therefore
\begin{gather*}
\Ha_{\Sigma,\Hr}^2\cong L^2(\Hr).
\end{gather*} 
On the other hand, let $\Ka_\sigma$, $\sigma\in\Sigma'$, be given. Since $\pi(\Gr)$ is proper, one has thanks to Lemma \ref{propersigma} $\Sigma=\Sigma'$ since $\tilde{\pi}_M\big|_{L^2(\Gr)}$ is injective, and therefore one can uniquely extend the mapping $\sigma\mapsto\Ka_\sigma$, $\sigma\in\irrep(\Gr)$, to the one on $\rep(\Gr)$ fulfilling the properties \ref{ksigma}. Thus there exists a closed subgroup $\Hr\subset\Gr$ such that
\begin{gather*}
\Ka_\sigma=\big\{\xi\in\Ha_\sigma|\;\sigma(h)\xi=\xi,\;h\in\Hr\big\},\quad\forall\sigma\in\Sigma.
\end{gather*}
It remains to be proved that each family $\Ka_\sigma$, $\sigma\in\Sigma$, defines a family $\Ka_\sigma$, $\sigma\in\rep(\Gr)$. We call two of such families of subspaces $\{\Ka_\sigma|\;\sigma\in\rep(\Gr)\}$ and $\{\tilde{\Ka}_\sigma|\;\sigma\in\rep(\Gr)\}$ equivalent, if $\Ka_\sigma=\tilde{\Ka}_\sigma$ for all irreducible $\sigma\in\Sigma$. If $\Ka_\sigma$, $\sigma\in\Sigma$ is given, then we define $\Ka_\sigma$, $\sigma\in\rep(\Gr)$, to be the maximal elements in the equivalence class generated by $\{\Ka_\sigma|\;\sigma\in\Sigma\}$. Thanks to the lemma of Zorn, the existence of these elements are ensured and they are unique.
\qed

\begin{rem}
If $\pi(\Gr)$ is not proper, i.e. $\Sigma\neq\Sigma'$, then different subspaces $\Ka_{\sigma_1}$ and $\Ka_{\sigma_2}$ may be assigned to the same subgroup $\Hr$. 
\end{rem}



\begin{theo}
Let $\Ma$ be a factor and $\Gr$ a full group, then there exists a bijective map between the set of full subgroups of $\Gr$ and the set of subalgebras $\Ba$ of $\big(\Ma^\Gr\big)\rc$ with $\Ba\rc\rc=\Ba$ and hence to the set of subalgebras of $\Ma$ containing $\Ma^\Gr$.
\end{theo}

\proof
Let $\Hr$ be a full subgroup of $\Gr$, then $\pi(\Hr)\subseteq\big(\Ma^\Gr\big)\rc$ and we define the map 
\begin{align*}
p_\pi:\Fu(\Gr)&\longrightarrow\Alg\big(\big(\Ma^\Gr\big)\rc\big)\\
\Hr&\mapsto p_\pi(\Hr):=(\Ma^\Hr\big)\rc
\end{align*}
Property $\big(\Ma^\Hr\big)\rc\rc=\Ma^\Hr$ is given by Lemma \ref{vonNeumannsubalg}.
We state for the inverse map: 
\begin{align*}
p_\pi^{-1}(\Ba)=\pi^{-1}(\Ba),\quad\forall\Ba\in\Alg\big(\big(\Ma^\Gr\big)\rc\big).
\end{align*}
It remains to be shown that $p_\pi^{-1}(\Ba)$ is a full subgroup of $\Gr$. Let us assume $\Ba\rc\rc=\Ba\subset(\Ma^\Gr\big)\rc$. Since $\Ma$ is a factor, we may write due to Lemma \ref{StrukturInnerFactor} $(\Ma^\Gr\big)\rc\cong\im\widehat{\pi}_2$. We define $\Hr:=p_\pi^{-1}\big(\Ua(\Ba)\big)\cong\Ua(\Ba)$. Since $\Hr$ is closed and as the subset of a compact group it is also compact and $\big(\Ma^\Hr\big)\rc=\Ba$. Lemma \ref{StrukturInner} ensures the isomorphism $\big(\Ma^\Hr\big)\rc\cong\Lainf_{\Sigma',\Hr}$ and therefore we obtain $\Hr\cong\Ua(\Ba)\cong\prod_{\sigma\in\Sigma_\Hr'}\Ua\big(\Ha_{\sigma,\Hr}\big)$, i.e. $\Hr$ is a full subgroup of $\Gr$.

\qed

\begin{rem}
The subalgebras of $\Lainf_{\Sigma',\Hr}$ are isomorphic to algebras of the following structure:
\begin{gather*}
l^\infty\left(\bigoplus_{\sigma\in\Sigma''}\L\big(\Ha_{\sigma,\Gr}\big/\Ja_\sigma\big)\right), \quad\Ja_\sigma\subset\Ha_{\sigma,\Gr},
\end{gather*}
for some $\Sigma''\subset\Sigma'$.
\end{rem}

\begin{theo}
There is a one-to-one correspondence between equivalence classes $\widehat{\Hr}$ of closed subgroups $\Hr$ of $\Ge$ and von Neumann subalgebras $\Aa$ of $\Ma$, such that $\Aa\rc\rc=\Aa$ and $\Aa$ is the fix-point algebra under each subgroup of the equivalence class of $\Hr$.
\end{theo}

\proof
If $\Hr$ is given, then the fix-point algebra fulfills thanks to Lemma \ref{vonNeumannsubalg} $(\Ma^\Hr)\rc\rc=\Ma^\Hr$. Let $\Hr_1$ be in the same equivalence class of $\Hr$, then $\text{span}\sigma({\Hr_1})=\text{span}\sigma({\Hr})$ for all $\sigma\in\Sigma'$ and hence $(\Ma^{\Hr_1})\rc=(\Ma^\Hr)\rc$, i.e. $\Ma^{\Hr_1}=\Ma^\Hr$.\\
Contrariwise, let $\Aa=\Aa\rc\rc$ be a von Neumann subalgebra with $\Ma^\Gr\subset\Aa\subset\Ma$, i.e. $(\Ma^\Gr)\rc\supset\Aa\rc$, then, since $(\Ma^{\Ge})\rc=(\Ma^\Gr)\rc$ and $(\Ma^\Gr)\rc\cong\prod_{\sigma\in\Sigma'}\L\big(\Ha_{\sigma,\Gr}\big)$, $\Aa$ is as a subalgebra isomorphic to $\prod_{\sigma\in\Sigma'}\L\big(\Ha_{\sigma,\Gr}\big/\Ja_\sigma\big)$. The subspaces $\Ja_\sigma$, $\sigma\in\rep(\Ge)$, constructed for $\Ge$ as in the proof of Theorem \ref{compactinnergal} obviously satisfy the conditions \ref{ksigma}. We define $\Hr$ as the product of the group of unitary operators on $\Ha_{\sigma,\Gr}$, which are the identity on $\Ka_\sigma$. Therefore, $\Hr$ is a closed subgroup of $\Ge$ with $\Aa=\Ma^\Hr$.
\qed

\subsection{Groups of Spatial Automorphisms}

\begin{lem}
Let $\Ma^\Gr:=\{A\in\Ma|\;\alpha^g(A)=A,\;\forall g\in\Gr\}$, then the image of $\pi_\alpha$ is in $(\Ma^\Gr)'$.
\end{lem}

\proof
By definition one obtains $\pi_\alpha(g)A=A\pi_\alpha(g)$ for all $A\in\Ma^\Gr$.
\qed

\begin{rem}
The center $\Za$ of $\Ma$ lies in $\Ma^\Gr\cap(\Ma^\Gr)\rc$ and $\pi(\Gr)\cap\Za=\S^1$.
\end{rem}

\begin{lem}
There exists a homomorphism $\Gr\longrightarrow\Ua(\Ha)$ such that 
\begin{gather*}
\alpha^g(A)=\pi_gA\pi_g^{-1},\qquad \forall g\in\Gr.
\end{gather*}
\end{lem}

\proof
Since each $g\in\Gr$ acts as an inner automorphism on $\LH$ we apply Lemma \ref{mappingpi}.
\qed

\begin{theo}
There exists a set $\Sigma':=\Sigma\backslash\Sigma^0,\quad\Sigma^0:=\ker(\hat{\pi}_M\circ\iota_1\circ\iota_2)$, such that the following diagram is commutative:

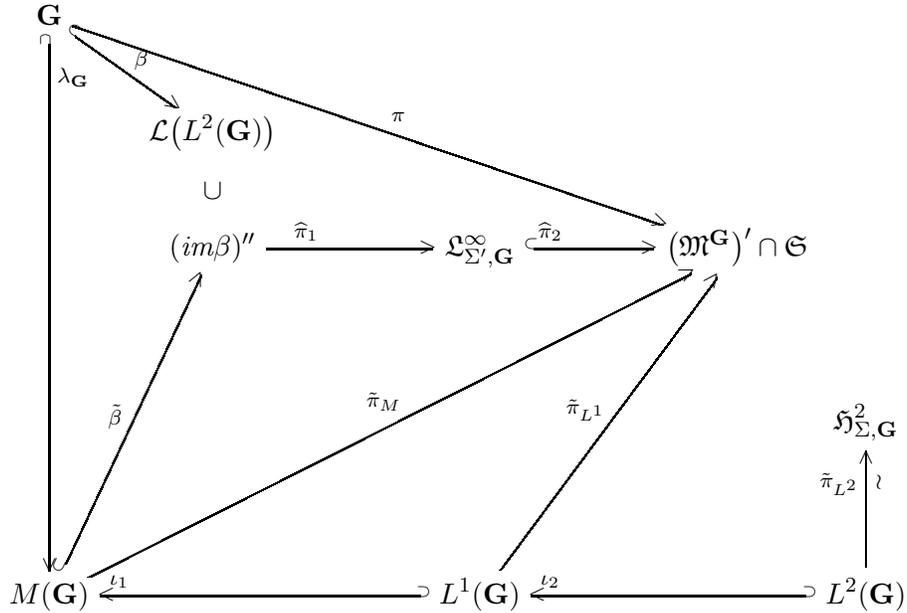
\begin{diagram}
\Gr       &     &      &     &&& &          &&&&&\\
\dInto(0,12)_{\lambda_\Gr}&\rdInto(2,2)^\beta\rdTo(10,4)^{\pi}&&&&&\\ 
         &          &\L\big(L^2(\Gr)\big) \\
         &          &\cup \\
         &          &(\im\beta)''    &\rTo^{\widehat{\pi}_1} &&&\Lainf_{\Sigma',\Gr}&\rInto^{\widehat{\pi}_2}      &&&\big(\Ma^\Gr\big)'\cap\Sa&&\\
         &\ruInto(2,6)^{\tilde{\beta}}&      &     &&& &    &&\ruTo(10,6)^{\tilde{\pi}_M}\ruTo(4,6)^{\tilde{\pi}_{L^1}}&&&\\
\\ 
         &          &      &     &&& &          &&&&&\Ha_{\Sigma,\Gr}^2\\
         &          &      &     &&& &          &&&&&\uTo(0,4)_\wr^{\tilde{\pi}_{L^2}}\\
\\
M(\Gr)    &\lInto^{\iota_1}     &      &     &&&L^1(\Gr)&\lInto^{\iota_2}      &&&&&L^2(\Gr)
\end{diagram}
where $\Sa:=\{S\in\LH|\;S\Ma=\Ma S\}$.
\end{theo}

\proof
Since $\Gr$ acts as an inner automorphism group on $\L(\Ha)$, one may refer to the proof of the inner case, Theorem \ref{compactinner}.
\qed

\begin{lem}\label{structurespatial1}
The following statements hold:
\begin{itemize}
\item[(i)] $\big[\big(\L(\Ha)^\Gr\big)'\cap\Sa\big]\rc=\Ma^\Gr\cup\Sa'$;
\item[(ii)] $\big(\L(\Ha)^\Gr\big)'=\big(\im\pi\big)''$;
\item[(iii)] $\Ma^\Gr=(\im\pi)\rc$;
\item[(iv)] $\left(\Ma^\Gr\right)'=\Big(\big(\L(\Ha)^\Gr\big)'\cup\Ma'\Big)''$;
\item[(v)] $\im\pi\cap\Ma'=\C\1$.
\end{itemize}
\end{lem}

\proof
\begin{itemize}
\item[(i)] Due to $\Ma'\subseteq\Sa$, therefore $\Ma=\Ma''\supseteq\Sa'$, and $\big[\big(\L(\Ha)\big)^\Gr\big]''=\big(\L(\Ha)\big)^\Gr$ one obtains $\big[\big(\L(\Ha)^\Gr\big)'\cap\Sa\big]\rc=\big[\L(\Ha)^\Gr\cup\Sa'\big]\cap\Ma=\Ma^\Gr\cup\Sa'$.
\item[(ii)] This relation follows from the fact that $\pi$ maps $\Gr$ onto the unitary elements of $\Ma$.
\item[(iii)] This is a consequence of $(ii)$.
\item[(iv)] $\left(\Ma^\Gr\right)'=\Big(\big(\L(\Ha)^\Gr\cap\Ma\big)'\big)''=\Big(\big(\L(\Ha)^\Gr\big)'\cup\Ma'\Big)$;
\item[(v)] This is a direct consequence of $\pi$ mapping $\Gr$ on the unitary elements of $\Ma$.
\end{itemize}

\qed

\begin{rem}
Since $\Sa'$ lies in $\Za$, statement $(i)$ of the lemma reduces for factors $\Ma$ to $\Big(\big[\L(\Ha)^\Gr\big]'\cap\Sa\Big)\rc=\Ma^\Gr$. 
\end{rem}

\begin{theo}\label{structurespatial2}
One has the following identity:
\begin{gather*}
\left(\Ma^\Gr\right)'=(\im\widehat{\pi}_2\cup\Za)\rc{}'.
\end{gather*}
\end{theo}

\proof
First we prove $\im\widehat{\pi}_2\subseteq\left(\Ma^\Gr\right)'$. Thanks to the commutative property of the diagram it suffices to show $\tilde{\pi}_M(\Gr)\subseteq\left(\Ma^\Gr\right)'$. Let $A\in\Ma^\Gr$ and $B\in\tilde{\pi}_M(\Gr)$, i.e.
\begin{align*}
BA&=\int\pi(g)d\mu(g)A\\
&=\int\pi(g)Ad\mu(g)\\
&=\int A\pi(g)d\mu(g)\\
&=AB.
\end{align*}
Hence $B\in\left(\Ma^\Gr\right)'$ and finally we conclude:
\begin{align*}
&\big(\Ma^\Gr\big)'\supseteq\im\widehat{\pi}_2\cup\Za\\
\Longleftrightarrow\quad&\Ma^\Gr\subseteq\big(\im\widehat{\pi}_2\cup\Za\big)'\\
\Longleftrightarrow\quad&\Ma^\Gr\subseteq\big(\im\widehat{\pi}_2\cup\Za\big)\rc\\
\Longleftrightarrow\quad&\big(\Ma^\Gr\big)'\supseteq\big(\im\widehat{\pi}_2\cup\Za\big)\rc{}'.
\end{align*}
On the other hand one has $\Ma^\Gr\supseteq\big(\tilde{\pi}_M(\Gr)\big)'\cap\,\Ma$, because for $A\in\big(\tilde{\pi}_M(\Gr)\big)'\cap\Ma$ and $B\in\tilde{\pi}_M(\Gr)$ one has
\begin{gather*}
AB=A\int\pi(g)d\mu(g)=\int\pi(g)d\mu(g)A=BA
\end{gather*}
for all measures $\mu\in M(\Gr)$. We obtain in particular for one-point measures:
\begin{gather*}
A\pi(g_0)=\pi(g_0)A,\\
\Longleftrightarrow\quad\pi(g_0)A\pi(g_0)^{-1}=A,\quad\forall g_0\in\Gr,
\end{gather*}
i.e. $A\in\Ma^\Gr$. Finally, we derive:
\begin{align*}
&&\Ma^\Gr&\supseteq\big(\tilde{\pi}_M(\Gr)\big)'\cap\Ma\\
\Longleftrightarrow&&\Ma^\Gr&\supseteq\big(\tilde{\pi}_M(\Gr)\big)'\cap\Za'\cap\Ma\\
\Longleftrightarrow&&\Ma^\Gr&\supseteq\big(\tilde{\pi}_M(\Gr)\cup\Za\big)'\cap\Ma\\
\Longleftrightarrow&&\Ma^\Gr&\supseteq\big(\tilde{\pi}_M(\Gr)\cup\Za\big)\rc\\
\Longleftrightarrow&&\big(\Ma^\Gr\big)'&\subseteq\big(\tilde{\pi}_M(\Gr)\cup\Za\big)\rc{}'.
\end{align*}
\qed

\begin{lem}
Let $J$ be the modular conjugation, then one has the following identity:
\begin{gather}
\Big[\big(\Ma'\big)^\Gr\Big]'=J\big(\Ma^\Gr\big)'J.
\end{gather}
\end{lem}

\proof
First of all, $\Ma'$ is invariant under the action of the automorphism group $\alpha^\Gr$, since for each $A\in\Ma'$ there exists an element $B\in\Ma$ with $A=JBJ$, such that
\begin{align*}
U_gAU_g^{-1}=U_gJBJU_g^{-1}=JU_gBU_g^{-1}J\in\Ma'
\end{align*} 
for all $g\in\Gr$. The commutation of the operators is ensured by Theorem \ref{standardimplementation}. \\
The statement of the lemma is equivalent to
\begin{gather*}
\big(\Ma'\big)^\Gr=\Big[J\big(\Ma^\Gr\big)'J\Big]'=J\Ma^\Gr J.
\end{gather*}
If $A\in J\Ma^\Gr J$, then there exists an element $B$ of $\Ma^\Gr$ with $A=JBJ$ and one obtains:
\begin{align*}
U_gAU_g^{-1}=U_gJBJU_g^{-1}=JU_gBU_g^{-1}J=JBJ=A,\quad\forall g\in\Gr.
\end{align*}
Thus $A$ lies also in $\big(\Ma'\big)^\Gr$. Contrariwise, let us suppose $A\in\big(\Ma'\big)^\Gr$, then there is an element $B\in\Ma$ with $A=JBJ$ and one concludes:
\begin{gather*}
A=U_gAU_g^{-1}=U_gJBJU_g^{-1}\\
\Longleftrightarrow\quad JBJ=JU_gAU_g^{-1}J,\qquad\forall g\in\Gr.
\end{gather*}
Therefore $B\in\Ma^\Gr$ and $A$ is contained in $J\Ma^\Gr J$.
\qed



\begin{theo}\label{spatial2}
If $\pi(\Gr)$ is proper, then there exists an injective map from the set of all closed subgroups $\Hr$ of $\Gr$ and $\Gr$-invariant von Neumann subalgebras $\Aa=\Aa'{}\rc$ of $\Ma$, $\Ma^\Gr\subset\Aa$, (which are constructed by $\Ka_\sigma$) according to the following diagram:
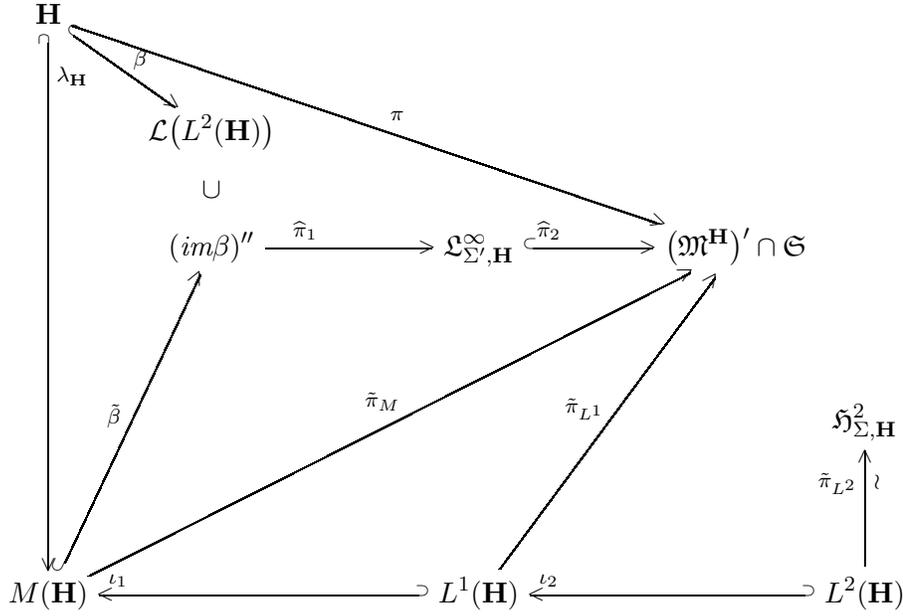
\begin{diagram}
\Hr       &     &      &     &&& &          &&&&&\\
\dInto(0,12)_{\lambda_\Hr}&\rdInto(2,2)^\beta\rdTo(10,4)^{\pi}&&&&&\\ 
         &          &\L\big(L^2(\Hr)\big) \\
         &          &\cup \\
         &          &(\im\beta)''    &\rTo^{\widehat{\pi}_1} &&&\Lainf_{\Sigma',\Hr}&\rInto^{\widehat{\pi}_2}      &&&\big(\Ma^\Hr\big)'\cap\Sa&&\\
         &\ruInto(2,6)^{\tilde{\beta}}&      &     &&& &    &&\ruTo(10,6)^{\tilde{\pi}_M}\ruTo(4,6)^{\tilde{\pi}_{L^1}}&&&\\
\\ 
         &          &      &     &&& &          &&&&&\Ha_{\Sigma,\Hr}^2\\
         &          &      &     &&& &          &&&&&\uTo(0,4)_\wr^{\tilde{\pi}_{L^2}}\\
\\
M(\Hr)    &\lInto^{\iota_1}     &      &     &&&L^1(\Hr)&\lInto^{\iota_2}      &&&&&L^2(\Hr)
\end{diagram}
where the $\Ka_\sigma$ are subspaces of $\Ha_\sigma$ satisfying the conditions \ref{ksigma}.
\end{theo}

\proof
Thanks to properness of $\pi$, we get as for the inner case $(\im\pi_\Hr)''\neq(\im\pi_\Gr)''$ via the construction of the $\Ka_\sigma$. By Lemma \ref{structurespatial1} this implies $\big[\L(\Ha)^\Hr\big]'\neq\big[\L(\Ha)^\Gr\big]'$. Furthermore due to Lemma \ref{structurespatial1}, one has $\left(\Ma^\Hr\right)'=\Big(\big(\L(\Ha)^\Hr\big)'\cup\Ma'\Big)''=\Ma'\rtimes_\alpha\Hr$ acting on $\Ha$. But $\Ma'\rtimes_\alpha\Hr\neq\Ma'\rtimes_\alpha\Gr$ holds since $\C\rtimes_\alpha\Hr\neq\C\rtimes_\alpha\Gr$.

\qed


\begin{theo}\label{spatialfactor}
Let $\Ma$ be a factor and $\Gr$ a full group, then there exists a bijective map between the set of full subgroups of $\Gr$ and the set of subalgebras $\Ba$ of $\big(\Ma^\Gr\big)'\cap\Sa$ with $\Ba\rc{}'\cap\Sa=\Ba$ and hence to the set of subalgebras of $\Ma$ containing $\Ma^\Gr$.
\end{theo}

\proof
Let $\Hr$ be a full subgroup of $\Gr$, then $\pi(\Hr)\subseteq\big(\Ma^\Gr\big)'\cap\Sa$ and we define the map 
\begin{align*}
p_\pi:\Fu(\Gr)&\longrightarrow\Alg\big(\big(\Ma^\Gr\big)'\cap\Sa\big)\\
\Hr&\mapsto p_\pi(\Hr):=(\Ma^\Hr\big)'\cap\Sa
\end{align*}
We state for the inverse map: 
\begin{align*}
p_\pi^{-1}(\Ba)=\pi^{-1}(\Ba),\quad\forall\Ba\in\Alg\big(\big(\Ma^\Gr\big)'\cap\Sa\big).
\end{align*}
It remains to be shown that $p_\pi^{-1}(\Ba)$ is a full subgroup of $\Gr$. \\
Let $\Ba\rc{}'\cap\Sa=\Ba\subset(\Ma^\Gr\big)'\cap\Sa$. Since $\Ma$ is a factor, we may write due to Lemma \ref{structurespatial1} and Theorem \ref{structurespatial2} $(\Ma^\Gr\big)'\cap\Sa\cong\im\widehat{\pi}_2$. We define $\Hr:=p_\pi^{-1}\big(\Ua(\Ba)\big)\cong\Ua(\Ba)$. Since $\Hr$ is closed and as the subset of a compact group it is also compact and $\big(\Ma^\Hr\big)'\cap\Sa=\Ba$. Because $\big(\Ma^\Hr\big)'\cap\Sa\cong l^\infty\left(\bigoplus_{\sigma\in\Sigma_\Hr'}\L(\Ha_{\sigma,\Hr})\right)$ we obtain $\Hr\cong\Ua(\Ba)\cong\prod_{\sigma\in\Sigma_\Hr'}\Ua\big(\Ha_{\sigma,\Hr}\big)$, i.e. $\Hr$ is a full subgroup of $\Gr$.

\qed

\begin{theo}\label{spatial-one-to-one}
There is an one-to-one correspondence between equivalence classes $\H$ of closed subgroups $\widehat{\Hr}$ of $\Ge$ and von Neumann subalgebras $\Aa$ of $\Ma$, such that $\Aa$ is the fix-point algebra under each subgroup of the equivalence class of $\Hr$..
\end{theo}

\proof
Since $\Gr$ acts as an inner automorphism group on $\L(\Ha)$, we refer to the proof of the inner case.
\qed

\subsection{General Case}

\begin{lem}\label{general1}
The map 
\begin{gather*}
\tilde{\pi}_M:M(\Gr)\longrightarrow\big((\Ma\rtimes_\alpha\Gr)^\Gr\big)'\cap\Sa,
\end{gather*}
where $\Sa:=\big\{S\in\L(L^2(\Gr,\Ha)|\;S(\Ma\rtimes_\alpha\Gr)=(\Ma\rtimes_\alpha\Gr) S\big\}$, is injective.
\end{lem}

\proof
The map $\tilde{\pi}_M$ is injective since for $0\neq\mu\in M(\Gr)$, as defined by
\begin{gather*}
\tilde{\pi}_M(\mu)=\int_\Gr U_gd\mu(g),
\end{gather*} 
there exists at least one $f\in L^2(\Gr)$ such that 
\begin{gather*}
\tilde{\pi}_M(\mu)(\eta\cdot f)=\eta(\mu\star f)\neq 0
\end{gather*} 
for all $\eta\in\Ha$.
\qed

\begin{theo}\label{crossproduct1}
There exists an map from the set of all closed subgroups $\Hr$ of $\Gr$ and $\Gr$-invariant von Neumann subalgebras $\Na=\Na''$ of $\Ma\rtimes_\alpha\Gr$, $(\Ma\rtimes_\alpha\Gr)^\Gr\subset\Na$.
\end{theo}

\proof
Due to Lemma \ref{general1} the map $\tilde{\pi}_M$ is injective and thanks to Theorem \ref{generalTakesaki} we may apply Theorem \ref{spatial2}.
\qed

\begin{cor}
There exists an injective map from the set of all closed subgroups $\Hr$ of $\Gr$ and $\Gr$-invariant von Neumann subalgebras $\Aa=\Aa''$ of $\Ma$, $\Ma^\Gr\subset\Aa$.
\end{cor}

\proof
Set $\Aa:=\Na\cap\Ma$.
\qed

\begin{theo}
Let $\Ma\rtimes_\alpha\Gr$ be a factor and $\Ge$ the enveloping group, then there exists a bijective map between the set of full subgroups of $\Ge$ and the set of subalgebras $\Na$ of $\big(\Ma\rtimes_\alpha\Gr\big)'\cap\Sa$ with $\Na\rc{}'\cap\Sa=\Na$ and hence to the set of subalgebras of $\Ma\rtimes_\alpha\Gr$ containing $\big(\Ma\rtimes_\alpha\Gr\big)^\Gr$.
\end{theo}

\proof
Apply Theorem \ref{spatialfactor}.
\qed

\begin{examp}
If the action of a countably infinite discrete group $Gr$ on an abelian von Neumann algebra $\Ma$ is free and ergodic, then the corresponding crossed product is a factor \cite[V.7.8]{Takesaki1:1970}.
\end{examp}

\begin{defi}
The restriction of the enveloping group to $\Ma$ is defined as $\Ge_0:=\{g\in\Ge|\;\alpha^g(\Ma)=\Ma\}$.
\end{defi}

\begin{theo}
If $\alpha:\Gr\longrightarrow\Aut(\Ma)$ contains all irreducible representations, then the map of the set of all closed subgroups of $\Ge_0$ to the set of all subalgebras $\Aa$ with $\Ma^\Gr\subset\Aa\subset\Ma$ is bijective.
\end{theo}

\proof
Apply Theorem \ref{spatial-one-to-one}.
\qed

\subsection{Minimal Action}

The aim of this section is to localise the case investigated in \cite{Izumi:1998} in the framework of this thesis. There, the authors analyse the minimal or outer action of a compact or discrete group $\Gr$, respectively, on a factor $\Ma$ with separable predual. \\
Minimal action, i.e. $\big(\Ma^\Gr\big)'\cap\Ma=\C\1$, is a crucial assumption in their investigation since it implies central ergodicity of $\alpha$ which is equivalent to $\Na\rtimes_\alpha\R$ being a factor for all fixed-point von Neumann subalgebras $\Na$ of $\Ma$, confer Remark 4.5 of \cite{Izumi:1998}. Moreover, their approach relies heavily on the existence of an unique normal conditional of $\Ma$ onto all intermediate subfactor which is ensured in case of minimal action by Theorem \ref{condexpminimalaction}. Consequently, they obtains a one-to-one correspondence between closed subgroups of $\Gr$ and subfactors of $\Ma$. \\
The ansatz of our analysis is more general because we do only demand of the crossed product to be a factor or equivalently central ergodicity of $\alpha$. Minimal action restricts our formalism to only spatial automorphisms as in \cite{Izumi:1998}.
\begin{lem}
The set $\big(\Ma^\Gr\big)'\cap\Sa$ consists in the case of minimal action only of 'spatial' operators.
\end{lem}

\proof 
\begin{align*}
\big(\Ma^\Gr\big)'\cap\Sa&=\big(\Ma^\Gr\big)'\cap\big[\Ma\cup\Sa\backslash\Ma\big]\\
&=\C\1\cup\Big[\big(\Ma^\Gr\big)'\cap\Sa\backslash\Ma\Big]\\
&=\Big\{S\in\LH|\;SA=AS\quad\forall A\in\Ma^\Gr,\,S\notin\Ma\Big\}
\end{align*}
\qed

\section{Non-Commutative Martingales}

As announced, this section contains first steps towards a full analysis of the interplay between Galois correspondence for von Neumann algebras and non-commutative probability. We want to conclude this chapter with the application of our results to non-commutative probability, namely we begin by identifying non-commutative martingales in our apparatus and then show for them a non-abelian version of Theorem \ref{convergencemartingale}.

\begin{prop}
Let $\Gr$ be a compact group then for faithful, normal and semifinite state $\varphi:\Ma\longrightarrow[0,1]$ on the von Neumann algebra $\Ma$ 
\begin{gather*}
\phi:=\int_\Gr(\alpha^g)^t(\varphi)d\mu(g)
\end{gather*}
is also a faithful, normal and semifinite state on $\Ma$. Moreover, $\phi$ is $\Gr$-invariant.
\end{prop}

\proof
Let $\varphi$ be faithful, see Definition \ref{definitionsforstates}, then one obtains for all positive element $A^*A$ of $\Ma$:
\begin{align*}
\phi(A^*A)&=\int_\Gr(\alpha^g)^t(\varphi)(A^*A)d\mu(g)\\
&=\int_\Gr\varphi\big(\alpha^g(A^*A)\big)d\mu(g)\\
&=\int_\Gr\varphi\big(\alpha^g(A)^*\alpha(A)\big)d\mu(g)\\
&>0.
\end{align*}
Normality of the state $\varphi$, confer Definition \ref{definitionnormal}, transfers also to $\phi$:
\begin{align*}
\text{l.u.b.}_n\phi(A_n)&=\text{l.u.b.}_n\int_\Gr(\alpha^g)^t(\varphi)(A_n)d\mu(g)\\
&=\int_\Gr\text{l.u.b.}_n\,\varphi\big(\alpha^g(A_n)\big)d\mu(g)\\
&=\int_\Gr\varphi\big(\text{l.u.b.}_n\,\alpha^g(A_n)\big)d\mu(g)\\
&=\int_\Gr(\alpha^g)^t(\varphi)(\text{l.u.b.}_n\,A_n)d\mu(g)\\
&=\phi(\text{l.u.b.}_nA_n).
\end{align*}
Finally, if we suppose $\varphi$ to be semifinite, then by Definition \ref{definitionsemifinite} the state $\phi$ is so, too, because $\varphi(A)<\infty$ implies $\phi(A)<\infty$ and the set
\begin{gather*}
\Ma_+^\phi:=\big\{A\in\Ma_+|\;\phi(A)<\infty\big\}
\end{gather*}
\\ \\has to be dense in $\Ma$.
\qed

\begin{theo}
Let $(\Ma,\varphi)$ be a non-commutative probability space and $\Hr$ a closed subgroup of the compact group $\Gr$ then the mapping 
\begin{align*}
\E_\Hr:\Ma&\longrightarrow\Ma^\Hr\\
A&\mapsto\E_\Hr(A):=\int_\Hr h(A)d\mu(h)
\end{align*}
is a conditional expectation of $\Ma$ onto $\Ma^\Hr$ with respect to $\phi$.
\end{theo}

\proof
Obviously, $\E_\Hr$ is a mapping onto $\Ma^\Hr$ and satisfies the conditions of Definition \ref{noncommcondexp}, namely:
\begin{itemize}
\item[(i)] $\|\E_\Hr(A)\|=\|\int_\Hr h(A)d\mu(h)\|\leq\int_\Hr\|h(A)\|d\mu(h)=\|A\|$;
\item[(ii)] $\E_\Hr(A)=\int_\Hr h(A)d\mu(h)=A$ since  $h(A)=A$;
\item[(iii)] The third property is a consequence of translation invariance of the Haar measure:
\begin{align*}
\phi\circ\E_\Hr(A)&=\int_\Gr g^t\left(\varphi\left[\int_\Hr h(A)d\mu(h)\right]\right)d\mu(g)\\
&=\int_\Gr\int_\Hr g^t\big(\varphi\big[h(A)\big]\big)d\mu(h)d\mu(g)\\
&=\int_\Gr\int_\Hr\varphi\big[g\circ h(A)\big]d\mu(h)d\mu(g)\\
&=\int_\Hr\int_\Gr\varphi\big[g\circ h(A)\big]d\mu(g)d\mu(h)\\
&=\int_\Hr\int_\Gr\varphi\big[g(A)\big]d\mu(g)d\mu(h)\\
&=\int_\Gr\varphi\big[g(A)\big]d\mu(g)\\
&=\phi(A).
\end{align*}
\end{itemize}
\qed




\begin{theo}
Let $(\Ma,\varphi)$ be a non-commutative probability space, $\Gr$ a compact group and $(\Hr_t)_{t\in\TT}$ a sequence of closed subgroups of $\Gr$ with $\Hr_s\supset\Hr_t$ for $s<t$ and $\Hr^0\equiv\Gr$. Then the family $(X_t)_{t\in\TT}$ of random variables on $(\Ma,\varphi)$ adapted to the filtration $\big(\Ma^{\Hr_t}\big)_{t\in\TT}$ constitutes a non-commutative martingale.
\end{theo}

\proof
According to the following diagram for $t>s$
\begin{diagram}
\Aa  & \rTo^{\hspace{1,5cm} X}              &&& \Ma\\
     &\rdTo(4,2)^{X_{t}}\rdTo(4,4)^{X_{s}}  &&& \dTo_{\E_{\Hr_{t}}}  \\
     &                                      &&& \Ma^{\Hr_{t}}\\
     &                                      &&& \dTo_{\E_{\Hr_{s}}}   \\ 
     &                                      &&& \Ma^{\Hr_{s}}\\
\end{diagram}
we obtain:
\begin{align*} 
\E\big(X_t|\Ma^{\Hr_s}\big)&\equiv\E_{\Hr_s}(X_t)\\
&=\int_{\Hr_s} h(X_t)d\mu(h)\\
&=\int_{\Hr_s}\int_{\Hr_t} h\circ h'(X)d\mu(h')d\mu(h)\\
&=\int_{\Hr_t}\int_{\Hr_s} h\circ h'(X)d\mu(h)d\mu(h')\\
&=\int_{\Hr_t}\int_{\Hr_s} h(X)d\mu(h)d\mu(h')\\
&=\int_{\Hr_s} h(X)d\mu(h)\\
&=\E_{\Hr_s}(X)=X_s.
\end{align*}
\qed
We want to give now a non-abelian version for the convergence theorem, Theorem \ref{convergencemartingale}.

\begin{theo}
Let $(X_t)_{t\in\TT}$ be a non-commutative martingale satisfying
\begin{gather*}
E\big(|X_1^*X_1|\big)\leq E\big(|X_2^*X_2|\big)\leq...
\end{gather*}
and
\begin{gather*}
\emph{weak-}\lim_{i\rightarrow\infty}E\big(|X_i^*X_i|\big)=:c<\infty,
\end{gather*}
then $X_\infty:=\emph{strong-}\lim_{i\rightarrow\infty}X_i$ exists and the family $(X_1,X_2,...,X_\infty)$ constitutes a martingale.
\end{theo}

\proof
Since the sequence of expectations is increasing and bounded, its weak limit and, therefore, the weak limit of the sequence $\big(X_t^*X_t)_{t\in\TT}$ are given. But the existence of the latter one is equivalent to the existence of the strong limit of the sequence $(X_t)_{t\in\TT}$, confer for more details \cite[Lemma II.2.5]{Takesaki:1970}.\\
Furthermore, the martingale property of $(X_t)_{t\in\TT}$
\begin{align*}
\E\big(X_t|\Ma^{\Hr_s}\big)=X_s,\qquad s<t,s,t\in\TT,
\end{align*}
leads, due to continuity of the conditional expectation, to
\begin{align*}
&&\text{strong-}\lim_{i\rightarrow\infty}\E\big(X_t|\Ma^{\Hr_s}\big)&=X_s&\\
&\Longleftrightarrow&\hspace*{-2cm}\E\big(\text{strong-}\lim_{i\rightarrow\infty}X_t|\Ma^{\Hr_s}\big)&=X_s&\\
&\Longleftrightarrow&\E\big(X_\infty|\Ma^{\Hr_s}\big)&=X_s\,.&
\end{align*}
Thus, the martingale property for the family $(X_1,X_2,...,X_\infty)$ is also ensured.
\qed


\setcounter{secnumdepth}{2}


\chapter{Summary and Outlook}

\begin{flushright}
 \emph{There is something fascinating about science.\\ One gets such wholesale returns of conjecture\\ out of such a trifling investment of fact.}\\
\vspace{0,5cm}
Mark Twain
\end{flushright}
\vspace{0,5cm}
This thesis has been concerned with the interaction of Galois theory, operator algebras and non-commutative probability theory. \\
It has been shown that classical Galois theory is a powerful tool in the analysis not only of fields but may play a decisive r\^ole in operator algebras as well. We have given a new approach for its application to a von Neumann algebra $\Ma$ for which we analysed one-to-one correspondences between subgroups of a compact group $\Gr$ and von Neumann subalgebras of $\Ma$ without demanding neither further restrictions on the nature of the action of $\Gr$ nor additional properties for the subalgebras, contrary to the existing literature. \\
Based on the ansatz of Coja-Oghlan and Michali\v cek \cite{Coja:2005}, we have continued the introduction of non-commutative probability theory, namely we proposed non-abelian analogues of stochastic notions and objects such as conditional expectation, stochastic independence, stochastic processes and martingales.\\ \\
Future investigations on the subject of this thesis may consider two different directions, namely the generalisation of its mathematical foundation and the quest of applications to other disciplines.\\
The generalisation and implementation of this approach to locally compact (abelian) groups is under investigation, but it is not clear if one may benefit  from more general formulations of Peter-Weyl theorem. Abelian locally compact groups have been analysed by Connes and Takesaki \cite{Connes:1978}. Since we have emphasised here the appearance of non-commutative martingales in our framework of Galois correspondence, non-commutative probability in general and stochastic processes in particular may be dealt with in a more systematically manner.\\
As mentioned in the introduction, both subjects Galois correspondence for von Neumann algebras as well as non-abelian probability theory may be deployed on several fields. One of them is obviously mathematical physics, especially the algebraic formulation of quantum field theory where von Neumann algebras play a central part and any information on the nature of their substructure is of decisive significance. Non-commutative stochastic processes will obtain increasing influence, because they are tailored perfectly for the description of quantum effects in classical processes such as the Brownian motion. Another not so apparent field of implementation is financial mathematics where the first steps are already made, to wit formulating a non-commutative analogue of Black-Scholes equations \cite{Forgy:2002}, \cite{Chen:2001}.












\cleardoublepage
\addcontentsline{toc}{chapter}{\protect\numberline{Notation}}
\chapter*{Notation}

\begin{tabular}{ll}
$\mathcal{C}^{l}(\Omega)$    & vector space of $l$-times continuously differentiable \\
                             & functions on $\Omega$\\
$\Cnn(\Omega)$               & vector space of continuously differentiable \\
                             & functions on $\Omega$ vanishing at infinity\\
$\mathcal{C}(\Omega)$        & $=\mathcal{C}^{0}(\Omega)$ continuous functions on $\Omega$\\
$\mathcal{C}^{\infty}(\Omega)$ & $=\cap\{\mathcal{C}^{l}(\Omega)|\;l\in\N_{0}\}$\\
$\mathcal{B}(\R^n)$          & $=\big\{f\in\Cg(\R^n)|\;\forall\alpha\in\N_0^n:\;\sup\{|D^\alpha f(x)|\,|x\in\R^n\}<\infty\big\}$\\
$\mathcal{E}(\Omega)$        & $=\mathcal{C}^{\infty}(\Omega)$\\
$\mathcal{D}(\Omega)$        & $=\mathcal{C}_{0}^{\infty}(\Omega)$\\
$\mathcal{E}'(\Omega)$       & dual space of $\mathcal{E}(\Omega)$\\
$\mathcal{D}'(\Omega)$       & dual space of $\mathcal{D}(\Omega)$\\
$\mathcal{S}'(\R^{n})$       & dual space of $\mathcal{D}(\R^{n})$\\
\\
$\X$                         & normed space\\
$\X^*$                       & dual space of $\X$\\
$\langle u,\phi\rangle$ & Application of the distribution $u$ on $\phi\in\mathcal{X}\quad\big(\mathcal{X}\in\{\mathcal{E}(\Omega),\mathcal{D}(\Omega),\mathcal{S}(\R^{n})\}\big)$\\
$(u,\phi)$ & $=\langle u,\overline{\phi}\rangle$\\
$\widehat{u}$ & Fourier transform of $u$\\
\\
\\
$\Gr$                        & (locally) compact group\\
$\Ge$                        & $:=\prod_{\sigma\in\Sigma}\U(\Ha_{\sigma,\Gr}\big)$, enveloping of the group $\Gr$\\
$\Hr$                        & subgroup of $\Gr$\\
$\Ha$                        & Hilbert sapce\\
$\Ka$                        & subspace of $\Ha$\\
$\Ka^\perp$                  & orthogonal complement of $\Ka$\\
$\LH$                        & set of bounded operators on $\H$\\
$\alpha$                     & $:\Gr\longrightarrow\LH$, representation of $\Gr$ on $\LH$\\
$K_u$                        & Weyl operator\\
$D_{ij}(g)$                 & matrix elements with respect to $\alpha$\\
$\Da(\Gr)$                   & coefficient algebra\\
$\Sigma\equiv\Sigma(\Gr)$    & equivalence class of finite-dimensional, irreducible and\\
                             & unitary representations\\
$\chi(g)$                    & characters with respect to $\alpha$\\
\end{tabular}

\newpage

\begin{tabular}{ll}
$(\Omega,\SA,P)$             & classical probability space\\
$(\Ma,\varphi)$              & non-commutative probability space\\
$\SA,\SB,\SF$                & $\sigma$-algebras\\
$\mathbb{B}$                 & Borel-$\sigma$-algebra\\
$X,Y$                        & random variables\\
$(\SF_t)_{t\in I}$           & classical filtration\\
$(\Ma_t)_{t\in I}$           & non-commutative filtration\\
$\E$                         & conditional expectation\\
$\Ua$                        & $U^{*}$-algebra\\
$\Aa,\Ba$                    & $C^{*}$-algebras\\
$\Ma,\Na$                    & von Neumann algebras\\
$(\Aa,\mathbf{G},\alpha)$    & $C^*$-dynamical system\\
$(\Ma,\mathbf{G},\alpha)$    & $W^*$-dynamical system\\
$\Za$                        & center\\
$\Aa^\sim$                      & universal enveloping of $\Aa$\\
$\Aa^*$                      & dual of $\Aa$\\
$\Aa_*$                      & predual of $\Aa$\\
$\Sa$                        & state space\\
$\Fa$                        & factor state space\\
$\Pa$                        & pure state space\\
$\Aa'$                       & commutant of $\Aa$\\
$\U(\Aa)$                    & unitary elements of $\Aa$\\
$[A,B]$                      & $:=AB-BA$, commutator of $A,B\in\LH$\\
$\Aa\rc$                     & $:=\Aa'\cap\Aa$, relative commutant of $A$\\
$\Ma^\Gr$                    & $\Gr$-invariant subalgebra\\
$\xi,\eta$                        & $\in\Ha$\\
$(\Ha_\o,\pi_\o,\xi_\o)$     & cyclic representation with respect to $\o$\\
$\o,\o_1,\o_2$               & states on $\Ma$\\
\\ \\
$S$, $J$, $\Delta$           & Tomita oerator, modular conjugation, modular operator\\
$\sigma^t(A)$                & $:=\Delta^{it}A\Delta^{-it}$ modular group of automorphisms on $\Ma$\\
$\sigma_\o^t$                & modular automorphism group on $\Ma$ with respect to the state $\o$\\     
$\Gamma_t$                   & Connes' cocycle\\
$\delta$                     & derivation, infinitesimal generator\\
$\Sigma\equiv\Sigma(\Gr)$    & equivalence class of finite-dimensional, irreducible and\\
                             & unitary representations
\end{tabular}

\newpage      

\begin{tabular}{ll}
$\autMa$                     & group of automorphisms of $\Ma$\\
$\innMa$                     & group of inner automorphisms of $\Ma$\\
$\spaMa$                     & group of spatial automorphisms of $\Ma$\\
$\outMa$                     & group of outer automorphisms of $\Ma$\\
$\alpha$                     & $:\Gr\longrightarrow\autMa$\\
$\beta$                      & $:\Gr\longrightarrow\L\big(L^2(\Gr)\big)$\\
$\tilde{\beta}$              & $:M(\Gr)\longrightarrow\L\big(L^2(\Gr)\big)$\\
$\lambda$                    & $:\Gr\longrightarrow M(\Gr)$\\
$\iota_1$                    & $:\L^1(\Gr)\longrightarrow M(\Gr)$\\
$\iota_2$                    & $:\L^2(\Gr)\longrightarrow L^1(\Gr)$\\
$\pi_\alpha$                 & $:\Gr\longrightarrow\L\big(L^2(\Gr,\Ma)\big)$\\
$\tilde{\pi}_{L^1}$          & $:L^1(\Gr)\longrightarrow\left(\Ma^\Gr\right)\rc$\\
$\tilde{\pi}_{L^2}$          & $:L^2(\Gr)\longrightarrow\Ha_{\Sigma,\Gr}^2$\\
$\tilde{\pi}_M$              & $:M(\Gr)\longrightarrow\left(\Ma^\Gr\right)\rc$\\
$\Ma\rtimes_\alpha\Gr$       & crossed product of $\Ma$ by $\alpha$\\
$\Sigma^0$                   & $:=\ker(\tilde{\pi}_M\circ\iota_1\circ\iota_2)$\\
$\Sigma'$                    & $:=\Sigma\backslash\Sigma^0$\\
$\La_{\Sigma,\Gr}^p$         & $:=\left\{(T_\sigma)_{\sigma\in\Sigma}\Big|\;\Big(\sum_\sigma\|T_\sigma\|^p_{\text{uniform}}\Big)^{1/p}\hspace{-3mm}<\infty,T_\sigma\in\L\big(\Ha_{\sigma,\Gr}\big)\right\}$\\
$\Lainf_{\Sigma,\Gr}$        & $:=\Big\{(T_\sigma)_{\sigma\in\Sigma}\big|\;\sup_\sigma\|T_\sigma\|_{\text{uniform}}<\infty,\;T_\sigma\in\L\big(\Ha_{\sigma,\Gr}\big)\Big\}$\\
$\La_{\Sigma,\Gr}^0$         & $:=\left\{(T_\sigma)_{\sigma\in\Sigma}\big|\;\lim_\sigma\|T_\sigma\|_{\text{uniform}}=0,\;T_\sigma\in\L\big(\Ha_{\sigma,\Gr}\big)\right\}$\\
$\Ha_{\Sigma,\Gr}^2$         & $:=\left\{\big((\xi_\sigma^1,...,\xi_\sigma^{n_\sigma})\big)_{\sigma\in\Sigma}\big|\;\sum_{k=1}^{n_\sigma}\|\xi_\sigma^k\|^2<\infty\right\}$
\end{tabular}


\cleardoublepage

\addcontentsline{toc}{chapter}{\protect\numberline{Bibliography}}
\bibliographystyle{cmp}
\bibliography{dissmgesamt}
\nocite{*}



\clearpage

\addcontentsline{toc}{chapter}{\protect\numberline{Acknowledgements}}
\chapter*{Acknowledgements}

Sincere thanks are given to those who have guided, supported and accompanied me during my studies on this dissertation.\\ \\
This thesis would not have been possible without Prof Dr J. Michali\v cek who not only has been my supervisor but also has guided me from the very first semester. His encouragement after setbacks, his patience during hours-long discussions and his freely sharing of ideas were of immense help. \\
Financial security and complete freedom of research are indispensable ingredients for a fulfilled academic life and, therefore, I owe a debt of gratitude to Prof Dr W. Krumbholz for supporting me during the last eight years. Herr Krumbholz, 1. d4\hspace{0,5mm}!\\
I am also deeply grateful to Prof Dr H. Hebbel and Prof Dr U. T\"ushaus for their immediate readiness to serve as examiners.\\
Moreover, I would like to thank all friends and colleagues for the warm and friendly working atmosphere, the lively discussions at our Teerunde and, of course, the blitz chess matches.  I shall always treasure this memory. \\ \\
Special thanks go to my family, in particular to my parents who accepted and mastered for their children countless burdens and obstacles during the last three difficult decades. Tashakur.

\thispagestyle{empty}

\end{document}